\documentclass [12pt, a4paper]{article}
\usepackage[utf8]{inputenc}
\usepackage[T1]{fontenc}
\usepackage{amsmath}
\usepackage{amsfonts}
\usepackage{amsthm}
\usepackage{color}

\newlength{\oldparindent}
\newcommand{\bpf}{\begin{preuve}}
\newcommand{\epf}{ \end{preuve} \medskip}

\newcommand{\benum}{\begin{enumerate}}
\newcommand{\eenum}{\end{enumerate}}

\newcommand{\bitem}{\begin{itemize}}
\newcommand{\eitem}{\end{itemize}}

\newcommand{\brmq}{\begin{rmq}}
\newcommand{\ermq}{\end{rmq}}

\newcommand{\brmqs}{\begin{rmqs}}
\newcommand{\ermqs}{\end{rmqs}}

\newcommand{\bapp}{\begin{application}}
\newcommand{\eapp}{\end{application}}

\newcommand{\bapps}{\begin{applications}}
\newcommand{\eapps}{\end{applications}}

\newcommand{\bdefi}{\begin{definition}}
\newcommand{\edefi}{\end{definition}}

\newcommand{\beq}{\begin{equation}}
\newcommand{\eeq}{\end{equation}}

\def\bpm{\begin{pmatrix}}
\def\epm{\end{pmatrix}}

\newcommand{\bcas}{\begin{cases}}
\newcommand{\ecas}{\end{cases}}

\newcommand{\bex}{\begin{exemp}}
\newcommand{\eex}{\end{exemp}}

\newcommand{\bexs}{\begin{exemps}}
\newcommand{\eexs}{\end{exemps}}

\newcommand{\bprop}{\begin{proposition}}
\newcommand{\eprop}{\end{proposition}}

\newcommand{\bthm}{\begin{theoreme}}
\newcommand{\ethm}{\end{theoreme}}

\newcommand{\bcor}{\begin{corollaire}}
\newcommand{\ecor}{\end{corollaire}}

\newcommand{\blem}{\begin{lemme}}
\newcommand{\elem}{\end{lemme}}

\newcommand{\beqna}{\begin{eqnarray}}
\newcommand{\eeqna}{\end{eqnarray}}

\newcommand{\beqnas}{\begin{eqnarray*}}
\newcommand{\eeqnas}{\end{eqnarray*}}

\definecolor{green}{rgb}{0,.7,.2}
\definecolor{orange}{rgb}{0.9,.5,0}

\def\Ric{ {\rm{Ric}}} %parentheses have been corrected
\def\Hess{ {\rm{Hess}}}
\def\Ent {{\rm{Ent}}}   %parentheses have been corrected

\def\Id{{\rm{Id}}} %parentheses have been corrected

\def\cC{{\mathcal C}}

\def\cP{{\mathcal P }}

\newcommand{\bbN}{{\mathbb {N}}}% entiers
\newcommand{\bbR}{{\mathbb {R}}}%reels
\def\bf\Sigma{{\mathbf{\Sigma}}}

\newtheorem{theoreme}{Theorem}[section]
\newtheorem{lemme}[theoreme]{Lemma}
\newtheorem{definition}[theoreme]{Definition}
\newtheorem{proposition}[theoreme]{Proposition}
\newtheorem{corollaire}[theoreme]{Corollary}

\newenvironment{exemp}{\noindent{\bf Example. --- }}{\par}
\newenvironment{exemps}{\noindent{\bf Examples}\benum}{\eenum\par}
\newtheorem{rmq}[theoreme]{Remark}
\newtheorem{rmqs}[theoreme]{Remarks}

\newenvironment{preuve}{\noindent{\it Proof. --- }}
{\hfill\rule{1.3mm}{2mm}\par}
\newenvironment{application}{\noindent{\bf Application. --- }}{\par}
\newenvironment{applications}{\noindent{\bf Applications. ---
}\benum}{\eenum\par}

\theoremstyle{definition}

\date{\today}

\begin{document}

\title{ \textbf{\large $W$-entropy formulas and Langevin deformation of flows on Wasserstein space over Riemannian manifolds}}
\author{\ \ Songzi Li\thanks{Research of S. Li has been supported by NSFC No. 11901569.},\ \ \ \ \  Xiang-Dong Li
\thanks{ Research of X.-D. Li has been supported by National Key RD Program of
China (No. 2020YF0712700), NSFC No. 11771430 and No. 11688101, and Key Laboratory RCSDS, CAS,
No. 2008DP173182.} }
  
\maketitle

\begin{abstract} We introduce Perelman's  $W$-entropy and prove the $W$-entropy formula along the geodesic flow on the $L^2$-Wasserstein space  over compact Riemannian manifolds equipped with Otto's Riemannian metric, which allows us to recapture a previous result due to Lott and Villani 
on the displacement convexity of $s{\rm Ent}+ns\log s$ on $P^\infty_2(M)$ over Riemannian manifolds with non-negative Ricci curvature.  To better understand the similarity  between 
the $W$-entropy formula for the geodesic flow on the Wasserstein space and the $W$-entropy formula for the heat equation of the Witten Laplacian on the underlying manifolds,
we introduce the Langevin deformation of  flows on the Wasserstein space over Riemannian manifold, which interpolates  the gradient flow and the geodesic flow on the Wasserstein space over Riemannian manifolds, and can be regarded as 
the potential flow of the compressible Euler equation with damping on manifolds. We prove the existence, uniqueness and regularity of the Langevin deformation  on  the Wasserstein space over the Euclidean space and compact Riemannian manifolds,  and prove the convergence of the Langevin deformation for $c\rightarrow 0$ and $c\rightarrow \infty$ respectively.
We prove an analogue of the Perelman type 
$W$-entropy formula along the Langevin deformation on the Wasserstein space on  Riemannian manifolds. A rigidity theorem is proved for the $W$-entropy for the geodesic flow, and a rigidity model is also provided for the Langevin deformation on the Wasserstein space over complete Riemannian manifolds with the $CD(0, m)$-condition.
\end{abstract}


\medskip


%\tableofcontents

\section{Introduction}

In his seminal paper~\cite{P1}, Perelman pointed out that the Ricci flow can be realized as the 
gradient flow of the $\mathcal{F}$-functional and proved the monotonicity of its $\mathcal{W}$-entropy, which plays   
a key role in his proof of the Poincar\'e conjecture. Meanwhile in the field of optimal transport, the infinite dimensional Riemannian geometry and the theory of the gradient flow on the Wasserstein space 
have been developed by Otto, Lott, McCann, Villani and Sturm \cite{Ot, OtV, LoV, Lo2, V1, V2, EKS, St1, St2, St3} among others, which leads to the reveal of the deep relation between the gradient flow and the entropy formula. Moreover, the convexity of the Boltzmann-Shannon entropy or the R\'enyi entropy along geodesics on the Wasserstein space has been a key tool in Lott-Villani  \cite{LoV, Lo2, V1, V2} and Sturm \cite{St1, St2, St3}  to develop a synthesis of comparison geometry on metric measure spaces with the extended notion of the curvature-dimension $CD(K, m)$-condition. See \cite{AGS} for the further development on the theory of gradient flows on metric measure spaces.

We start with a brief description of the key concepts of the topic and then present our main results. Let $(M, g)$ be a complete Riemannian manifold equipped with a weighted volume measure $d\mu=e^{-f}d\nu$, where $f\in C^2(M)$ and $d\nu$ denotes the volume measure on $(M, g)$. The Witten Laplacian, which is self-adjoint and non-positively definite on $L^2(M, \mu)$,  is defined as follows
\begin{eqnarray*}
L:=\Delta-\nabla f\cdot\nabla.
\end{eqnarray*} 
Let $m>n$ be a constant. Following Bakry and Emery \cite{BE}, we introduce the $m$-dimensional Bakry-Emery Ricci curvature associated with the Witten Laplacian $L$,
\begin{eqnarray*}
\Ric_{m, n}(L):=\Ric+\nabla^2 f-{\nabla f\otimes \nabla f\over m-n}.
\end{eqnarray*} 
We also make the convention that when $m=n$, $\Ric_{n, n}(L)=\Ric$  and is only defined when  $f$ is identically a constant.  
By \cite{Bes, Lot}, when $m>n$ is an integer,  $\Ric_{m, n}(L)$ is the horizontal projection of the Ricci curvature on the product manifold $\widetilde{M}=M\times N$ with the warped product metric $\widetilde g=g\otimes e^{-{2f\over m-n}}g_{N}$, where $(N, g_N)$  is any  $(m-n)$-dimensional complete Riemannian  manifold. 
Following \cite{BE}, we say that the weighted Riemannian manifold 
$(M, g, \mu)$ satisfies the $CD(K, m)$-condition for some constant $K\in \mathbb{R}$ and $m\in [n, \infty]$ if and only if 
$$\Ric_{m, n}(L)\geq K.$$ 

The Boltzmann-Shannon entropy of the probability measure $\rho d\mu$ with respect to the reference measure $\mu$ is defined by
\begin{eqnarray*}
{\rm Ent}(\rho):= \int_M \rho \log \rho d\mu.
\end{eqnarray*}

Let $P_2(M, \mu)$ (resp. $P_2^\infty(M, \mu)$) be the Wasserstein 
space (reps. the smooth Wasserstein space) of all probability measures $\rho(x)d\mu(x)$ with density function (resp. with smooth density function) $\rho$ on $M$ such that $\int_M d^2(o, x)\rho(x)d\mu(x)<\infty$,  where $d(o, \cdot)$ denotes the distance function from a fixed point $o\in M$. 

By Otto \cite{Ot}, the tangent space $T_{\rho d\mu}P_2^\infty(M, \mu)$ is identified as follows 
\begin{eqnarray*}
T_{\rho d\mu}P_2^\infty(M, \mu)=\{s=\nabla_\mu^*(\rho \nabla\phi): \phi\in C^\infty(M), \ \ \int_M |\nabla\phi|^2\rho d\mu<\infty\},
\end{eqnarray*}
where
$\nabla_\mu^*$ denotes the $L^2$-adjoint of the Riemannian gradient $\nabla$ with respect to the weighted volume measure $d\mu$  on $(M, g)$.
For $s_i=\nabla_\mu^*(\rho\nabla\phi_i)\in T_{\rho d\mu} P_2^\infty(M, \mu)$, $i=1, 2$, Otto \cite{Ot} introduced the following infinite dimensional Riemannian metric on $P_2^\infty(M, \mu)$ 
\begin{eqnarray*}
\langle\langle s_1, s_2\rangle\rangle:=\int_M \langle \nabla \phi_1, \nabla\phi_2 \rangle \rho d\mu,
\end{eqnarray*}
provided that 
\begin{eqnarray*}
\|s_i\|^2:=\int_M |\nabla\phi_i|^2\rho d\mu<\infty, \ \ \ i=1, 2.
\end{eqnarray*}
Let $T_{\rho d\mu}P_2(M, \mu)$ be the completion of $T_{\rho d\mu}P_2^\infty (M, \mu)$ with Otto's Riemannian metric. Then $P_2(M, \mu)$ is a formal infinite dimensional Riemannian manifold.

By  Benamou and Brenier \cite{BB}, for any given $\mu_i=\rho_i d\mu\in P_2(M, \mu)$, $i=0, 1$, the $L^2$-Wasserstein distance between $\mu_0$ and $\mu_1$ coincides with the geodesic distance between $\mu_0$ and $\mu_1$ on $P_2(M, \mu)$ equipped with Otto's infinite dimensional Riemannian metric, i.e.,  
\begin{eqnarray*}
W_2^2(\mu_0, \mu_1)=\inf\limits\left\{{1\over 2}\int_0^1 |\nabla\phi(x, t)|^2\rho(x, t)d\mu(x): \partial_t \rho=\nabla_\mu^*(\rho \nabla\phi), \ \rho(0)=\rho_0, \ \rho(1)=\rho_1\right\}.
\end{eqnarray*}
Given $\mu_0=\rho(\cdot, 0)\mu, \mu_1=\rho(\cdot, 1)\mu\in P_2^\infty(M, \mu)$, it is known that there is a unique minimizing Wasserstein geodesic $\{\mu(t), t\in [0, 1]\}$  of the form $\mu(t) =(F_t)_*\mu_0$  joining $\mu_0$ and $\mu_1$ in $P_2(M, \mu)$, where $F_t \in {\rm Diff}(M)$  
is given by $F_t(x) = \exp_x(-t \nabla \phi(\cdot, 0))$  for an appropriate Lipschitz function $\phi(\cdot, t)$ (see \cite{Mc}). 
 If the Wasserstein 
geodesic in $P_2(M, \mu)$  belongs entirely to $P_2^\infty(M, \mu)$,  then the geodesic flow $(\rho, \phi)\in TP_2^\infty(M, \mu)$ satisfies the transport equation  and the Hamilton-Jacobi equation 
\begin{eqnarray}
{\partial_t} \rho-\nabla_\mu^*(\rho \nabla \phi)&=&0,\label{TA}\\
{\partial_t}\phi+{1\over 2}|\nabla \phi|^2&=&0, \label{HJ}
\end{eqnarray}
with the boundary condition $\rho(0)=\rho_0$ and $\rho(1)=\rho_1$.   When  $\rho_0 \in C^\infty(M, \mathbb{R}^+)$ and $\phi_0\in C^\infty(M)$,  defining $\phi(\cdot, t)\in C^\infty(M)$ by the Hopf-Lax solution
 \begin{eqnarray}
\phi(x, t)=\inf\limits_{y\in M}\left(\phi_0(y)+{d^2(x, y)\over 2t}\right),\label{HLS}
 \end{eqnarray} 
 and solving the transport equation $(\ref{TA})$ by the characteristic method,  it is known that $(\rho, \phi)$ satisfies $(\ref{TA})$ and $(\ref{HJ})$ with $\rho(0)=\rho_0$ and $\phi(0)=\phi_0$. See  \cite{V1}  Sect. 5.4.7. See also \cite{Lo1, Lo2}.
   In view of this, the transport equation $(\ref{TA})$ 
 and the Hamilton-Jacobi equation $(\ref{HJ})$ describe the geodesic flow on the tangent bundle $TP_2^\infty(M, \mu)$ over the Wasserstein space $P_2(M, \mu)$. 
% Note that the Hamilton-Jacobi equation $(\ref{HJ})$ is also called the eikonnal equation in geometric optics, see e.g. \cite{Car}. 

Our first result of this paper is the following $W$-entropy formula for the geodesic flow on the Wasserstein space $P_2(M, \mu)$. 
\bthm\label{MT2}
Let $(M, g)$  be a compact Riemannian manifold or a complete Riemannian manifold with suitable bounded geometry condition,  $f\in C^2(M)$, $d\mu=e^{-f}dv$.
Let $(\rho, \phi): M\times [0, T]\rightarrow \mathbb{R}^+\times \mathbb{R}$  be a geodesic in $P_2^\infty(M, \mu)$. 
For any $m\geq n$, define the $H_m$-entropy and $W_m$-entropy for the geodesic flow $(\rho, \phi)$ on $TP^\infty_2(M, \mu)$ as follows
\begin{eqnarray*}
H_m(\rho, t)={\rm Ent}(\rho(t))+{m\over 2}\left(1+\log(4\pi t^2)\right),
\end{eqnarray*}
and
\begin{eqnarray*}
W_m(\rho, t)={d\over dt}(tH_m(\rho, t)).
\end{eqnarray*}
Then for all $t>0$, we have
\begin{eqnarray}
{d\over dt}W_m(\rho, t)&=&t\int_M \left[\left|{\rm
Hess}~\phi-{g\over t}\right|^2+\Ric_{m,
n}(L)(\nabla \phi, \nabla \phi) \right]\rho d\mu\nonumber\\
& &\ \ \ \ \ \ \ \ \ \ +{t \over m-n}\int_M \left|\nabla \phi\cdot
\nabla f-{m-n\over t}\right|^2 \rho d\mu.\label{Wgeo}
\end{eqnarray}
In particular, if $\Ric_{m, n}(L)\geq 0$, then $W_m(\rho, t)$ is increasing in time $t$ along the geodesic flow on $TP^\infty_2(M, \mu)$.  
\ethm

Note that
\begin{eqnarray*}
{d\over dt}W_m(\rho, t)={d^2\over dt^2}\left(t{\rm Ent}(\rho(t))+mt\log t \right).
\end{eqnarray*}
As a corollary of Theorem \ref{MT2} , we recapture the following result due to Lott-Villani \cite{LoV, Lo2}. 

\bcor\label{Th-LV}  (\cite{LoV, Lo2}) Let $M$ be a compact Riemannian manifold. Suppose that $\Ric_{m, n}(L)\geq 0$. Then $t{\rm Ent}(\rho(t))+mt\log t$ is convex in time $t$ along the geodesic on $P_2(M, \mu)$.
\ecor

%Moreover,  we can also extend  the $W$-entropy formula $(\ref{Wgeo})$ in Theorem \ref{MT2}  to complete Riemannian manifolds with bounded geometry condition, and prove a rigidity theorem for the $W$-entropy for the geodesic flow on the Wasserstein space $P_2(M, \mu)$ over complete Riemannian manifolds with the $CD(0, m)$-condition.  
%
%
%\bthm \label{MT4} ($W$-entropy formula and rigidity theorem for geodesic flow on Wasserstein space over complete Riemannian manifolds)  Let $(M, g)$ be a complete Riemannian manifold with bounded geometry condition, $f\in C^4(M)$ such that $\nabla^k f$ are uniformly bounded on $M$, $1\leq k\leq 4$.  Let $(\rho, \phi)$ be a  be a geodesic in $P_2^\infty(M, \mu)$ satisfying suitable  growth condition as required in Theorem \ref{entropy-geo-noncompact}. Then the $W$-entropy formula  $(\ref{Wgeo})$ in Theorem \ref{MT2} still  holds. In particular,  if the  $CD(0, m)$-condition holds, i.e., $\Ric_{m, n}(L)\geq 0$, then  ${d\over dt}W_m(\rho, t)\geq 0$  on $[0, \infty)$, and $t{\rm Ent}(\rho(t))+mt\log t$ is convex in $t$.  Moreover,  under the assumption $\Ric_{m, n}(L)\geq 0$, ${d\over dt}W_{m}(\rho, t)$ holds at some $t=t_0>0$ if and only if $(M, g)$ is isomeric to $\mathbb{R}^n$, $n=m$, and $(\rho, \phi)=(\rho_n, \phi_m)$. 
%\ethm

Recall that, in our previous papers \cite{Li07, Li12, Li16, LL15, LL17}, inspired by the work of Perelman \cite{P1} and Ni \cite{N1}, we have proved the following $W$-entropy formula for the heat equation associated with the Witten Laplacian on $(M, g, \mu)$. 
 
\bthm\label{MT3} (\cite{Li07, Li12, Li16, LL15, LL17}) Let $M$ be a compact or a complete Riemannian manifold with suitable bounded geometry condition. Let $u$ be a positive solution of the heat equation
\begin{eqnarray}
\partial_t u=Lu.   \label{HeatEqu}
\end{eqnarray}
Define the $H_m$-entropy and the $W_m$-entropy as follows
\begin{eqnarray*}
H_m(u, t)={\rm Ent}(u(t))+{m\over 2}(1+\log(4\pi t)),
\end{eqnarray*}
and
\begin{eqnarray*}
W_m(u, t)={d\over dt}(tH_m(u, t)).
\end{eqnarray*}
Then 
\begin{eqnarray}
{d\over dt}W_m(u, t)&=&2t\int_M \left[\left|{\rm
Hess}~\log u+{g\over 2t}\right|^2+\Ric_{m, n}(L)(\nabla \log u, \nabla \log u) \right] u d\mu\nonumber\\
& &\hskip1cm +{2t\over m-n}\int_M \left|\nabla \log u\cdot \nabla f-{m-n\over t}\right|^2 u d\mu.\label{Wm2}
\end{eqnarray}
In particular, if $\Ric_{m, n}(L)\geq 0$, then $W_{m}(u, t)$ is decreasing  in time $t$ along the heat equation $\partial_t u=Lu$. 
\ethm

As a corollary of Theorem \ref{MT3}, we have the following

\bcor \label{Th-Li}  Let $M$ be a compact or a complete Riemannian manifold with suitable bounded geometry condition. Suppose that $\Ric_{m, n}(L)\geq 0$.  Then $t{\rm Ent}(u(t))+{m\over 2}t\log t$  is convex in time $t$ along the heat equation $\partial_t u=Lu$. 
\ecor

The $W$-entropy formula $(\ref{Wm2} )$ can be regarded as  an analogue of Perelman's $W$-entropy formula for the Ricci flow, and extends Ni's  $W$-entropy formula for the heat equation of the Laplace-Beltrami operator on Riemannian manifolds with non-negative Ricci curvature \cite{N1}. 
%In \cite{LL15}, when $m\in \mathbb{N}$, we gave an alternative proof of the $W$-entropy formula  $(\ref{Wm2} )$  by applying Ni's result to $\widetilde{M}=M\times N$ equipped with the warped product  metric $\widetilde g=g\otimes e^{-{2f\over m-n}}g_{N}$, where $(N, g_N)$ is a $(m-n)$ dimensional compact Riemannian manifold. We also extended the $W$-entropy formula $(\ref{Wm2})$  to the heat equation of the time dependent Witten Laplacian on compact Riemannian manifolds equipped with time dependent metrics and potentials, and to Witten Laplacian on complete Riemannian manifolds with the $CD(K, m)$-condition, see \cite{LL13a, LL14}.

%Note that
%\begin{eqnarray*}
%{d\over dt}W_m(\rho, t)={d^2\over dt^2}\left(t{\rm Ent}(\rho(t))+{m\over 2}t\log t \right).
%\end{eqnarray*}
%
%Thus, as a corollary of Theorem \ref{MT3}, we have the following
%
%\bcor \label{Th-Li}  Let $M$ be a compact Riemannian manifold. Suppose that $Ric_{m, n}(L)\geq 0$.  Then $t{\rm Ent}(u(t))+{m\over 2}t\log t$  is convex in time $t$ on $(0, T]$ for all $T>0$. 
%\ecor

Theorem \ref{MT2} and Theorem \ref{MT3}, Corollary \ref{Th-LV} and Corollary \ref{Th-Li},  have very similar features. We first look at the heat equation case. When $m\in \mathbb{N}$, let $u_m(x, t)={1\over (4\pi t)^{m\over 2}}e^{-{\|x\|^2\over 4t}}$ be the heat kernel of the heat equation $\partial_t u=\Delta u$  on $\mathbb{R}^m$. The Boltzmann-Shannon entropy of the Gaussian heat kernel measure $u_m(x, t)dx$ 
is given by
$${\rm Ent}(u_m(t))=-{m\over 2}(1+\log(4\pi t)).$$
Thus the $H_m$-entropy for the heat equation of the Witten Laplacian is given by\footnote{Following Villani \cite{V1, V2}, we call $H_m(u(t))$ the {\it relative entropy} even though it is slightly different from the classical definition of the relative entropy in probability theory.}
\begin{eqnarray*}
H_m(u(t))={\rm Ent}(u(t))-{\rm Ent}(u_m(t)),
\end{eqnarray*}
and the $W_m$-entropy for the heat equation of the Witten Laplacian is given by the Boltzmann entropy formula in statistical mechanics
\begin{eqnarray}
W_m(u, t):={d\over dt}\left(t[{\rm Ent}(u(t))-{\rm Ent}(u_m(t))]\right).\label{Wm}
\end{eqnarray}
This gives a natural probabilistic interpretation of the $W$-entropy for the heat equation of the Witten Laplacian on Riemannian manifolds. See also \cite{Li12} for the probabilistic interpretation of the Perelman $W$-entropy for the Ricci flow. In \cite{Li12, Li16}, Theorem \ref{MT3} has been extended to complete Riemannian manifolds with bounded geometric condition and a rigidity theorem for the $W_m$-entropy has been proved on complete Riemannian manifolds with  the $CD(0, m)$-condition. More precisely, we have the following

\bthm\label{MT3A}  
%($W$-entropy formula and rigidity theorem for the heat equation of Witten Laplacian on complete Riemannian manifolds \cite{Li12, Li16})  Let $(M, g)$ be a complete Riemannian manifold with bounded geometry condition\footnote{A complete Riemannian manifold $(M, g)$ is called to satisfy  the bounded geometry condition if the Riemannian curvature tensor as well its covariant derivatives are uniformly bounded upto the $4$-th order.} , $f\in C^4(M)$ such that $\nabla^k f$ are uniformly bounded on $M$ for $1\leq k\leq 4$.  Let $u$ be the fundamental solution to the heat equation $(\ref{HeatEqu})$. Then the $W$-entropy formula $(\ref{Wm2})$  in Theorem \ref{MT3} still holds.  In particular,  if the  $CD(0, m)$-condition holds, i.e., $Ric_{m, n}(L)\geq 0$, then ${d\over dt}W_m(u, t)\geq 0$ on $[0, \infty)$,  and $t{\rm Ent}(\rho(t))+{m\over 2}t\log t$ is convex in $t$. 
Suppose that the assumption in Theorem \ref{MT3} holds and $Ric_{m, n}(L)\geq 0$. Then  ${d\over dt}W_m(u, t)\geq 0$ on $[0, \infty)$. Moreover, ${d\over dt}W_m(u, t)=0$  holds at some $t=t_0>0$ if and only if $(M, g)$ is isometric to Euclidean space $\mathbb{R}^n$, $m=n$, and $u(x, t)=u_{m}(x, t)$ is the heat kernel of the heat equation $\partial_t u=\Delta u$  on $\mathbb{R}^m$. 
\ethm

For the geodesic flow on the Wasserstein space, the W-entropy formula can be formulated in the same way. When $m\in \mathbb{N}$, it is easy to check that 
\begin{eqnarray}
\rho_m(x, t)&=&{1\over (4\pi t^2)^{m/2}}e^{-{\|x\|^2\over 4t^2}},\label{rhominfty}\\
\phi_m(x, t)&=&{\|x\|^2\over 2t},\label{phiminfty}
\end{eqnarray}
 is a  geodesic flow on $TP^\infty_2(\mathbb{R}^m)$, where $t>0, x\in \mathbb{R}^m$. Moreover, it serves as the rigidity model. Indeed, the Boltzmann-Shannon entropy of the probability measure $\rho_m(t, x)dx$ is given by 
\begin{eqnarray*}
{\rm Ent}(\rho_m(t))=-{m\over 2}(1+\log(4\pi t^2)),
\end{eqnarray*}
and  the $H_m$-entropy for the geodesic flow on the Wasserstein space $P_2(M, \mu)$ are defined as
\begin{eqnarray}
H_m(\rho(t))={\rm Ent}(\rho(t))-{\rm Ent}(\rho_m(t)).\label{NHW-2a}
\end{eqnarray}
The Boltzmann formula leads us to introduce 
\begin{eqnarray}
W_m(\rho, t):={d\over dt}\left(t[{\rm Ent}(\rho(t))-{\rm Ent}(\rho_m(t))]\right).\label{NHW-3a}
\end{eqnarray} 
Similarly to the case of Theorem \ref{MT3},  we can extend  the $W$-entropy formula $(\ref{Wgeo})$ in Theorem \ref{MT2}  to complete Riemannian manifolds with bounded geometry condition. Moreover, a rigidity theorem can be also proved  for the $W$-entropy for the geodesic flow on the Wasserstein space $P_2(M, \mu)$ over complete Riemannian manifolds with the $CD(0, m)$-condition.  More precisely, we have the following (see also Theorem \ref{Rigidity2} in Section $3$)

\bthm \label{MT4} Under the same condition as in Theorem \ref{MT2}, assume that $(\rho, \phi)$ is a smooth solution (with suitable growth condition) to the transport equation  $(\ref{TA})$ and the Hamilton-Jacobi equation  $(\ref{HJ})$, and suppose that $Ric_{m, n}(L)\geq 0$. Then  ${d\over dt}W_m(\rho, t)\geq 0$  on $[0, \infty)$.  Moreover, ${d\over dt}W_{m}(\rho, t)=0$ holds at some $t=t_0>0$ if and only if $(M, g)$ is isomeric to $\mathbb{R}^n$, $n=m$, and $(\rho, \phi)=(\rho_n, \phi_m)$.  
\ethm

%\bthm\label{MT4}
%Let $(M, g)$ be a complete Riemannian manifold with bounded geometry condition, $f\in C^4(M)$ such that $\nabla^k f$ are uniformly bounded on $M$, $1\leq k\leq 4$.  Let $(\rho, \phi)$ be a smooth solution to the transport equation  $(\ref{TA})$ and the Hamilton-Jacobi equation  $(\ref{HJ})$ satisfying suitable  growth condition as required in Theorem \ref{entropy-geo-noncompact}. Then the $W$-entropy formula  $(\ref{Wgeo})$ in Theorem \ref{MT2} still  holds. In particular,  if the  $CD(0, m)$-condition holds, i.e., $Ric_{m, n}(L)\geq 0$, then  ${d\over dt}W_m(\rho, t)\geq 0$  on $[0, \infty)$, and $t{\rm Ent}(\rho(t))+mt\log t$ is convex in $t$.  Moreover,  under the assumption $Ric_{m, n}(L)\geq 0$, ${d\over dt}W_{m}(\rho, t)$ holds at some $t=t_0>0$ if and only if $(M, g)$ is isomeric to $\mathbb{R}^n$, $n=m$, and $(\rho, \phi)=(\rho_n, \phi_m)$. 
%\ethm

It is then natural to ask the following question:  How to understand  the similarity and the difference between the $W$-entropy formulas for the heat equation on a Riemannian manifold $M$ and for the geodesic flow on the Wasserstein space over $M$? Can we pass through one of them to another one? 

One of possible approaches to answer this question is to use the vanishing viscosity limit method from the heat equation to the Hamilton-Jacobi equation. However,  it seems that  this approach does not work in our situation.  

In this paper, inspired by J.-M.Bismut's work (see \cite{Bis05, Bis10}) on the deformation of hypoelliptic Laplacians on the tangent bundle over Riemannian manifolds,  we introduce a deformation of geometric flows  $(\rho, \phi): [0, T]\rightarrow TP_2(M, \mu)$  by solving the following equations on $TP_2(M, \mu)$ (the tangent bundle over the Wasserstein space $P_2(M, \mu)$), 
\begin{eqnarray}
\partial_t \rho-\nabla_\mu^*(\rho\nabla \phi)&=&0,\label{MF1}\\
c^2\left({\partial_t\phi}+{1\over 2}|\nabla \phi|^2\right)&=&-\phi-\log \rho-1,\label{MF2}
\end{eqnarray}
where  $c\geq 0$. We call $(\rho, \phi)$ the Langevin deformation of flows and prove its analogue of  $W$-entropy formula. 

The Langevin deformation of flows interpolates the geodesic flow on the tangent bundle of  the Wasserstein space $P_2(M, \mu)$ and the heat equation of the Witten Laplacian on the underlying manifold $M$, regarded as the gradient flow of the Boltzmann-Shannon entropy on the Wasserstein space $P_2(M, \mu)$.  Indeed, we can derive  the heat equation and the geodesic flow as the  limit of \eqref{MF1} and \eqref{MF2} in a proper sense 
by taking $c \rightarrow 0$ and $c \rightarrow \infty$ respectively.  See the precise statement of this result 
in Section $5$ and its proof in Section $7$.  

It turns out that the Langevin deformation of flows has a close connection with the compressible Euler equation with damping, such that its wellposedness and the rigorous proof of convergence can be obtained through the results on the corresponding compressible Euler equation with damping. In Section $6$, using the Kato-Majda theory of the hyperbolic quasi-linear systems and the Hamilton-Jacobi theory, we prove that, for any given $c>0$, there exists $T=T_c>0$ such that the Cauchy problem of the system $(\ref{MF1})$ and $(\ref{MF2})$ has a unique smooth solution $(\rho, \phi)\in C^1([0, T], C^\infty(M, \mathbb{R}^+)\times C^\infty(M))$  with given initial data $(\rho_0, \phi_0)\in C^\infty(M, \mathbb{R}^+)\times C^\infty(M)$. 

The following theorem provides us a dissipation formula for the Hamiltonian along the Langevin deformation of flows on $TP_2(M, \mu)$.

\bthm\label{MT0}   Let $(M, g)$  be a compact Riemannian manifold, $f\in C^2(M)$ and $d\mu=e^{-f}dv$. For any $c\in (0, \infty)$ , let $(\phi, \rho)$ be a smooth solution to $(\ref{MF1})$ and $(\ref{MF2})$. 
Define the Hamiltonian and the Lagrangian as follows
%\footnote{We may interpret $H$ as the Hamiltonian of the Langevin deformation in an external potential $V$, and $L$ as the Lagrangian of the Langevin deformation in an external potential $V$. We may also interpret $L$ as the Hamiltonian of the Langevin deformation in an external potential $-V$, and $H$ as the Lagrangian of the Langevin deformation in an external potential $-V$.}
\begin{eqnarray*}
H(\rho, \phi)&=&{c^2\over 2}\int_M |\nabla\phi|^2\rho d\mu+\int_M \rho\log \rho d\mu,\\
L(\rho, \phi)&=&{c^2\over 2}\int_M |\nabla\phi|^2\rho d\mu-\int_M \rho\log \rho d\mu.
\end{eqnarray*}
Then
\begin{eqnarray}
{d\over dt}H(\rho, \phi)&=&-\int_M |\nabla \phi|^2\rho d\mu,\label{monotoH}\\
{d^2\over dt^2}L(\rho, \phi)&=&2\int_M \left[c^{-2}|\nabla \phi+\nabla\log\rho|^2+|{\rm Hess}\phi|^2+\Ric(L)(\nabla\phi, \nabla\phi)\right]\rho d\mu.\label{2ndH}
\end{eqnarray}
In particular, $H$ is always monotonically nonincreasing, and if the $CD(0, \infty)$-condition holds, i.e., $\Ric(L)=\Ric+\nabla^2 f\geq 0$, then $L(\rho, \phi)$ is convex along the Langevin deformation flow $(\rho, \phi)$ defined by $(\ref{MF1})$ and $(\ref{MF2})$. 
\ethm

To state the $W$-entropy type formula for the Langevin deformation of  flows on $TP_2(M, \mu)$, we need first to introduce the reference model in order to define the relative entropy functional as  in $(\ref{NHW-2a})$.  In Section $6$, we prove that there is a special solution to the transport equation $(\ref{MF1})$ and the deformed Hamilton-Jacobi equation $(\ref{MF2})$ on Euclidean space $(\mathbb{R}^m, dx)$ when $m\in \mathbb{N}$. More precisely, let $u: (0, T)\rightarrow  (0, \infty)$ be a smooth solution to the ODE 
\begin{eqnarray}
c^2u''+u'={1\over 2u} \label{u-1}
\end{eqnarray}
where $T>0$ is the lifetime of the solution $u$. Let $\alpha(t)={u'(t)\over u(t)}$, and let $\beta(t)\in C((0, T), \mathbb{R})$  be the unique solution to  the ODE \begin{eqnarray*}
c^2\dot \beta(t)=-\beta(t)-m\log u(t)-{m\over 2}\log(4\pi)+1,\label{bbb}
\end{eqnarray*}
with any given initial data $\beta(0)\in \mathbb{R}$. For $x\in \mathbb{R}^m$ and $t\in (0, T)$, define

\begin{eqnarray*}
\rho_m(x, t)&=&{1\over (4\pi u^2(t))^{m/2}}e^{-{\|x\|^2\over 4u^2(t)}},\\
\phi_m(x, t)&=&{\alpha(t)\over 2}\|x\|^2+\beta(t).
\end{eqnarray*}
Then $(\rho_m, \phi_m)$ is a smooth solution of $(\ref{MF1})$ and $(\ref{MF2})$ on $(\mathbb{R}^m, dx)$.
The above $(\rho_m, \phi_m)$ can be regarded as a reference model to $(\ref{MF1})$ and $(\ref{MF2})$.

\bthm\label{MT1} Let $(M, g)$  be a compact Riemannian manifold, $f\in C^2(M)$, $d\mu=e^{-f}dv$. For any $c\geq 0$,  let  $\alpha(t)$ be as above and 
$(\phi, \rho)$ be a smooth solution to  $(\ref{MF1})$ and  $(\ref{MF2})$. 
Then
\begin{eqnarray}\label{went.2b}
\ \ \ \ & &{d^2\over dt^2} {\rm Ent}(\rho(t))+\left(2\alpha(t)+{1\over c^2}\right){d\over dt}{\rm Ent}(\rho(t))+{1\over c^2}\int_M {|\nabla\rho(t)|^2\over \rho(t)}d\mu+m\alpha^2(t)\nonumber\\
&=& \int_M \left[\left|{\rm Hess}\phi-\alpha(t)g\right|^2+\Ric_{m, n}(L)(\nabla \phi, \nabla \phi) \right]\rho d\mu\nonumber\\
& &\hskip4cm +(m-n)\int_M \left|\alpha(t)+{\nabla \phi\cdot \nabla f\over m-n}\right|^2 \rho d\mu.
%\label{MF3}
\end{eqnarray}
In particular, if the $CD(0, m)$-condition holds, i.e., $\Ric_{m, n}(L)\geq 0$, we have
\begin{eqnarray*}
{d^2\over dt^2} {\rm Ent}(\rho(t))+\left(2\alpha(t)+{1\over c^2}\right){d\over dt}{\rm Ent}(\rho(t))+{1\over c^2}\int_M {|\nabla\rho(t)|^2\over \rho(t)}d\mu+m\alpha^2(t)\geq 0.
\end{eqnarray*}
\ethm

\medskip

\brmq
The limiting cases $c\rightarrow 0$ and $c\rightarrow \infty$ can be specified as follows. 

\begin{itemize}

\item When $c=0$, from $(\ref{MF2})$ we have
\begin{eqnarray}
\phi=\log\rho+1={\delta {\rm Ent}(\rho)\over \delta\rho}, \label{gradEnt}
\end{eqnarray}
which is the $L^2$-derivative of the Boltzmann-Shannon entropy ${\rm Ent}$  on $P_2^\infty(M, \mu)$ equipped with Otto's infinite dimensional Riemannian metric (\cite{Ot, V1, V2}). In this case, $\rho$ satisfies the  heat equation 
\begin{eqnarray}
{\partial_t} \rho=L\rho,\label{BHE1}.
\end{eqnarray}
Equivalently, when $c=0$, $(\rho, \phi)$ can be regarded as the gradient flow of the Boltzmann-Shannon entropy on $P_2^\infty(M, \mu)$ 
equipped with Otto's infinite dimensional Riemannian metric. Inn this case, we have $u(t)=\sqrt{t}$, $\alpha(t)=-{1\over 2t}$, and
\begin{eqnarray*}
\phi_m(x, t)&=&-{\|x\|^2\over 4t}-{m\over 2}\log(4\pi t)+1,\\
\rho_m(x, t)&=&{1\over (4\pi t)^{m/2}}e^{-{\|x\|^2\over 4t}}.
\end{eqnarray*}
The second order entropy dissipation formula reads
\begin{eqnarray*}
{d^2\over dt^2}{\rm Ent}(\rho(t))=2\int_M \left[|{\rm Hess}\log \rho|^2+\Ric(L)(\nabla\log\rho, \nabla\log\rho)\right]\rho d\mu,
\end{eqnarray*}
and we have 
\begin{eqnarray}
& &{d^2\over dt^2} {\rm Ent}(\rho(t))+{2\over t}{d\over dt}{\rm Ent}(\rho(t))+{m\over 2t^2}\nonumber\\
&=&2 \int_M \left[\left|{\rm Hess}\log\rho-{g\over 2t}\right|^2+\Ric_{m, n}(L)(\nabla \log\rho, \nabla \log\rho) \right]\rho d\mu\nonumber\\
& &\hskip1cm +{1\over m-n}\int_M \left|\nabla \log\rho\cdot \nabla f+{m-n\over 2t}\right|^2 \rho d\mu.\label{MF4}
\end{eqnarray}
This is an equivalent form of the $W$-entropy formula in Theorem \ref{MT3}.

\item When $c=\infty$, to make the sense of the equation $(\ref{MF2})$, $\rho$ and 
$\phi$ must satisfies  the transport equation $(\ref{TA})$ and the Hamilton-Jacobi equation $(\ref{HJ})$, i.e., $(\rho, \phi)$  is the geodesic flow on the tangent bundle over the Wasserstein space 
 $P_2^\infty(M, \mu)$. In this case, we have $u(t)=t$, $\alpha(t)={1\over t}$ and

\begin{eqnarray*}
\phi_m(x, t)&=&{\|x\|^2\over 2t},\\
\rho_m(x, t)&=&{1\over (4\pi t^2)^{m/2}}e^{-{\|x\|^2\over 4t^2}}.
\end{eqnarray*}
The second order entropy dissipation formula reads
\begin{eqnarray*}
{d^2\over dt^2}{\rm Ent}(\rho(t))=\int_M \left[|{\rm Hess}\phi|^2+\Ric(L)(\nabla\phi, \nabla\phi)\right]\rho d\mu,
\end{eqnarray*}
and we have
\begin{eqnarray}
& &{d^2\over dt^2} {\rm Ent}(\rho(t))+{2\over t}{d\over dt}{\rm Ent}(\rho(t))+{m\over t^2}\nonumber\\
&=& \int_M \left[\left|{\rm Hess}\phi-{g\over t}\right|^2+\Ric_{m, n}(L)(\nabla \phi, \nabla \phi) \right]\rho d\mu\nonumber\\
& &\hskip1cm +{1\over m-n}\int_M \left|\nabla \phi\cdot \nabla f+{m-n\over t}\right|^2 \rho d\mu.\label{MF4}
\end{eqnarray}
This is an equivalent form of the $W$-entropy formula in Theorem \ref{MT2}.

\end{itemize}
\ermq

In general case $0<c<\infty$, Theorem \ref{MT1} suggests us to introduce a variant of the $W$-entropy as follows: for any $0<t_0<t<\infty$, 
\begin{eqnarray}
W_c(\rho(t))-W_c(\rho(t_0))&:=& {d\over dt}{\rm Ent}(\rho(t))+{1\over c^2}{\rm Ent}(\rho(t))
+2\int_{t_0}^t \alpha(s) {d\over ds} {\rm Ent}(\rho(s))ds\nonumber\\
& &+{1\over c^2}\int_{t_0}^t\int_M {|\nabla\rho(s)|^2\over \rho(s)}d\mu ds.\label{wentc}
\end{eqnarray}
By direct calculation we can verify that
\begin{eqnarray*}
{d\over dt}W(\rho_m(t))=-m\alpha^2(t).
\end{eqnarray*}
In view of this, Theorem \ref{MT1} can be reformulated as follows

\bthm \label{MainMM}  Let $M$ be a compact Riemannian manifold. Under the above notations, we have
\begin{eqnarray}\label{went.2}
\nonumber {d\over dt}(W_{c}(\rho(t))-W_{c}(\rho_m(t)))&=&\int_M  \left|{\rm Hess}~\phi-\alpha(t)g\right|^2\rho d\mu+\int_M \Ric_{m, n}(L)(\nabla\phi, \nabla \phi)\rho d\mu\\
& &\hskip1cm +{1\over m-n}\int_M \left|\nabla f\cdot\nabla\phi+(m-n)\alpha(t)\right|^2\rho d\mu.
\end{eqnarray}
In particular, if $\Ric_{m, n}(L)\geq 0$, then for all $t>0$, we have the comparison inequality
\begin{eqnarray}
{d\over dt}W_{c}(\rho(t))\geq {d\over dt}W_{c}(\rho_m(t)).
\end{eqnarray}
\ethm

%Now we introduce the $W$-entropy for the Langevin deformation of flows as follows: for any $0<t_0<t<\infty$,
%\begin{eqnarray*}
%W_{c}(\rho(t))-W_{c}(\rho(t_0)) :&=& {d\over dt}{\rm Ent}(\rho(t))+\int_{t_0}^t \left( 2\alpha(s)+   {1\over c^2}\right) {d\over ds} {\rm Ent}(\rho(s))ds\\
%& &\hskip2cm  + {1\over c^2}\int_{t_0}^{t}\|\nabla \Ent (\rho(t))\|^2 ds,
%\end{eqnarray*}
%where 
%$$\|\nabla \Ent (\rho(t))\|^2=\int_M {|\nabla \rho(t)|^2\over \rho(t)}d\mu.$$
%
%
%The following result can be viewed as an analogue  of the $W$-entropy formula for the Langevin  deformation of geometric flows on $TP^\infty_2(M, \mu)$, which  
%interpolate the $W$-entropy formula for the geodesic flow on $TP^\infty_2(M, \mu)$ and the $W$-entropy formula for the heat equation of the Witten Laplacian on $(M, \mu)$. 
%

In particular, when $m=n$ and $\mu=v$, we have the following result on compact Riemannian manifolds with standard volume measure. 
\bthm\label{MT1b} Let $(M, g)$  be a compact Riemannian manifold. For any $c\geq 0$,  let $\alpha(t)$ be as above and 
$(\phi, \rho)$ be a smooth solution to the transport equation and the deformed Hamilton-Jacobi equation
\begin{eqnarray}
\partial_t\rho+\nabla\cdot(\rho\nabla\phi)&=&0,\label{TAA}\\
c^2\left(\partial_t\phi+{1\over 2}|\nabla\phi|^2\right)&=&-\phi-\log\rho-1.\label{HJJ}
\end{eqnarray} 
%Then 
%\begin{eqnarray*}
%\ \ \ \ & &{d^2\over dt^2} {\rm Ent}(\rho(t))+\left(2\alpha(t)+{1\over c^2}\right){d\over dt}{\rm Ent}(\rho(t))+n\alpha^2(t)\\
%&=& \int_M \left[\left|{\rm Hess}\phi-\alpha(t)g\right|^2+Ric(\nabla \phi, \nabla \phi) \right]\rho dv+{1\over c^2}\int_M {|\nabla\rho|^2\over \rho}dv.
%\end{eqnarray*}
%In particular, if $Ric\geq 0$, then along the deformed flow $(\rho, \phi)$ defined by $(\ref{TAA})$ and $(\ref{HJJ})$, we have
%\begin{eqnarray*}
%{d^2\over dt^2} {\rm Ent}(\rho(t))+\left(2\alpha(t)+{1\over c^2}\right){d\over dt}{\rm Ent}(\rho(t))+n\alpha^2(t)\geq {1\over c^2}\int_M {|\nabla\rho|^2\over \rho}dv.
%\end{eqnarray*}
\begin{eqnarray}
%\label{went.3}
\nonumber {d\over dt}(W_{c}(\rho(t))-W_{c}(\rho_n(t)))=\int_M  \left|{\rm Hess}~\phi-\alpha(t)g\right|^2\rho d\nu+\int_M \Ric(\nabla\phi, \nabla \phi)\rho d\nu.
\end{eqnarray}
In particular, if $\Ric\geq 0$, then for all $t>0$, we have the comparison inequality
\begin{eqnarray*}
{d\over dt}W_{c}(\rho(t))\geq {d\over dt}W_{c}(\rho_n(t)).
\end{eqnarray*}
\ethm

\brmq \label{remrigidity}
 The entropy formulas in Theorem \ref{MT2}, Theorem \ref{MT3} and Theorem \ref{MainMM} suggest that the specific solution $(\rho_m, \phi_m)$  should play as the reference model for the rigidity theorem of the $W$-entropy on complete 
Riemannian manifolds with the $CD(0, m)$-condition.  Similarly to Theorem \ref{MT3A} and Theorem \ref{MT4}, for any $c\in (0, \infty]$,  we can extend Theorem \ref{MainM}  to a class of smooth solutions $(\rho, \phi)$ on weighted complete Riemannian manifolds $(M, g, f)$ with natural bounded geometry condition.
%\footnote{Here we say that $(M, g)$ satisfies the bounded geometry condition if the Riemannian curvature tensor ${\rm Riem}$ and its covariant derivatives $\nabla^k {\rm Riem}$ are uniformly bounded on $M$, $k=1, 2, 3$.}, and $f\in C^4(M)$ such that $\nabla f\in C_b^3(M)$. 
We can therefore expect that the following rigidity theorem holds:  Let  
$M$ be a complete Riemannian manifold with 
natural bounded geometry condition and with $CD(0, m)$-condition, i.e., $Ric_{m, n}(L)\geq 0$. Then ${d\over dt}W_c(\rho(t))={d\over dt}W_c(\rho_m(t))$ holds at some $t=t_0>0$, if and only if $M$ is 
isometric to $\mathbb{R}^n$, $m=n$, $f$ is a constant, $(\rho, \phi)=({\rho}_m, {\phi}_m)$. To save the length of the paper, we omit the detail of the technical part of the proof of this result and will do it in a future work. 
%
%Theorem \ref{MT2}, Theorem \ref{MT3} and Theorem \ref{MainM} can be extended to compact Riemannian manifolds with $CD(K, m)$-condition, where $m\geq n$ and $K\in \mathbb{R}$  are two constants.  To save the length of paper, we will do this in a forthcoming paper, in which we will also introduce the Langevin  deformation of flows on the tangent bundle over the Wasserstein space which interpolate the  gradient flow of the R\'enyi entropy (i.e., the  nonlinear porous media equation) and the geodesic flow on the Wasserstein space. Moreover we will extend the $W$-entropy formula to the new Langevin deformation of flows. 
%

\ermq

As we have already mentioned, the Langevin deformation of flows converges in a proper sense to the gradient flow of the Boltzmann-Shannon entropy on $TP_2^\infty(M, \mu)$ when $c\rightarrow 0$,  and converges in a proper sense  to the geodesic flow on $TP_2^\infty(M, \mu)$ when $c\rightarrow \infty$.  
In view of this, the  $W$-entropy formula~\eqref{went.2} indeed interpolates the $W$-entropy formula~\eqref{Wgeo} for the geodesic flow on $TP_2^\infty(M, \mu)$ and the $W$-entropy formula~\eqref{Wm2} for the heat equation of the Witten Laplacian on $(M, \mu)$.

In \cite{LoV, St1, St2, St3}, Lott-Villani and Sturm proved that the Boltzmann entropy ${\rm Ent}$ is $K$-convex along the geodesic on the Wasserstein space $P_2(M, \mu)$ if and only if the $CD(K, \infty)$-condition holds, i.e., $\Ric(L)\geq K$, and the R\'enyi entropy 
$S_N(\rho\mu)=-\int_M \rho^{-1/N}d\mu$ is convex along the geodesic on $P_2(M, \mu)$ for all $N\geq m$ if and only if the $CD(0, m)$-condition holds. In  view of this, we would like to raise the following conjecture for the characterization of the $CD(0, m)$-condition on complete Riemannian manifolds, which can be regarded  as the converse of Corollary \ref{Th-LV} 
and Corollary \ref{Th-Li}. 

\textbf{Conjecture}.  Let $(M, g)$ be a compact Riemannian manifold or a complete Riemannian manifold with bounded geometry condition, $f\in C^\infty(M)$ with $\nabla f\in C_b^\infty(M)$. Suppose that the $W_m$-entropy associated to the heat equation of the Witten Laplacian or the optimal transport problem is non-decreasing in $t$. Then the $CD(0, m)$-condition holds, i.e., $\Ric_{m, n}(L)\geq 0$.

 Now we mention some related work in the literature. In \cite{MT}, McCann and Topping proved the contraction property of the $L^2$-Wasserstein distance between solutions of the backward heat equation on closed manifolds equipped with the Ricci flow, which extends a previous result for the Fokker-Planck equation on Euclidean space  (see Otto \cite{Ot})  and on complete Riemannian manifolds with suitable Bakry-Emery curvature condition (see Sturm and von Renesse \cite{StR} ). See also \cite{T1, T2}. In \cite{Lo2}, Lott further proved two convexity results of the Boltzmann-Shannon type entropy along the geodesics on the Wasserstein space over closed manifolds equipped with Ricci flow, which are  closely related to Perelman's results on the monotonicity of the $\mathcal{F}$ and $\mathcal{W}$-entropy functionals for Ricci flow. In \cite{LL13b}, the authors extended Lott's convexity results to the Wasserstein space on compact Riemannian manifolds equipped with Perelman's Ricci flow.

%\begin{remark} Let $U_\mu(\rho)=\int_M U(\rho)d\mu$, where $U: [0, \infty)\rightarrow \mathbb{R}$ is a continuous convex function with $U(0)=0$. By Remark 2 and Remark 3 in \cite{Lo2}, if the $CD(0, m)$-condition holds, then 
% $tU_\mu(\rho(t))+mt\log t$ is convex along the geodesics flow $\rho(t)d\mu$ on $P_2^\infty(M, \mu)$. Similarly to Theorem \ref{MT2} and Theorem \ref{MT3}, we can introduce the $W$-entropy with respect to the R\'enyi entropy $U_\mu$ and prove the $W$-entropy formula for the geodesic flow and the gradient flow of $U_\mu(\rho)=\int_M U(\rho)d\mu$ on the Wasserstein space, i.e., the porous media equation $\partial_t \rho=\nabla_\mu^*(\rho \nabla U'(\rho))$. Moreover, we can also  extend the $W$-entropy formula to  the Langevin deformation of geometric flows which interpolates the 
% porous media  equation on $(M, g, \mu)$  and  the geodesic flow on $P_2^\infty(M, \mu)$.  Due to the limit of the length of the paper, we will develop this in a forthcoming paper. 
%\end{remark}

The paper is organized as follows. In Section $2$,  we recall some elementary facts about the infinite dimensional Riemannian geometry on the Wasserstein space over Riemannian manifolds. In Section $3$ we introduce the $W$-entropy for the geodesic flow on the Wasserstein space and prove the $W$-entropy formula  $(\ref{Wgeo})$ for the geodesic flow, i.e., Theorem \ref{MT2}. A rigidity theorem is proved. In Section $4$ we introduce the Langevin deformation 
of geometric flows on the tangent bundle of a complete Riemannian manifold. In Section $5$ we introduce the Langevin deformation of geometric flows on the Wasserstein space over Riemannian manifolds, and prove its corresponding $W$-entropy formula in Section $6$. In Section $7$, we provide the rigorous proof of the convergence of  Langevin deformation for $c\rightarrow 0$ and $c\rightarrow \infty$ respectively. 

To end this section, let us mention that the previous version \cite{LL-flow17} of the present paper  has been posted on arXiv (arxiv1604.02596, 2016) but has not been submitted, while Theorem \ref{main.convergence} and its proof in Section $7$ are new. Theorem \ref{MT2}  has been also stated (without giving the proof) and the Langevin deformation  of flows (with the positive sign in front of $\nabla V$) has been  introduced in our previous paper \cite{LL18SCM}. 

\section{Otto's calculus on Wassertsien space over weighted Riemannian manifolds}\label{prelimi}

Let $(M, g)$ be a complete Riemannian manifold, $f\in C^2(M)$, and $d\mu=e^{-f}d\nu$, where $d\nu$ denotes the volume measure on $(M, g)$. Integration by parts formula shows that, for all $u\in C_0^\infty(M)$ and $X\in C_0^\infty(M, TM)$, we have
$$\int_M \langle X, \nabla u\rangle d\mu=\int_M \nabla_\mu^* X ud\mu,$$
where $\nabla_\mu^*$ denotes the $L^2$-adjoint of $\nabla$ with respect to $\mu$, and is given by
$$
\nabla_\mu^* X=-\nabla\cdot X+\langle\nabla f, X\rangle.$$
The Witten Laplacian on $(M, g)$ with respect to $\mu$ is defined by
$$
L=-\nabla_\mu^*\nabla.$$
More precisely, we have
$$L=\Delta-\nabla f\cdot\nabla.$$

By Bakry-Emery \cite{BE}, the Bochner-Weitzenb\"ock formula holds
\begin{eqnarray*}
L|\nabla u|^2-2\langle \nabla u, \nabla Lu\rangle=2|{\rm Hess} u|^2+2\Ric(L)(\nabla u, \nabla u),\  \  \ \forall\ u\in C^\infty(M),
\end{eqnarray*}
where 
\begin{eqnarray*}
\Ric(L):=\Ric+\nabla^2 f
\end{eqnarray*}
is the infinite dimensional  Bakry-Emery Ricci curvature associated with the Witten Laplacian $L$.  Following \cite{BE}, we say that $CD(K, \infty)$-condition holds if and only if $\Ric(L)\geq K$.

Let
$$
P^{\infty}(M, \mu) = \{\rho d\mu: \rho \in C^{\infty}(M), \rho\geq 0, \int_M \rho d\mu = 1 \}.
$$
For all $\rho d\mu \in
P^{\infty}(M, \mu)$, the tangent space at $\rho d\mu$ is given by
\begin{eqnarray*}
T_{\rho d\mu}P^\infty(M, \mu)=\{s\in C^\infty(M): \int_M sd\mu=0\}.
\end{eqnarray*}
By solving the Poisson equation $-\rho L\phi+\nabla \phi\cdot\nabla \rho=s$, there exists a unique function $\phi\in C^\infty(M)$ (up to a constant) such that
$$
s=V_\phi:=\nabla_\mu^*(\rho \nabla\phi).
$$
Following \cite{Lo2}, $V_\phi$ can be identified as the  vector field on $P_2^\infty(M, \mu)$ defined by
\begin{eqnarray*}
(V_\phi F)(\rho d\mu)=\left.{\partial\over \partial \varepsilon}\right|_{\varepsilon=0} F(\rho d\mu+\varepsilon \nabla^*_{\mu}(\rho\nabla \phi)d\mu),
\end{eqnarray*}
where $F\in C^\infty(P_2^\infty(M), \mathbb{R})$.

Let $P_2^\infty(M, \mu)$ be the Wasserstein space of probability measures $\rho d\mu\in P^\infty(M)$ with finite second moment
$$
\int_M d^2(o, \cdot)\rho d\mu<\infty,$$
where $o\in M$ is a fixed point. Similarly to Otto \cite{Ot}, we can introduce the infinite dimensional Riemannian meric on $T_{\rho d\mu}P_2^\infty(M, \mu)$ as follows
$$
\langle\langle s_1, s_2\rangle\rangle=\int_M \langle\nabla \phi_1, \nabla \phi_2\rangle \rho d\mu, \ \ \ \forall s_i=V_{\phi_i}\in T_{\rho d\mu}P_2^\infty(M, \mu),\ i=1, 2.
$$
The tangent space of $P_2(M, \mu)$ at $\rho d\mu$, denoted by $T_{\rho d\mu}P_2(M, \mu)$, is defined as the completion of $T_{\rho d\mu}P_2^\infty(M, \mu)$ with respect to the norm 
$\|s\|^2:=\int_M |\nabla \phi|^2 \rho d\mu$ for $s=\nabla_\mu^*\cdot(\rho \nabla\phi)$. 

The Wasserstein distance between $\mu_1$ and $\mu_2$ is defined by
\begin{eqnarray*}
W_2^2(\mu_1, \mu_2)=\inf\limits_{\pi}\int_{M\times M} d^2(x, y)d\pi(x, y),
\end{eqnarray*}
where $\pi\in \prod(\mu_1, \mu_2)$, i.e., $\pi$ is a probability measure on $M\times M$ such that
\begin{eqnarray*}
\int_M \pi(\cdot, dy)=\mu_1, \ \ \ \ \int_M \pi(dx, \cdot)=\mu_2.
\end{eqnarray*}

The following result extends  Benamou and Brenier's result  \cite{BB} from Euclidean space to complete Riemannian manifolds. 

\bthm\label{Th-BBM} Let $(M, g)$ be a complete Riemannian manifold. Let $\mu_0=\rho_0d\mu$,  $\mu_1=\rho_1d\mu$ be two probability measures in $P_2(M, \mu)$. Then
\begin{eqnarray*}
W_2^2(\mu_0, \mu_1)=\inf\int_{M}\int_0^1 |\nabla \phi(x, t)|^2\rho(x, t) d\mu(x)dt,
\end{eqnarray*}
where
$\rho: M\times [0, 1]\rightarrow [0, \infty)$ and $\phi:  M\times [0, 1]\rightarrow \mathbb{R}$ satisfy
\begin{eqnarray}
\label{geod1}\partial_t\rho-\nabla_\mu^*(\rho \nabla \phi)=0,\\
\label{geod2}\partial_t\phi+{1\over 2}|\nabla \phi|^2=0,\\
\label{in.value.geo}\rho(\cdot, 0)=\rho_0, \ \ \ \rho(\cdot, 1)=\rho_1.
\end{eqnarray}
\ethm
\bpf
The proof is analogue of the one in \cite{BB}.
\epf

 The function $\phi$ in Theorem \ref{Th-BBM} is called the potential function, and $v=\nabla \phi$ can be considered as the velocity of the curve $\rho(\cdot, t)dv$ in ${P}_2^\infty(M, \mu)$.  The transport equation~\eqref{geod1} and the Hamilton-Jacobi equation~\eqref{geod2} describe the geodesic $\rho_s d\mu$ which links $\rho_0d\mu$ and $\rho_1 d\mu$ in $P_2^\infty(M, \mu)$ equipped with Otto's infinite dimensional Riemannian metric.

Following Lott~\cite{Lo2}, we can prove the entropy dissipation formula along the geodesics on the Wasserstein space over a compact Riemannian manifold
with weighted volume measure. When the potential function is constant,  it is due to Lott \cite{Lo2}. To save length of the paper we omit the details of the proof.

\bprop\label{thm-A}
Let $(M, g)$  be a compact Riemannian manifold, $f\in C^2(M)$.
Let $\rho: [0, 1]\times M\rightarrow\mathbb{R} $ be a positive solution of the
transport equation
\begin{eqnarray}
\partial_t\rho-\nabla_\mu^*(\rho\nabla\phi)=0.\label{TA1}
\end{eqnarray}
%where for fixed $t\in [0, 1]$, 
%$\phi(\cdot, t)\in C^\infty(M)$ can be viewed as the velocity vector of
%the smooth curve $t\in [0, 1]\rightarrow \rho(t)d\mu\in P_2^\infty(M, \mu)$ in $T_{\rho(t)
%d\mu}P_2^\infty(M, \mu)$. 
Let ${\rm Ent}(\rho(t)):= \int_M \rho \log
\rho d\mu$, then we have
\begin{eqnarray}
{d\over dt}{\rm Ent}(\rho(t))=\int_M \langle \nabla \rho, \nabla \phi \rangle d\mu,\label{A-1}
\end{eqnarray}
and
\begin{eqnarray}
\nonumber \frac{d^2}{dt^2}{\rm Ent}(\rho(t)) = -\int_M L \rho
(\partial_t\phi + \frac{1}{2}|\nabla \phi|^2) d\mu + \int_M (|{\rm Hess}
\phi|^2 + Ric(L)(\nabla \phi, \nabla \phi))\rho d\mu. 
\\ \label{A-2}
\end{eqnarray}
\eprop

As a direct consequence of Proposition~\ref{thm-A}, we have the following
\bthm \label{Ent-Hess} For the geodesic $(\rho, \phi)$ on $P_2^\infty(M, \mu)$,  we have
\begin{eqnarray*}
{d^2\over dt^2} {\rm Ent}(\rho(t))=\int_M (|{\rm Hess} ~\phi|^2 + \Ric(L)(\nabla \phi, \nabla \phi))\rho d\mu.
\end{eqnarray*}
\ethm

In view of Theorem $\ref{Ent-Hess}$, the Hessian of the Boltzmann-Shannon entropy functional ${\rm Ent}$ on $P_2^\infty(M,
\mu)$ equipped with Otto's infinite dimensional Riemannian metric is given by
\begin{eqnarray}\label{hess.ent}
{\rm Hess}_{P_2^\infty(M, \mu)} {\rm Ent}(V_\phi, V_\phi)=\int_M (|{\rm Hess}~\phi|^2+\Ric(L)(\nabla\phi, \nabla \phi))\rho d\mu.
\end{eqnarray}

As a corollary of Theorem $\ref{Ent-Hess}$, we have the following result due to Lott-Villani \cite{LoV, Lo2}, Sturm-von Renesse \cite{StR} and Sturm \cite{St2, St3}.

\bcor\label{Convex-1} If $\Ric(L)\geq 0$, then ${d^2\over dt^2} {\rm Ent}(\rho(t))\geq 0$, i.e., ${\rm Ent}$ is convex along geodesic in $P_2^\infty(M, \mu)$.
\ecor

\section{The $W$-entropy formula for the geodesic flow on Wasserstein space}

In this section, we introduce the $W$-entropy and prove its variational formula along the geodesic flow on the Wasserstein space over compact Riemannian manifolds with weighted volume measure, i.e., Theorem \ref{MT2}. We will also compare the $W$-entropy formula in Theorem \ref{MT2}  with  the $W$-entropy formula for the  heat equation of the Witten Laplacian on compact Riemannian manifolds (i.e., Theorem \ref{MT3} ), and then introduce the $W$-entropy for the optimal transport problem on compact or complete Riemannian manifolds with weighted volume measure.

\subsection{Proof of Theorem~\ref{MT2}}
\bpf
By Theorem~\ref{Ent-Hess}, we have
%\begin{eqnarray*}
%{d\over dt}{\rm Ent}(\rho(t))&=&-\int_M L\phi \rho d\mu,\\
%{d^2\over dt^2} {\rm Ent}(\rho(t))
%&=&\int_M (|{\rm Hess}\phi|^2+Ric(L)(\nabla\phi, \nabla \phi))\rho d\mu.
%\end{eqnarray*}
%Thus
\begin{eqnarray*}
{d\over dt}W_m(\rho, t)&=&t{d^2\over dt^2} {\rm Ent}(\rho(t))+2{d\over dt}{\rm Ent}(\rho(t))+{m\over t}\\
&=&t\int_M (|{\rm Hess}~\phi|^2+\Ric(L)(\nabla\phi, \nabla \phi))\rho d\mu-2\int_M L\phi \rho d\mu+{m\over t}\\
%&=&t\int_M \left(\left|{\rm Hess}\phi-{g\over t}\right|^2+{2\over t}\Delta \phi+Ric(L)(\nabla\phi, \nabla \phi)\right)\rho d\mu-t\int_M L\phi \rho d\mu+{m-n\over t}\\
&=&t\int_M \left[\left|{\rm Hess}~\phi-{g\over t}\right|^2+\Ric(L)(\nabla \phi, \nabla \phi)\right] \rho d\mu+2\int_M \nabla \phi \cdot \nabla f \rho d\mu+{m-n\over t}.
\end{eqnarray*}
Note that
\begin{eqnarray*}
\ \ \ & &\Ric(L)(\nabla \phi, \nabla\phi)+{2\over t}\nabla \phi\cdot\nabla f+{m-n\over t^2}\\
%&=& Ric_{m, n}(\nabla\phi, \nabla \phi)+{|\nabla\phi \cdot\nabla f|^2\over m-n}+{2\over t}\nabla \phi\cdot\nabla f+{m-n\over t^2}\\
&=& \Ric_{m, n}(L)(\nabla\phi, \nabla \phi)+{1\over m-n}\left|\nabla \phi\cdot \nabla f+{m-n\over t}\right|^2.
\end{eqnarray*}
Thus
\begin{eqnarray*}
& &{d\over dt}W_m(\rho, t) \\
&= &t\int_M \left[\left|{\rm Hess}~\phi-{g\over t}\right|^2+\Ric_{m, n}(L)(\nabla \phi, \nabla \phi) \right]\rho d\mu+{t
\over m-n}\int_M \left|\nabla \phi\cdot \nabla f+{m-n\over t}\right|^2 \rho d\mu.
\end{eqnarray*}
This proves  the $W$-entropy formula $(\ref{Wgeo})$ in Theorem \ref{MT2}. 
\epf

As a corollary of Theorem \ref{MT2}, we can recapture the following result due to Lott-Villani \cite{LoV}. See also Lott \cite{Lo2}.

\bcor (i.e., Corollary \ref{Th-LV}) If $\Ric_{m, n}(L)\geq 0$, then $t{\rm Ent}+m t\log t$ is convex in $t$ along the geodesic in $P_2^\infty(M, \mu)$.
\ecor

\subsection{The $W$-entropy for  the geodesic flow on Wasserstein space}

Let $\rho_n(t, x)={e^{-{\|x\|^2\over 4t^2}}\over (4\pi t^2)^{n/2}}$ be the special solution  of the transport equation~\eqref{geod1} with the velocity $\phi_n(t, x)={\|x\|^2\over 2t}$ on the Euclidean space $\mathbb{R}^n$. By calculus,   the Boltzmann-Shannon entropy of $\rho_n(t, x)dx$ with respect to the Lebesgue measure on $\mathbb{R}^n$ is given by 
\begin{eqnarray*}
{\rm Ent}(\rho_n(t))=-{n\over 2}(1+\log(4\pi t^2)).
\end{eqnarray*}

Let $(\rho(t), \phi(t))$ be smooth solution to $(\ref{TA})$ and $(\ref{HJ})$. Define the Boltzmann-Shannon entropy
\begin{eqnarray}
{\rm Ent}(\rho(t))=\int_M \rho \log \rho d\mu,\label{NHW-0}
\end{eqnarray}
By Theorem \ref{Ent-Hess}, we have
\begin{eqnarray}
{d\over dt}{\rm Ent}(\rho(t))&=&-\int_M \langle\nabla\rho, \nabla\phi\rangle d\mu=\int_M L\phi \rho d\mu,\label{NHW-1}\\
{d^2\over dt^2}{\rm Ent}(\rho(t))&=&-\int_M [|\nabla^2\phi|^2+\Ric(L)(\nabla\phi, \nabla\phi)]\rho d\mu.\label{NHW-2}
\end{eqnarray}
%
%Now, for the optimal transport problem on compact Riemannian manifolds $(M, g)$ with the weighted measure $d\mu=e^{-f}dv$, the geodesic flow on $P_2^\infty(M, \mu)$ equipped with Otto's infinite dimensional Riemannian metric is given by the transport equation together with the Hamilton-Jacobi equation
%\begin{eqnarray}
%& &\partial_t \rho+\nabla_\mu^*(\rho\nabla \phi)=0,\label{HJ2}\\
%& &\partial_t \phi+{1\over 2}|\nabla \phi|^2=0.\label{HJ1}
%\end{eqnarray}
%Define the Boltzmann-Shannon entropy
%\begin{eqnarray}
%{\rm Ent}(\rho(t))=\int_M \rho \log \rho d\mu,\label{NHW-0}
%\end{eqnarray}
%By Theorem \ref{Ent-Hess}, we have
%\begin{eqnarray}
%{d\over dt}{\rm Ent}(\rho(t))&=&-\int_M \langle\nabla\rho, \nabla\phi\rangle d\mu=\int_M L\phi \rho d\mu,\label{NHW-1}\\
%{d^2\over dt^2}{\rm Ent}(\rho(t))&=&-\int_M [|\nabla^2\phi|^2+Ric(L)(\nabla\phi, \nabla\phi)]\rho d\mu.\label{NHW-2}
%\end{eqnarray}
Following Perelman \cite{P1}, we introduce the $W$-entropy for the geodesic flow $(\rho, \phi)$ on $TP_2(M, \mu)$ as 
\begin{eqnarray}
W_m(\rho, t):={d\over dt}\left(tH_m(\rho, t)\right),\label{NHW-3}
\end{eqnarray}
where
\begin{eqnarray}
H_m(\rho, t)={\rm Ent}(\rho(t))-{\rm Ent}(\rho_m(t))\label{NHW-4}
\end{eqnarray}
is the difference between the Boltzmann-Shannon entropy of the probability measure $\rho d\mu$ with respect to the weighted volume measure $\mu$ on $(M, g)$ and  the Boltzmann-Shannon entropy of the probability measure $\rho_m(t, x)dx$ with respect to the Lebesgue measure $dx$ on $\mathbb{R}^m$. 

Substituting $(\ref{NHW-1})$ into $(\ref{NHW-4})$ and $(\ref{NHW-3})$ we have
\begin{eqnarray}
W_m(\rho, t)=\int_M \left[tL \phi-\log \rho-{\rm Ent}(\rho_m(t))\right]\rho d\mu.\label{NHW-5}
\end{eqnarray}
Moreover, we can reformulate Theorem \ref{MT2} as follows
\begin{eqnarray*}
{d\over dt} W_m(\rho, t)&=&t \int_M \left[\left|{\rm Hess}~\phi-{g\over t}\right|^2+\Ric_{m, n}(L)(\nabla \phi, \nabla \phi) \right]\rho d\mu\\
& &\hskip2cm +{t
\over m-n}\int_M \left|\nabla \phi\cdot \nabla f+{m-n\over t}\right|^2 \rho d\mu.
\end{eqnarray*}

In view of this,  Theorem~\ref{MT2} can be reformulated as follows.
\bthm\label{Th-Main2}  Let $M$ be a compact Riemannian manifold. Let $(\rho, \phi)$ be a smooth geodesic flow on $TP_2^\infty(M, \mu)$, i.e., $\rho$ is a smooth solution to the transport equation~\eqref{TA} and $\phi$ is a smooth solution to the Hamilton-Jacobi equation~\eqref{HJ}.
Let
\begin{eqnarray*}
W_m(\rho, t)={d\over dt}\left(t[{\rm Ent}(\rho(t))-{\rm Ent}(\rho_m(t))]\right)
\end{eqnarray*}
be the $W$-entropy associated to the optimal transport on $(M, g, \mu)$.
Then
\begin{eqnarray*}
{d \over d t}W_m(\rho, t)&=&t\int_M \left[\left|{\rm
Hess}~\phi-{g\over t}\right|^2+\Ric_{m,
n}(L)(\nabla \phi, \nabla \phi) \right]\rho d\mu\\
& &\ \ \ \ \ \ \ \ \ \ +{t \over m-n}\int_M \left|\nabla \phi\cdot
\nabla f-{m-n\over t}\right|^2 \rho d\mu.
\end{eqnarray*}
In particular, if $\Ric_{m, n}(L)\geq 0$, then the Helmholtz free energy $S_m=t({\rm Ent}(\rho(t))-{\rm Ent}(\rho_m(t)))$ associated with \eqref{TA} and \eqref{HJ} is convex in time $t$.
\ethm

\bcor (i.e., Corollary  \ref{Th-LV},   \cite{LoV, Lo2}) Let $M$ be a compact Riemannian manifold. Suppose that $\Ric_{m, n}(L)\geq 0$. Then $t{\rm Ent}(\rho(t))+mt\log t$ is convex in time $t$ along the geodesic flow $(\rho(t), \phi(t))$ on $TP_2(M, \mu)$.
\ecor
\bpf
 By Theorem \ref{MT2}, if $\Ric_{m, n}(L)\geq 0$, $t({\rm Ent}(\rho(t))-{\rm Ent}(\rho_m(t))=t{\rm Ent}(\rho(t))+{mt\over 2}[\log (4\pi t^2)+1]$ is convex in $t$ along the geodesic $\rho(t)$  on $P_2^\infty(M, \mu)$. Note that
\begin{eqnarray*}
{d^2\over dt^2}\left( t({\rm Ent}(\rho(t))-{\rm Ent}(\rho_m(t))\right)={d^2\over dt^2} (t{\rm Ent}(\rho(t))+mt\log t).
\end{eqnarray*}
Hence $t{\rm Ent}+mt\log t$ is convex in $t$ along  the geodesic $\rho(t)$  on $P_2^\infty(M, \mu)$.  For the general case of non smooth geodesic on $P_2(M, \mu)$, see Lott \cite{Lo2}. 
\epf

In particular, taking $f=0$, $m=n$ and $L=\Delta$, we have the following 

\bthm Let $M$ be a compact Riemannian manifold. Let $\phi$ and $\rho$ be a smooth solution to the Hamilton-Jacobi equation and the transport equation
\begin{eqnarray}
\partial_t\phi+{1\over 2}|\nabla\phi|^2&=&0,\label{HJ-3}\\
\partial_t \rho+\nabla\cdot(\rho\nabla\phi)&=&0.\label{HJ-4}
\end{eqnarray}
Let
\begin{eqnarray}
W_n(\rho, t)={d\over dt}\left(t[{\rm Ent}(\rho(t))-{\rm Ent}(\rho_n(t))]\right).\label{Sn}
\end{eqnarray}
Then
\begin{eqnarray}
{d\over d t}W_n(\rho, t)=t\int_M \left[\left|{\rm
Hess}~\phi-{g\over t}\right|^2+\Ric(\nabla \phi, \nabla \phi) \right]\rho d\mu.\label{Sn2}
\end{eqnarray}
In particular, if $\Ric\geq 0$, then the $W$-entropy $W_n$ associated with $(\ref{HJ-3})$ and $(\ref{HJ-4})$ is increasing in time $t$, and $t{\rm Ent}+nt\log t$ is convex in $t$ along the geodesic flow $(\rho(t), \phi(t))$  on $TP_2(M, \nu)$.
\ethm

\subsection{The case of complete noncompact Riemannian manifolds}

To extend Theorem \ref{MT2} to complete Riemannian manifolds with bounded geometry condition, we need the following 
\bthm\label{entropy-geo-noncompact} Let $M$ be a complete Riemannian manifold, and $f\in C^2(M)$. Suppose that $\Ric(L)=\Ric+\nabla^2f$ is uniformly bounded on $M$, i.e., there exists a constant $C>0$ such that $|\Ric(L)|\leq C$. 
Let $\rho$ and $\phi$ be smooth solutions to the transport equation~\eqref{TA} and the Hamilton-Jacobi equation~\eqref{HJ}, and satisfying the following growth condition
\begin{eqnarray*}
 \int_M \left[|\nabla\log \rho|^2+|\nabla\phi|^2+|\nabla^2\phi|^2+|L\phi|^2+|\nabla L\phi|^2\right]\rho d\mu<\infty,
\end{eqnarray*}
and there exist  a point $o\in M$, and some functions $C_i\in C([0, T], \mathbb{R}^+)$ and $\alpha_i\in C([0, T], \mathbb{R}^+)$ such that
\begin{eqnarray*}
C_1(t)e^{-\alpha_1(t) d^2(x, o)}\leq \rho(x, t)\leq C_2(t)e^{\alpha_2 (t)d^2(x, o)}, \ \ \forall x\in M, t\in [0, T],
\end{eqnarray*} 
and
\begin{eqnarray*}
\int_M d^4(x, o)\rho(x, t)d\mu<\infty, \  \ \ \forall t\in [0, T].
\end{eqnarray*}
Then the entropy dissipation formulas hold
\begin{eqnarray*}
\partial_t \int_M \rho \log \rho d\mu&=&\int_M \nabla \phi \cdot \nabla \rho d\mu=-\int_M L\phi \rho d\mu,\\
\partial_t^2\int_M \rho \log \rho d\mu&=&-\int_M [|\nabla^{2}\phi|^2+\Ric(L)(\nabla\phi, \nabla\phi)] \rho d\mu.
\end{eqnarray*}
\ethm
\bpf
 Let $\eta_k$ be an increasing sequence of functions in $C_0^\infty(M)$ such that $0\leq \eta_k\leq 1$, $\eta_k=1$ on $B(o, k)$ , $\eta_k=0$ on $M\setminus B(o, 2k)$, and $\|\nabla\eta_k\|\leq {1\over k}$. By standard argument and integration by parts, we have 
\begin{eqnarray*}
\partial_t \int_M \rho \log\rho \eta_k d\mu
%&=&\int_M \partial_t (\rho\log \rho)\eta_k d\mu\\
%&=&\int_M \partial_t \rho(\log \rho+1)\eta_k d\mu\\
&=&\int_M \nabla_\mu^*(\rho\nabla\phi)(\log \rho+1)\eta_k d\mu\\
&=&\int_M \rho\nabla\phi \cdot \nabla \log \rho \eta_k d\mu+\int_M \rho\nabla\phi \cdot(\log \rho+1)\nabla\eta_k d\mu\\
&=:&I_1+I_2.
\end{eqnarray*}
Under the assumption of theorem, the Lebesgue dominated convergence theorem yields
\begin{eqnarray*}
I_1&=&\int_M \rho\nabla\phi \cdot \nabla \log \rho \eta_k d\mu\rightarrow \int_M \nabla\phi\cdot \nabla \rho d\mu.
\end{eqnarray*}
%if $\nabla\phi\cdot \nabla \rho \in L^1(M, \mu)$, which is the case if $\int_M |\nabla\phi|^2\rho d\mu<\infty$ and $\int_M {|\nabla\rho|^2\over \rho}d\mu<\infty$.  
and 
\begin{eqnarray*}
I_1=\int_M \nabla_\mu^*(\eta_k\nabla\phi)\rho d\mu=-\int_M \eta_k L\phi  \rho d\mu+\int_M \nabla\eta_k\cdot\nabla\phi \rho d\mu\rightarrow -\int_M L\phi \rho.
\end{eqnarray*}
%if $\int_M |L\phi|^2\rho d\mu<\infty$ and $\int_M |\nabla\phi|^2\rho d\mu<\infty$. 
Thus
\begin{eqnarray*}
 \int_M \nabla\phi\cdot \nabla \rho d\mu=-\int_M L\phi \rho d\mu.
\end{eqnarray*}

On the other hand, since $\int_M |\nabla\phi|^2\rho d\mu<\infty$ and $\int_M |\log \rho+1|^2\rho d\mu<\infty$, we have $\int_M |\log \rho+1||\nabla\phi|\rho d\mu<\infty$.  By Lebesgue dominated convergence theorem, we have 
\begin{eqnarray*}
I_2=\int_M \rho\nabla\phi \cdot (\log\rho +1)\nabla \eta_k d\mu\rightarrow 0.
\end{eqnarray*}
This proves that \begin{eqnarray*}
\partial_t \int_M \rho \log\rho \mu=\int_M \nabla\phi \cdot \nabla \rho d\mu=-\int_M L\phi \rho d\mu.
\end{eqnarray*}

By  standard argument, we have

\begin{eqnarray*}
\partial_t \int_M L\phi \rho \eta_k d\mu 
%&=&\int_M \partial_t(L\phi \rho)\eta_kd\mu\\
%&=&\int_M L\partial_t \phi \rho \eta_k d\mu+\int_M L\phi \partial_t \rho \eta_kd\mu\\
&=&\int_M L\left(\partial_t \phi+{1\over 2}|\nabla\phi|^2\right)\rho \eta_k d\mu-{1\over 2}\int_M L|\nabla\phi|^2\rho \eta_k d\mu+\int_M L\phi \partial_t \rho \eta_kd\mu\\
&=&-{1\over 2}\int_M L|\nabla\phi|^2\rho \eta_k d\mu+\int_M L\phi \nabla_\mu^*(\rho \nabla\phi)\eta_kd\mu\\
&=&I_3+I_4.
\end{eqnarray*}

By the weighted Bochner formula and $|\Ric(L)|\leq C$, we have
\begin{eqnarray*}
\int_M |L|\nabla\phi|^2|\rho d\mu&=&2\int_M \left|\nabla\phi \cdot\nabla L\phi+|\nabla^2\phi|^2+\Ric(L)(\nabla\phi, \nabla\phi)\right|\rho d\mu\\
&\leq&2 \int_M\left[|\nabla\phi| |\nabla L\phi|+|\nabla^2\phi|^2+C|\nabla\phi|^2\right|\rho d\mu\\
&<&\infty
\end{eqnarray*}
provided that 
\begin{eqnarray*}
\int_M [|\nabla\phi|^2+|\nabla L\phi|^2+ |\nabla^2\phi |^2]\rho d\mu<\infty.
\end{eqnarray*}
Thus
\begin{eqnarray*}
I_3\rightarrow -{1\over 2}\int_M L|\nabla\phi|^2\rho d\mu.
\end{eqnarray*}

On the other hand, as $\int_M [ |L\phi|^2+|\nabla \phi|^2+ |\nabla L\phi|^2 ]\rho d\mu<\infty$, we have
\begin{eqnarray*}
I_4
&=&\int_M L\phi \nabla_\mu^*(\rho \nabla\phi)\eta_kd\mu\\
&=&\int_M \nabla (\eta_K L\phi)\cdot \rho \nabla\phi d\mu\\
&=&\int_M [L\phi\nabla\eta_k \cdot \nabla \phi+\eta_k \nabla L\phi\cdot\nabla\phi]\rho d\mu\\
&\rightarrow&\int_M \nabla L\phi\cdot\nabla \phi \rho d\mu.
\end{eqnarray*}

In summary, under the assumption of theorem, we have
\begin{eqnarray*}
\partial_t \int_M L\phi \rho d\mu&=& -{1\over 2}\int_M L|\nabla\phi|^2\rho d\mu+\int_M \nabla L\phi\cdot\nabla \phi \rho d\mu\\
&=&-\int_M [|\nabla^2\phi|^2+\Ric(L)(\nabla\phi, \nabla\phi)]\rho d\mu.
\end{eqnarray*}
This finishes the proof of Theorem \ref{entropy-geo-noncompact}. 
\epf

Based on Theorem \ref{entropy-geo-noncompact}, Theorem \ref{MT2} can be extended to complete Riemannian manifolds as follows. 

\bthm\label{MTT2} Let $M$ be a complete Riemannian manifold with bounded geometry condition. Under the same condition as in Theorem \ref{entropy-geo-noncompact}, the $W$-entropy formula in Theorem \ref{MT2} remains true. 
\ethm
\bpf The proof is similar to the one of Theorem \ref{MT2}. 
\epf

\subsection{The rigidity theorem for the $W$-entropy formula}

Similarly to the monotonicity  and rigidity theorem (i.e., Theorem \ref{MT3A})  of  the $W$-entropy for the heat equation associated with the Witten Laplacian over complete Riemannian manifolds with bounded geometry condition,  which is the gradient flow of the Botlzmann-Shannon entropy on the Wasserstein space, we have the following monotonicity  and rigidity theorem of  the $W$-entropy for the geodesic flow  on the Wasserstein space over complete Riemannian manifolds with bounded geometry condition.

\bthm \label{Rigidity2} Let $M$ be a complete Riemannian manifold with bounded geometry condition  and with  $Ric_{m, n}(L)\geq 0$.  Suppose that $(\rho, \phi)$ is a smooth  geodesic flow on $TP_2^\infty(M, \mu)$, i.e., $(\rho(t), \phi(t), t\in [0, \infty))$ is a smooth solution to the transport equation~\eqref{TA} and the Hamilton-Jacobi equation~\eqref{HJ} satisfying the growth condition as required in Theorem \ref{entropy-geo-noncompact}. Then $W_m(\rho, t)$ is increasing in $t$, and $t{\rm Ent}(\rho(t))+mt\log t$ is convex in $t$. Moreover, ${d\over dt}W_{m}(\rho, t)=0$ at some $t=\tau>0$ if and only if $M$ is isomeric to $\mathbb{R}^n$, $n=m$, and $(\rho(t), \phi(t))=(\rho_n(t), \phi_n(t))$ for all $t\geq 0$. That is to say, the Euclidean space $\mathbb{R}^n$ equipped with the Gaussian measure $N(0, t^2)$ is the rigidity model for the $W$-entropy associated with the geodesic flow on $TP_2^\infty(M, \mu)$ with natural  growth condition at infinity  over complete Riemannian manifold with bounded geometry condition  and with  $Ric_{m, n}(L)\geq 0$.   
\ethm
\bpf 
The proof is as the same as the one of  Theorem $2.5$ in \cite{Li12} (i.e., Theorem \ref{MT4} as stated in Section $1$).  For the convenience of the reader, we give the detail here.  Indeed, by the explicit expression of ${d\over dt}W_m(\rho, t)$ in Theorem \ref{MT2}, we see that  ${d\over dt}W_m(\rho, t)=0$ holds at some $t=\tau$ if and only if $\nabla^2 \phi(\cdot, \tau)={g\over \tau}$, $\Ric_{m, n}(L)(\nabla\phi, \nabla \phi)=0$ and $\nabla\phi\cdot\nabla f+{m-n\over \tau}=0$. Thus, $\phi$ is a strict convex function on $M$, which implies that $M$ is diffeomorphic to $\mathbb{R}^n$. Integrating along the shortest geodesics between $x_0$ and 
$x$  on $M$ shows that
\begin{eqnarray*}
2\tau (\phi (x,\tau)-\phi (x_0,\tau))=r^2(x_0,x),\ \ \ \ \forall x\in M,
\end{eqnarray*} where $x_0$ is the minimum point of $\phi(\cdot, \tau)$. 
This yields 
\begin{eqnarray*}
\Delta r^2(x_0, x)=2n, \ \ \ \forall \ x\in M,
\end{eqnarray*}
which implies that $(M, g)$ is isometric to the Euclidean space $(\mathbb{R}^n, (\delta_{ij}))$. 
By the generalized Cheeger-Gromoll splitting theorem  (see Theorem 1.3, p. 565, [6]), we can derive that $f$ must be a constant and $m=n$. 

Thus $\phi(\cdot, \tau)\in C^\infty(\mathbb{R}^n)$ satisfies $\nabla^2\phi(x, \tau)={\delta_{ij}\over \tau}$. This yields $\nabla\phi(x, \tau)={x\over \tau}$ under the assumption $\nabla\phi(0, \tau)=0$, and $\phi(x, \tau)={\|x\|^2\over 2\tau}$ up to an additional  constant.
By the Hopf-Lax formula for the solution to the Hamilton-Jacobi equation~\eqref{geod2} with $\phi(x, \tau)={\|x\|^2\over 2\tau}$, we have 
\begin{eqnarray*}
\phi(t, x)=\inf\limits_{y\in \mathbb{R}^n}\left\{{\|x\|^2\over 2\tau}+{\|x-y\|^2\over 2(t-\tau)}\right\}= {\|x\|^2\over 2t}, \ \ \ \ \ \forall \ t> \tau, \  x\in \mathbb{R}^n.
\end{eqnarray*}
By the uniqueness of the smooth solution  to the Hamilton-Jacobi equation~\eqref{geod2}, we see that $\phi(x, t)={\|x\|^2\over 2t}$ for all $t>0$. 
Solving the transport equation~\eqref{geod1} with the initial data $\lim\limits_{t\rightarrow 0}\rho(t, x)=\delta_0(x)$, we have 
$$\rho(x, t) = {e^{-{\|x\|2\over 4t^2}}\over (4\pi t^2)^{n/2}}, \ \ \ \ \forall t>0,\  x\in \mathbb{R}^n.$$
This finishes the proof of Theorem \ref{Rigidity2}.
\epf

\subsection{Comparison between Theorem \ref{MT2} and Theorem \ref{MT3}}

Theorem \ref{MT2} has the same feature as Theorem \ref{MT3}. Their proofs are also quiet similar. Both of them are based on the dissipation formulas of the first order and the second order derivatives of the Boltzmann-Shannon entropy along the heat equation (i.e., the gradient flow of the Boltzmann-Shannon entropy) and the geodesic flow on the Wasserstein space $P_2^\infty(M, \mu)$.  

Theorem \ref{MT2} and Theorem \ref{MT3} lead us to the following observation: 
On compact Riemannian manifolds $(M, g)$ with the weighted measure $d\mu=e^{-f}d\nu$, if the $CD(0, m)$-condition holds, then the relative Boltzmann-Shannon entropy
 $H_m(\rho, t)={\rm Ent}(\rho(t))-{\rm Ent}(\rho_m(t))$ is convex along the geodesic flow on the Wasserstein space $P_2^\infty(M, \mu)$, and the relative Boltzmann-Shannon entropy 
 $H_m(u, t)={\rm Ent}(u(t))-{\rm Ent}(u_m( t))$ is convex along the backward gradient flow of  ${\rm Ent}(u)=-\int_M u\log u d\mu$ on the Wasserstein space $P_2^\infty(M, \mu)$. This 
leads us to raise the question whether there is an essential reason for which the Boltzmann-Shannon entropy share the same convexity property along the geodesic flow and the gradient flow on $P_2^\infty(M, \mu)$.

On the other hand, there is a difference between the $W$-entropy formula for the heat equation of the Witten Laplacian on Riemmanian manifold and the W-entropy formula for the geodesic flow on the Wasserstein space, i.e., ${g\over 2t}$ appears in $(\ref{Wm2})$, while ${g\over t}$ appears in 
$(\ref{Wgeo})$, and their rigidity models are also different (see~\cite{Li16, LL15, LL17} and Theorem \ref{Rigidity2}). An intuitive interpretation for this difference can be given as follows. The heat kernel  of the Laplacian on $\mathbb{R}^m$ is the transition probability of Brownian motion starting from time $0$ to time $t$. The mean square displacement (the variance of the distance that the ``Brownian particle'' moving along its trajectory) on  $\mathbb{R}^m$ during the time interval $[0,t]$ is given by $\mathbb{E}[|B_t|^2] = mt$. While for the geodesic flow~\eqref{geod1} and \eqref{geod2}, if we assume that the velocity of the ``light particle'' has the unit speed along each direction, the distance (denoted by $|X_t|$) of the ``light particle'' moving along the geodesic during time interval $[0,t]$ is indeed $|X_t| = t$. Hence, the mean square displacement of the ``light particle'' moving along the geodesic during time interval $[0,t]$ is $\mathbb{E}[|X_t|^2] = t^2$. This explains intuitively why the rigidity model for the $W$- entropy for the heat equation of the Witten Laplacian on complete Riemannian manifolds with the $CD(0, m)$-condition is the Gaussian space $(\mathbb{R}^m, g_0, N(0, t{\rm Id}))$, while the rigidity model for the $W$-entropy for the geodesic flow on the Wasserstein space over complete Riemannian manifolds with the $CD(0, m)$-condition is the Gaussian space $(\mathbb{R}^m, g_0, N(0, t^2{\rm Id}))$. Here $g_0$ denotes the Euclidean metric on $\mathbb{R}^m$, ${\rm Id}$ is the unit matrix on $\mathbb{R}^m$, and $N(0,t{\rm Id})$ denotes the Gaussian distribution on $\mathbb{R}^m$ with mean zero and variance $t{\rm Id}$.

\section{Langevin deformation of flows on finite dimensional manifolds}\label{finite.flow}

In this section we introduce the Langevin deformation of geometric flows on the tangent bundle $TM$ over a complete Riemannian manifold $(M,g)$, which interpolates the geodesic flow on $TM$ and the gradient flow of 
a potential function on $M$. Our work has been inspired by Bismut \cite{Bis05, Bis10} who introduced the deformation of a family of hypoelliptic Laplacians on $TM$, which is the infinitesimal generator of the Langevin diffusion process 
interpolating the geodesic flow on $TM$ and the Brownian motion on $M$. The ideas and results in this section will be extended in Section $6$ to the infinite dimensional Wasserstein space over compact Riemannian manifolds.

\subsection{The construction of the deformation of flows}

We first describe J.-M. Bismut's  idea for the construction of a family of hypoelliptic Laplacians on the tangent bundle over Riemannian manifolds (\cite{Bis05, Bis10}). Let $c>0$ be a parameter, let  $(x_t, v_t)$ be the Langevin diffusion process on the tangent bundle $TM$ over a complete Riemannian manifold $M$ which solves the following stochastic differential equation 
\begin{eqnarray}
\dot x&=&{v\over c}, \label{Bis-1}\\
dv&=&-{v\over c^2}dt+{dw_t\over c},\label{Bis-2}
\end{eqnarray}
where  $dw_t$ denotes the Stratonovich  differential of Brownian motion $w_t$ on $M$. This is the stochastic differential equation for  the Langevin hypoelliptic  diffusion process $(x_t, v_t)$ on the tangent bundle $TM$ over $M$.  The position process $x_t$ satisfies the  Langevin stochastic differential equation 
\begin{eqnarray}
c^2\ddot x&=&-\dot x+\dot w_t,\label{Bis-3}
\end{eqnarray}
where  $\dot w_t$ denotes the  Stratonovich derivation of Brownian path $w_t$ on $M$.
As was pointed out by Bismut  \cite{Bis05, Bis10},  taking $c\rightarrow 0$, the limiting process $x_t$ is the Brownian motion on $M$, ie., 
\begin{eqnarray*}
\dot x=\dot w_t,
\end{eqnarray*}
and when $c\rightarrow \infty$,  to make sense the Langevin stochastic differential equation $(\ref{Bis-3})$, the limiting process $x_t$ must satisfy the geodesic equation 
\begin{eqnarray*}
\ddot x=0. 
\end{eqnarray*}
Thus the Langevin diffusion processes $(x_t, v_t)$ provide  a deformation of geometric flows which interpolate the geodesic flows $\ddot x=0$ on the tangent bundle $TM$ over $M$ and the Brownian motion $x_t=w_t$ on the underlying Riemannian manifold $M$.  

Let $V$ be a smooth function on $M$. Instead of the above Langevin diffusion processes on $TM$, let us introduce the Langevin deformation of geometric flows on $TM$
\begin{eqnarray}
\dot x&=& v, \label{Bis-4}\\
c^2\dot v&=&-v - \nabla V(x). \label{Bis-5}
\end{eqnarray}
Then $x_t$ satisfies the second order ordinary differential equation 
\begin{eqnarray}
c^2\ddot x = -\dot x -\nabla V(x). \label{Bis-6}
\end{eqnarray}
The equation $(\ref{Bis-6})$ is indeed the Newton-Langevin equation which describes the motion of particles with mass $c^2$ moving in a fluid with friction coefficient $1$ and with an external potential $V$.  In view of the classical ODE theory, for fixed $c > 0$ and under a suitable condition on $V$(for example $\nabla V$ is Lipschitz), the above equations admit a unique local solution. It defines a family of geometric flows on $TM$ which interpolates the geodesic flow $\ddot{x}=0$ 
and the gradient flow $\dot x=-\nabla V(x)$. Indeed,
similarly to Bismut's situation, we can rigorously prove that when $c\rightarrow 0$,  the limiting flow $x_t$ is the gradient flow of $V$, i.e., 
\begin{eqnarray*}
\dot x_t= - \nabla V(x_t),
\end{eqnarray*}
and when $c\rightarrow \infty$,  the limiting flow $x_t$ must satisfy the geodesic equation 
\begin{eqnarray*}
\ddot x=0. 
\end{eqnarray*}
In view of this, $(x_t, v_t)$ is a deformation of geometric flows on $TM$ which interpolate the geodesic flows $\ddot x=0$ on the tangent bundle $TM$ over $M$ and the gradient flow $\dot x_t=-\nabla V(x_t)$ on the underlying Riemannian manifold $M$.  For the statement and proof of the convergence result, see Section 7. 

Following Bismut \cite{Bis05, Bis10} and Villani \cite{V1}, we use a Hamiltonian point of view to give an interpretation of  the above deformation of flows.  Let
\begin{eqnarray}
H(x, v)={{c^{2} |v|^2}\over 2}+V(x) \label{Bis-H}
\end{eqnarray}
be the Hamiltonian energy of a particle moving in the tangent bundle, where $-V$ is the external potential. Then
\begin{eqnarray*}
\nabla H(x, v)=\left({\partial H\over \partial x}, {\partial H\over \partial v}\right)^\tau=\left(\nabla V(x), c^{2}v\right)^\tau.
\end{eqnarray*}
Let
\[ A=\left( \begin{array}{ccc}
0 & c^{-2}\Id \\
-c^{-2}\Id & -c^{-4}\Id \\
\end{array} \right).\]
Then
\[ A \nabla H(x, v)=\left(v, -c^{-2}\nabla V(x) - c^{-2}v\right)^\tau.\]
Thus the deformation flow $(x_t, v_t)$ can be regarded as the ``{\it  $A$-Hamiltonian flow}'', given by 
\begin{eqnarray}
{d\over dt}\left( \begin{array}{ccc}
x\\
 v\\
\end{array} \right)=A\nabla H\left( \begin{array}{ccc}
x\\
 v\\
\end{array} \right).\label{Bis-AH}
\end{eqnarray}

Therefore, $(x_t, v_t)$ admits the dissipative property as the gradient flow: the Hamiltonian 
is monotone along the Langevin deformation flow, while the Lagrangian is $K$-convex if $-V$ is $K$-convex, i.e., $\nabla^2 V\leq -K$. More precisely, we have

\bprop \label{Th-B-2} Let $(x_t, v_t)$ be a smooth solution to ~\eqref{Bis-4} and ~\eqref{Bis-5} on $TM$.
Let

\begin{eqnarray*}
H(x, v)={c^{2}|v|^2\over 2}+V(x),\\
L(x, v)={c^{2}|v|^2\over 2}-V(x).
\end{eqnarray*}
Then
\beqna
{d\over dt} H(x_t, v_t) =-|v|^2\leq 0.
\eeqna
and 
\begin{eqnarray*}
{d^2\over dt^2} L(x, v)= 2|\dot v |^{2} - 2\nabla^{2} V(v, v).
\end{eqnarray*}
In particular, if $-V$ is $K$-convex, i.e., $\nabla^2 V\leq -K$, we have
\begin{eqnarray*}
{d^2\over dt^2} L(x, v)\geq 2|\dot v |^{2} + 2K|v|^2.
\end{eqnarray*}
\eprop
\bpf
By direct calculation, we have
\begin{eqnarray*}
{d\over dt} H(x_t, v_t)&=& c^{2}v\cdot \dot v+\nabla V\cdot \dot x\\
&=&v \cdot \left(- v - \nabla V(x)\right)+\nabla V\cdot v\\
&=&-|v|^2.
\end{eqnarray*}
On the other hand, we have 
\begin{eqnarray*}
{d\over dt} L(x, v)&=& c^{2}v\cdot \dot v - \nabla V\cdot \dot x\\
&=&v \cdot \left(- v - \nabla V(x)\right) - \nabla V\cdot v\\
&=&-|v|^2 - 2\nabla V\cdot v,
\end{eqnarray*}
and
\begin{eqnarray*}
{d^2\over dt^2} L(x, v)&=&  -2v\cdot \dot v-2\nabla^2 V(v, v)-2\nabla V\cdot \dot v\\
&=&2c^{-2}v\cdot (v+\nabla V(x))-2\nabla^2 V(v, v)+2c^{-2}\nabla V\cdot  (v+\nabla V(x))\\
%&=&2c^{-2}|v|^{2}  + 4c^{-2}v \cdot \nabla V(x)+2c^{-2}|\nabla V|^2- 2\nabla^{2} V(v, v) \\
&=&  2c^{-2}|v +\nabla V |^{2} - 2\nabla^{2} V(v, v)\\
&=&2|\dot v|^2- 2\nabla^{2} V(v, v).
\end{eqnarray*}
%In particular, if $\nabla^2 V\leq -K$, we have
%\begin{eqnarray*}
%{d^2\over dt^2} L(x, v)\geq 2|\dot v |^{2} + 2K|v|^2.
%\end{eqnarray*}

\epf

\subsection{The $W$-entropy formula for the deformation flow on $TM$}

In this subsection, we introduce a variant of the $W$-entropy functional and to prove its monotonicity along the deformed flows $(x_t, v_t)$ on $TM$. Our first observation is the following 

\begin{proposition} \label{VHH}For any $c>0$, let $(x_t, v_t)$ be the Langevin deformation flow on $TM$ defined by $(\ref{Bis-4})$ and $(\ref{Bis-5})$. Then 
\begin{eqnarray}
{d^2\over dt^2}V(x)+  c^{-2}{d\over dt}V(x) + c^{-2}|\nabla V|^2 = \nabla^{2} V(v, v). \label{Bis-VH1}
%\left({d^2\over dt^2}+{1\over c^2}{d\over dt}\right)V(x)&=&{1\over c^2}\left[ \nabla^2V(v, v)+|\nabla V|^2\right]. \label{Bis-VH1}
%\left({d^2\over dt^2}+{2\over c^2}{d\over dt}\right)H(x, v)&=&{2\over c^2}\left[ \nabla^2V(v, v)+|\nabla V|^2\right].\label{Bis-VH2}
\end{eqnarray}
In particular, if $\nabla^2 V \geq K$, where $K\in \mathbb{R}$ is a constant, we have
\begin{eqnarray}\label{ent.ineq.finite.dim}
{d^2\over dt^2}V(x)+  c^{-2}{d\over dt}V(x) + c^{-2}|\nabla V|^2 \geq  K|v|^{2}. 
%\left({d^2\over dt^2}+{2\over c^2}{d\over dt}\right)H(x, v)&\geq &{2\over c^2}\left[ K|v|^2+|\nabla V|^2\right].\label{Bis-VH2b}
\end{eqnarray}
\end{proposition}
\bpf Indeed, a simple calculation yields 
\beqnas
{d\over dt} V(x(t)) &=& \langle \nabla V, \dot{x}\rangle, \\
{d^{2}\over d^{2}t} V(x(t)) &=&  \nabla^{2} V(\dot{x}, \dot{x}) + \langle \nabla V, \ddot{x}\rangle\\
&=&  \nabla^{2} V(v, v) -  c^{-2}\langle \nabla V, v + \nabla V \rangle,
%&=&  \nabla^{2} V(v, v) -  c^{-2}|v + \nabla V|^{2} +  c^{-2}\langle v, v + \nabla V \rangle.
\eeqnas
%\begin{eqnarray*}
%{d\over dt}V(x)&=&\nabla V(x)\cdot \dot x=\nabla V(x)\cdot {v\over c},\\
%{d^2\over dt^2}V(x)&=&\nabla^2 V(\dot x, \dot x)+\nabla V \cdot \ddot x\\
%&=&{1\over c^2}\nabla^2 V(v, v)+{1\over c^2}\nabla V\cdot (-\dot x+\nabla V). 
%\end{eqnarray*}
Hence
\begin{eqnarray*}
{d^2\over dt^2}V(x)+  c^{-2}{d\over dt}V(x) + c^{-2}|\nabla V|^2 = \nabla^{2} V(v, v) \geq K|v|^{2}.
\end{eqnarray*}
%Similarly, we have
%\begin{eqnarray*}
%\left({d^2\over dt^2}+{2\over c^2}{d\over dt}\right)H&=&2\left|{v\over c^2}-{\nabla V\over c}\right|^2+2\nabla^2 V\left({v\over c}, {v\over c}\right)+{2\over c^2}\left(-{|v|^2\over c^2}+2\nabla V\cdot {v\over c} \right)\\
%&=&{2\over c^2}\left[|\nabla V|^2+\nabla^2V(v, v)\right].
%\end{eqnarray*}
This finishes the proof.
\epf

Thus from~\eqref{Bis-VH1} we see that ${d^2\over dt^2}V(x)+  c^{-2}{d\over dt}V(x) + c^{-2}|\nabla V|^2$ is a quantity that is invariant for $c \in [0, \infty]$. This leads us to introduce the $W$-entropy for  the Langevin deformation flow on $TM$ such that its time derivative is given by
\beqna\label{went.finite.dim}
{d\over dt}{W}_{c}(x, v) :=  {d^2\over dt^2}V(x)+  c^{-2}{d\over dt}V(x) + c^{-2}|\nabla V|^2.
\eeqna
Thus if $\nabla^2V\geq K$, by~\eqref{ent.ineq.finite.dim}, we have for all $c \geq 0$,
$$
{d\over dt}{W}_{c}(x, v)  = \nabla^2V \geq K.
$$

%
%The following remark indicates that $W_{c}$ is indeed an interpolation between the $W$-entropy for the gradient flow and that for the geodesic flow. 
\brmq
In fact, by the convergence result (Theorem \ref{cthm1}) in Section $7$, we see that the quantity $c^{-2}{d\over dt}V(x) + c^{-2}|\nabla V|^2$ is well defined when $c$ approaches $0$ and $\infty$. 
More precisely, when $c =0$, $\dot x=-\nabla V(x)$, we have $\ddot x=-\nabla^2 V \dot x=\nabla^2 V v$. Hence
\begin{eqnarray*}
{d^2\over dt^2}V(x)=\nabla^{2} V(\dot x, \dot x)+\nabla V\cdot \ddot x=2\nabla^{2} V(v, v).
\end{eqnarray*}
Moreover, when $\nabla V$ are $\nabla^2V$ is $K$-Lipschitz and uniformly bounded, then there exists a constant $T>0$ such that, as $c \rightarrow 0$, we have
$$
\ddot{x} = -c^{-2}(\dot x + \nabla V) \rightarrow  - \nabla^{2} V \cdot v \ \ ({\rm in \ the \ uniform\ convergence\ on}\ [0, T]),
$$ we derive that, as $c\rightarrow 0$, then in the uniform convergence on $[0, T]$, we have
\beqnas
c^{-2}{d\over dt}V(x) + c^{-2}|\nabla V|^2 =  \langle \nabla V, c^{-2}(\dot{x} +  \nabla V)\rangle  \rightarrow -\nabla^{2} V(v, v),
\eeqnas
which implies
\beqna\label{finite.went.1}
\lim\limits_{c\rightarrow 0}{d\over dt}{W}_{c}(x, v)
% &=&\lim\limits_{c\rightarrow 0}\left[ {d^2\over dt^2}V(x)+  c^{-2}{d\over dt}V(x) + c^{-2}|\nabla V|^2\right]\nonumber\\
%&=&\lim\limits_{c\rightarrow 0}\left[ \nabla^{2} V(\dot x,\dot x)+\nabla V\cdot \ddot x+ \langle \nabla V, c^{-2}(\dot{x} +  \nabla V)\rangle\right] \nonumber\\
=\nabla^{2} V(v, v)={d\over dt}{W}_{0}(x, v). 
\eeqna
On the other hand,  when $c \rightarrow \infty$, we have in the uniform convergence on $[0, T]$
\begin{eqnarray*}
{d^2\over dt^2}V(x)&=&\nabla^2V(\dot x, \dot x)+\langle \nabla V, \ddot x\rangle\\
&=&\nabla^2V(v, v)-\langle \nabla V, c^{-2}(\dot x+\nabla V)\rangle\\
&\rightarrow& \nabla^2V(v, v),
\end{eqnarray*}
and
$$
c^{-2}{d\over dt}V(x) + c^{-2}|\nabla V|^2  =  \langle \nabla V, c^{-2}(\dot{x} +  \nabla V)\rangle \rightarrow 0,
$$
which implies
\beqna\label{finite.went.2}
\lim\limits_{c\rightarrow \infty}{d\over dt}{W}_{c}(x, v) 
%&=&\lim\limits_{c\rightarrow \infty}\left[ {d^2\over dt^2}V(x)+  c^{-2}{d\over dt}V(x) + c^{-2}|\nabla V|^2\right]\nonumber\\
=\nabla^{2} V(v, v)= {d\over dt}{W}_{\infty}(x, v).
\eeqna
The convergences \eqref{finite.went.1} and \eqref{finite.went.2} show that our definition of $W_{c}$ is indeed an  interpolation between the $W$-entropy for the gradient flow and that for the geodesic flow. 
\ermq

\section{Langevin deformation of flows on Wasserstein space}

In this section we extend the idea and results in Section~\ref{finite.flow} to the infinite dimensional Wasserstein space over compact Riemannian manifolds, and introduce the Langevin deformation of flows on  the Wasserstein space $P_2^\infty(M, \mu)$.  Based on the connection between the Langevin deformation of flows on $P_2^\infty(M, \mu)$ and the the compressible Euler equation with damping, we prove the well-posedness of the Langevin deformation of flows. Finally we show the convergence results of the Langevin deformation of flows when $c$ approaches $0$ and $\infty$ respectively.

\subsection{Langevin deformation of flows on $TP_2^\infty(M, \mu)$}

Let $M$ be a compact Riemannian manifold, and $P_2^\infty(M, \mu)$ be the smooth Wasserstein space over $M$ equipped with the weighted volume measure $d\mu=e^{-f}dv$. Let $\rho: [0, T]\rightarrow P_2^\infty(M, \mu)$ be a smooth curve.
Since $\dot\rho\in T_{\rho d\mu}P_2^\infty(M, \mu)$, there exists a function $\phi$ on $M$ such that
\begin{eqnarray*}
\dot\rho=\nabla_\mu^* (\rho\nabla \phi).
\end{eqnarray*}
Equivalently, $\rho$ satisfies the transport equation with velocity $\phi$ on $P_2^\infty(M, \mu)$. 
\begin{eqnarray}
\partial_t\rho-\nabla_\mu^* (\rho\nabla \phi)=0. \label{tran.eq}
\end{eqnarray}
By Otto's infinite dimensional Riemannian metric on $T_{\rho d\mu}P_2^\infty(M, \mu)$, we have
\begin{eqnarray*}
\|\dot\rho\|^2=\int_M |\nabla\phi|^2\rho d\mu.
\end{eqnarray*}
Let $V\in C^1(\mathbb{R})$. Define the Hamiltonian and Lagrangian on $TP_2^\infty(M, \mu)$ by
\begin{eqnarray}\label{hamiltonian.infinite}
H(\rho, \dot\rho) &=& {\|\dot \rho\|^2\over 2}+V(\rho),\label{hamiltonian.infinite}
\\
L(\rho, \dot\rho) &=& {\|\dot \rho\|^2\over 2}-V(\rho).\label{Langrangian.infinite}
%&=& {1\over 2}\int_M |\nabla \phi|^2\rho d\mu-V(\rho).
\end{eqnarray}
Extending the idea in Section~\ref{finite.flow}, let us introduce the following ODE on $TP_2^\infty(M, \mu)$
\begin{eqnarray}
\label{infinite.flow.1}\partial_t \rho&=&v,\\
\label{infinite.flow.2}\partial_t v&=&-v-\nabla V(\rho)
\end{eqnarray}
The first equation is indeed  the transport equation~\eqref{tran.eq}, 
%\begin{eqnarray*}
%\partial_t \rho+\nabla_\mu^*(\rho\nabla \phi)=0,
%\end{eqnarray*}
and the second equation can be written as
\begin{eqnarray*}
\nabla_{\dot \rho}\dot \rho=-\dot \rho-\nabla V(\rho).
\end{eqnarray*}
According to Otto \cite{Ot}, we have
\begin{eqnarray*}
\nabla V(\rho)=\nabla_\mu^*\left(\rho\nabla {\delta V\over \delta \rho}\right),
\end{eqnarray*}
where ${\delta V\over \delta \rho}$ is the $L^2$-derivative of $V$ with respect to $\rho$, and the result in Lott~\cite{Lo2} leads to
\begin{eqnarray}\label{LOT}
\nabla_{\dot \rho}\dot \rho=\nabla_\mu^*\left(\rho\nabla\left(\partial_{t}\phi +{1\over 2}|\nabla \phi|^2\right)\right).
\end{eqnarray}
Hence the second equation reads
\begin{eqnarray*}
\nabla_\mu^*\left(\rho\nabla\left(\partial_{t}\phi+{1\over 2}|\nabla \phi|^2\right)\right)=-\nabla_\mu^*(\rho\nabla \phi)-\nabla_\mu^*\left(\rho \nabla {\delta V\over \delta \rho}\right),
\end{eqnarray*}
which is equivalent to (up to a additive constant)
\begin{eqnarray}
\partial_{t}\phi + {1\over 2}|\nabla \phi|^2=-\phi-{\delta V\over \delta \rho}. \label{vvv}
\end{eqnarray}
Now for any parameter $c>0$, we can introduce the Langevin deformation of flows on $TP_2^\infty(M, \mu)$ as follows
\begin{eqnarray}
\partial_t \rho-\nabla_\mu^*(\rho\nabla \phi)&=&0, \label{ww} \\
c^2\left({\partial\phi\over \partial t}+{1\over 2}|\nabla \phi|^2\right)&=&-\phi - {\delta V\over \delta \rho}.\label{vvvv}
\end{eqnarray}

In Section 6.2, we will prove the well-posedness of the Langevin deformation of flows on $TP_2^\infty(M, \mu)$. From its definition, we can see that, when $c\rightarrow 0$,  we may formally get 
\begin{eqnarray}
\phi= - {\delta V\over \delta \rho}, \label{grad1}
\end{eqnarray}
which yields that $\rho$ is the gradient flow of $V$ on the Wasserstein space $P_2^\infty(M, \mu)$
\begin{eqnarray}
\partial_t \rho = -\nabla_\mu^*\left(\rho\nabla {\delta V\over \delta \rho}\right), \label{grad2}
\end{eqnarray}
and when $c \rightarrow\infty$, we may formally get the geodesic flow on the tangent bundle over the Wasserstein space $P_2^\infty(M, \mu)$, i.e., 
\begin{eqnarray}
\partial_t \rho-\nabla_\mu^*(\rho\nabla \phi)&=&0,\label{Geo-TE-1}\\
\partial_t \phi+{1\over 2}|\nabla\phi|^2&=&0.\label{Geo-HJ-1}
\end{eqnarray}
Indeed, the convergence in a precise sense will be rigorously proved. See Theorem~\ref{main.convergence} for the precise statement and see Section 7 for the proof.

\subsection{The compressible Euler equation with damping}

When $\nu=0$,  the compressible Euler equations and the deformed  Hamilton-Jacobi equation with the transport equation have been well studied in the literature at least in the Euclidean case. See e.g. Carles \cite{Car}  and reference therein. The case of compact Riemannian manifolds is as the same as in the Euclidean case.  On the other hand, the compressible Euler equation with damping on Euclidean space has been also well studied in the literature. See   Wang and Yang \cite{WY}. See also Sideris, Thomases ad Wang \cite{STW} and reference therein. In this subsection, we develop a little bit more detail on  the link between the deformed Hamilton-Jacobi equation and the compressible Euler equation with damping on compact Riemannian manifolds.

Let $u=\nabla \phi$, $\gamma={1\over c^2}$, and $p\in C^1(\mathbb{R})$ be such that $\nabla p(\rho)=\rho {\delta V\over \delta \rho}(\rho)$,  then the Langevin deformation of flows~\eqref{ww}, \eqref{vvvv} turns into
\begin{eqnarray}
\partial_t \rho-\nabla_\mu^* (\rho u) &=&0, \label{tran}\\
\partial_t u+u\cdot \nabla u &=&-\gamma u-{1\over c^2} {\nabla p(\rho)\over \rho}. \label{euler}
\end{eqnarray}

In the case  $M=\mathbb{R}^n$ and $\mu=dx$, the above system is the compressible Euler equation with damping, where $\rho$ is the density of fluid, $u$ is the velocity of the fluid, and $\gamma={1\over c^2}$  is the friction constant, $p(\rho)$ is the pressure of the fluid. 
% $ p(\rho) = \rho^{\gamma}$ with $\gamma >1$ refers to the isentropic fluid, while $\gamma = 1$ is the isothermal fluid. 
When $\gamma=0$,  the connections between the compressible Euler equations and the deformed  Hamilton-Jacobi equation with the transport equation have been well studied in the literature, at least in the Euclidean case. See e.g. Carles \cite{Car}  and reference therein. The case of compact Riemannian manifolds is the same as in the Euclidean case.  On the other hand, there are extensive studies on the compressible Euler equation with damping on Euclidean space. See Wang and Yang \cite{WY}, and also Sideris, Thomases and Wang \cite{STW} and reference therein. 

In the case $M$ is a compact Riemannian manifold, we can regard $(\ref{tran})$ and $(\ref{euler})$ as the compressible Euler equation with damping on the compact Riemannian manifold $(M, g)$ equipped with the reference measure $\mu$. 
In this case, the compressible Euler equation with damping can be rewritten as follows
\begin{eqnarray}
\partial_t \rho-\nabla_\mu^* (\rho u) &=&0, \label{tran1}\\
\partial_t u+\nabla_u u &=&-\gamma u-{1\over c^2}\nabla V'(\rho), \label{euler1}
\end{eqnarray}
where $\nabla_{u} u$ denotes the Levi-Civita covariant derivative of the vector field $u$ along the  velocity field $u$ of the trajectory of the fluid.

 \subsection{Well-posedness of the Langevin deformation of flows}

In this subsection, we prove the existence and uniqueness of  the Langevin deformation of flows \eqref{ww} and \eqref{vvvv} for any fixed $c \in (0, \infty)$. 

Recall the Sobolev inequalities on compact Riemannian manifolds. We  only consider the unweighted case $\mu=\nu$. The general case can be treated similarly. By \cite{Au}, there exists a constant $C_{\rm Sob}>0$ such that
\begin{eqnarray}
\|f\|_{2n\over n-2}\leq C_{\rm Sob}(\|\nabla f\|_2+\|f\|_2), \ \ \ \forall f\in C^\infty(M).
\end{eqnarray}
Moreover,  for any $\alpha\in (0, 1)$, if $k>\alpha+{n\over 2}$, the Kondrakov embedding theorem holds
\begin{eqnarray}
\|f\|_{C^{0, \alpha}}\leq C_\alpha \|f\|_{k, 2}.
\end{eqnarray}
In particular, we have 
\begin{eqnarray}
\|f\|_\infty\leq C_\alpha \|f\|_{k, 2}.
\end{eqnarray}

Let $H^s(M)$ denotes the Sobolev space equipped with the Sobolev norm 
 \begin{eqnarray*}
 \|f\|_{s, 2}=\left(\sum\limits_{|\alpha|\leq s}\|D^\alpha f\|_{L^2}^2\right)^{1/2}
 \end{eqnarray*}
 where $\alpha$ is a multi-index, and for any vector valued function $U=(u_1, \ldots, u_n)$ on $M$, 
 \begin{eqnarray*}
 \|U\|_{s, 2}=\sum\limits_{i=1}^n \|u_i\|_{s, 2}
 \end{eqnarray*}
 For any $s\in \mathbb{Z}$, we define the Banach space
 \begin{eqnarray*}
 X_s=\{(f, g)| f\in H^2(M), g\in \oplus^n H^s(M)\}
 \end{eqnarray*}
equipped with the norm $\|(f, g)\|_{s, 2}=\|f\|_{s, 2}+\|g\|_{s, 2}$.

The following result gives the local existence and uniqueness of the solution to the compressible Euler equation with damping on compact Riemannian manifolds.

\bthm\label{thmm1}  (Local existence and uniqueness of smooth solution) Let $M$ be $\mathbb{R}^n$ or a compact Riemannian manifold, $s=[{n\over 2}]+1$. Let $V(\rho)=\int_M \rho \log \rho d\mu$ or $V(\rho)={1\over m-1}\int_M \rho^m d\mu$ for $m>1$. Suppose that $(\rho_0, u_0)\in H^{s+1}(M)$  with $\rho_0>0$. Then, there exists a constant $T>0$ such that the Cauchy problem of the compressible Euler  with damping $(\ref{tran})$ and $(\ref{euler})$ has a unique smooth solution $(\rho, u)$ in $C([0, T], H^{s+1}(M))\times  H^{s}(M))$. 
\ethm
\bpf
This is a standard result which can be proved by the method in Kato \cite{Ka} and Majda \cite{Maj}.
\epf

The following result gives the global existence and uniqueness of the solution to the compressible Euler equation with damping on compact Riemannian manifolds. 

\bthm \label{thmm2} (Global existence and uniqueness of smooth solution with small initial data) Let $M$ be $\mathbb{R}^n$ or a compact Riemannian manifold.  
Let $V(\rho)=\int_M \rho \log \rho d\mu$ or $V(\rho)={1\over m-1}\int_M \rho^m d\mu$ for $m>1$. Let $s=[{n\over 2}]+1$, $l\geq 2$. Then there  exists $\delta_0>0$ such that if $(\rho_0-1, u_0)\in H^{s+l}(M)$ is a small smooth initial 
 value in the sense that 
$\|\rho_0-1\|_{s+l, 2}+\|u_0\|_{s+l, 2}\leq \delta_0$ is sufficiently small.  Then the Cauchy problem of the compressible Euler equation with damping $(\ref{tran})$ and $(\ref{euler})$ admits a unique global smooth solution $(\rho, u)\in 
 C([0, \infty), H^{s+l}(M)\times H^{s+l-1}(M))$ with initial value $(\rho_0, u_0)$ and satisfying the following  energy estimate 
 \begin{eqnarray}
 & &\|\partial_t\rho(t)\|^2_{s+l-1, 2}+\|\rho(t)-1\|^2_{s+l, 2}+\|u(t)\|^2_{s+l, 2}\nonumber\\
 & &\hskip1.5cm +\int_0^t (\|\partial_t\rho(r)\|_{s+l-1, 2}+\|\nabla\rho(r)\|^2_{s+l-1, 2}+\|u(s)\|^2_{s+l, 2})dr\nonumber\\
 & &\hskip1cm \leq C(\|\partial_t\rho(0)\|_{s+l-1, 2}+\|\rho_0-1\|_{s+l, 2}+\|u_0\|_{s+l, 2}).
 \end{eqnarray}
\ethm
\bpf
 In the case $M=\mathbb{R}^n$, $V(\rho)=\int_M \rho \log \rho d\mu$ or $V(\rho)={1\over m-1}\int_M \rho^m d\mu$ with $m>1$,  this is the well-established result due to Wang and Yang \cite{WY}. See also Sideris, Thomases ad Wang \cite{STW} for the case $M=\mathbb{R}^n$ and $V(\rho)={1\over m-1}\int_M \rho^m d\mu$ with $m>1$. In the case $M$ is a compact Riemannian manifold and $V(\rho)={1\over m-1}\int_M \rho^m d\mu$ with $m>1$, the proof of theorem is similar to the ones in \cite{WY, STW}. In the case $M$ is a compact Riemannian manifold and $V(\rho)=\int_M \rho\log \rho d\mu$, we can modify the proof of the main results in  \cite{WY, STW}. The main point here is that on compact Riemannian manifold, the positivity of the initial data $\rho_0>0$ implies that there exists a constant $\varepsilon_0>0$ such that $\rho_0\geq \varepsilon_0>0$, and the argument used in Wang and Yang \cite{WY} can be extended to the case $V(\rho)=\int_M \rho\log \rho d\mu$ on compact Riemannian manifolds. We can also modify the argument used in 
Sideris, Thomases ad Wang \cite{STW} by taking the sound speed to be $\sigma(\rho)=\log \rho$. To save the length of the paper, we omit the detail of the proof. 
\epf

Now we turn back to the Langevin deformation of flows~\eqref{ww} and \eqref{vvvv}. The key point here is that we need to prove that if the initial value $u_{0} = \nabla \phi_{0}$ for some function $\phi_{0}$, $u_{t}$ will keep the gradient form along $t > 0$.  To see this, we show the following result. 
\bthm\label{potential} Let $M=\mathbb{R}^n$ or a compact Riemannian manicold, $(\rho, u)$ be a smooth solution to the compressible Euler equation with damping, i.e., $(\ref{tran})$ and $(\ref{euler})$. Let $\omega=d u$. Suppose that ${\nabla p(\rho)\over \rho}=\nabla V'(\rho)$. Then 
\begin{eqnarray}
\partial_t \omega+u\cdot \nabla \omega+\omega \wedge \nabla u^*=-\gamma \omega. \label{omega-1} 
\end{eqnarray} Moreover, if $|\nabla u|_{L^\infty}\leq C$, then for all $t\in [0, T]$, we have \footnote{In the case $C=\max\limits |\nabla u|_{L^\infty}<\gamma$, the above estimates holds for all $t\in [0, \infty)$.}
\begin{eqnarray}
\|\omega(t)\|_{L^p}\leq \|\omega(0)\|_{L^p}e^{(C-\gamma)t}. \label{omega-2}
\end{eqnarray}
In particular, if $u_0$ is a closed form, so is $u(t, \cdot)$, i.e., $du_0=0$ implies $du(t, \cdot)=0$ on $[0, T]$.
\ethm

\bpf
We first prove that $\omega=d u$ satisfies \eqref{omega-1}. By identifying $u$ with its dual $u^*$, and identifying $\nabla V'(\rho)$ with its dual $d V'(\rho)$,  with respect to the Riemannian metric on $M$, we can rewrite $(\ref{euler1})$ (or \eqref{euler}) as follows
\begin{eqnarray}\label{euler2}
\partial_t u^*+\nabla_u u^*&=&-\gamma u^*-{1\over c^2} dV'(\rho). 
\end{eqnarray}

Taking exterior differentiation on the both sides of the compressible Euler equation $(\ref{euler2})$, letting $\omega=d u^*=\sum\limits_{i=1}^n du_i\wedge e_i^*\in \Gamma(\Lambda^2 TM)$, and using $ddV'(\rho)=0$, we have
\begin{eqnarray}\label{euler3}
\partial_t \omega+ d \nabla_{u} u^* =-\gamma \omega.
\end{eqnarray}
%
%{\blue
%\begin{eqnarray}
%\partial_t \omega+ \nabla_{u} \omega+\sum\limits_{i=1}^n (\omega (e_i))\wedge  \nabla_{e_i} u^*+\sum\limits_{k=1}^n e_k^*\wedge R(e_k, u)u^*=-\gamma \omega.
%\end{eqnarray}
%}
%
%
Let $(e_i)$ be an ONB and normal at $x\in M$, writing  $u=\sum\limits_{i=1}^n u_ie_i$ and $u^*=\sum\limits_{i=1}^n u_ie_i^*$, we derive that 
\begin{eqnarray}\label{duu}
d \nabla_{u} u^*=\nabla_{u} (du^*)+\sum\limits_{i=1}^n du_i \wedge  \nabla_{e_i} u^*+\sum\limits_{k=1}^n e_k^* \wedge R(e_k, u) u^*. 
\end{eqnarray}
Indeed, since $(e_i)$ is an ONB and normal at $x\in M$, we have
\begin{eqnarray*}
d\nabla_{u} u^* =\sum\limits_{k=1}^n e_k^* \wedge \nabla_{e_k}(\nabla_u u^*)
%&=&\sum\limits_{k=1}^n e_k^* \wedge \nabla_{e_k} (\sum\limits_{i=1}^n u_i\nabla_{e_i}u^*)\\
%&=&\sum\limits_{k, i=1}^n  \left(e_k^*\wedge \nabla_{e_k} u_i  \nabla_{e_i}u^*+u_i e_k^*\wedge \nabla_{e_k} \nabla_{e_i}u^*\right)\\
= \sum\limits_{i=1}^n du_i \wedge \nabla_{e_i}u^*+\sum\limits_{k, i=1}^n u_i e_k^*\wedge \nabla_{e_k} \nabla_{e_i}u^*,
\end{eqnarray*}
and 
\begin{eqnarray*}
\nabla_{u} du^* =\sum\limits_{i=1}^n u_i\nabla_{e_i} (du^*)
%&=&\sum\limits_{i=1}^n u_i \nabla_{e_i} (\sum\limits_{k=1}^n e_k^*\wedge \nabla_{e_k}u^*)\\
%&=&\sum\limits_{k, i=1}^n \left(u_i \nabla_{e_i}e_k^*\wedge \nabla_{e_k}u^*+u_i e_k^*\wedge \nabla_{e_i}  \nabla_{e_k} u^*\right)\\
=\sum\limits_{k, i=1}^n u_i e_k^*\wedge \nabla_{e_i}  \nabla_{e_k} u^*.
\end{eqnarray*}
Hence
\begin{eqnarray*}
d\nabla_{u} u^* -\nabla_{u}(d u^*)
%&=&\sum\limits_{i=1}^n du_i \wedge \nabla_{e_i}u^*+\sum\limits_{k, i=1}^n u_i e_k^*\wedge \nabla_{e_k} \nabla_{e_i}u^*- \sum\limits_{k, i=1}^n u_i e_k^*\wedge \nabla_{e_i}  \nabla_{e_k} u^*\\
&=&\sum\limits_{i=1}^n du_i \wedge \nabla_{e_i}u^*+\sum\limits_{k, i=1}^n u_i e_k^*\wedge R(e_k, e_i)u^*.
\end{eqnarray*}
Moreover, we claim that for all $u\in \Gamma(TM)$,
\beqna\label{extra.term}
\sum\limits_{k=1}^n e_k^*\wedge R(e_k, u)u^*=0. 
\eeqna
Indeed, for any $i, j=1, \ldots, n$, acting on  $(e_i,  e_j)$, we have 
\begin{eqnarray*}
\sum^{n}_{k=1}(e^{*}_{k}\wedge R(e_k, u)u^{*})(e_i, e_j)
%&=&\sum^{n}_{k=1} e^{*}_{k}(e_i)(R(e_k, u)u^{*})(e_j) - e^{*}_{k}(e_j)(R(e_k, u)u^{*})(e_i)\\
= (R(e_i, u)u^{*})(e_j) - (R(e_j, u)u^{*})(e_i).
\end{eqnarray*}
%Note that, for any one-form $\alpha$ and vector fields $X$ and $Y$,  $\nabla_X \alpha (Y)=(\nabla_X \alpha )(Y)-\alpha(\nabla_X Y)$. 
%Taking $\alpha=\nabla_u u^*$, we have
Notice that
%\begin{eqnarray*}
%\nabla_{e_i}((\nabla_u u^*)(e_j)) =(\nabla_{e_i}\nabla_u u^*)(e_j)-(\nabla_u u^*)(\nabla_{e_i}e_j) =(\nabla_{e_i}\nabla_u u^*)(e_j),
%\end{eqnarray*}
%and
%\begin{eqnarray*}
%\nabla_{e_i} ((\nabla_u u^*) (e_j))=\nabla_{e_i}(\nabla_u u^*(e_j)+u^*(\nabla_u e_j) )=\nabla_{e_i}\nabla_u u_j+\nabla_{e_i}(u^*(\nabla_u e_j)),
%\end{eqnarray*}
%which implies that   
\begin{eqnarray*}
(\nabla_{e_i}\nabla_u u^*)(e_j)=\nabla_{e_i}\nabla_u u_j+\nabla_{e_i}(u^*(\nabla_u e_j)),
\end{eqnarray*}
meanwhile we also have
\begin{eqnarray*}
(\nabla_u\nabla_{e_i}u^*) (e_j)
%&=&\nabla_u( (\nabla_{e_i}u^*) (e_j))+(\nabla_{e_j} u^*)(\nabla_u e_j)\\
%&=&\nabla_u (\nabla_{e_i}(u^* (e_j))+u^*(\nabla_{e_i}e_j))\\
= \nabla_u \nabla_{e_i}u_j+\nabla_u (u^*(\nabla_{e_i}e_j)),
\end{eqnarray*}
such that
\begin{eqnarray*}
 (R(e_i, u)u^{*})(e_j)&=&\nabla_{e_i}(u^*(\nabla_u e_j))-\nabla_u (u^*(\nabla_{e_i}e_j))\\
 %&=&\nabla_{e_i}\langle u, \nabla_u e_j\rangle -\nabla_u \langle u, \nabla_{e_i}e_j\rangle\\
 %&=&\langle \nabla_{e_i} u, \nabla_u e_j\rangle+\langle u, \nabla_{e_i}\nabla_u e_j\rangle-\langle \nabla_u u, \nabla_{e_i}e_j \rangle-\langle u, \nabla_u \nabla_{e_i}e_j\rangle\\
% &=&\langle u, \nabla_{e_i}\nabla_u e_j\rangle-\langle u, \nabla_u \nabla_{e_i}e_j\rangle\\
 &=&\langle u, R(e_i, u)e_j\rangle.
 \end{eqnarray*}
Exchanging $e_i$ and $e_j$ leads to 
\begin{eqnarray*}
 (R(e_j, u)u^{*})(e_i)=\langle u, R(e_j, u)e_i\rangle.
 \end{eqnarray*} 
Therefore we derive that
\begin{eqnarray*}
\sum^{n}_{k=1}(e^{*}_{k}\wedge R(e_k, u)u^{*})(e_i, e_j) = R(e_i, u, e_j, u)-R(e_j, u, e_i, u) = 0,
\end{eqnarray*}
which finishes the proof of \eqref{extra.term}. 

Then taking \eqref{duu} and \eqref{extra.term} into \eqref{euler3}, we obtain 
\begin{eqnarray*}
\partial_t \omega+ \nabla_{u} \omega+\sum\limits_{i=1}^n (\omega (e_i))\wedge  \nabla_{e_i} u^*=-\gamma \omega,
\end{eqnarray*}
which is written briefly as \eqref{omega-1}.
%\begin{eqnarray*}
%(R(e_i, u)u^{*})(e_j) - (R(e_j, u)u^{*})(e_i)
%%&=&\langle u, R(e_i, u)e_j\rangle-\langle u, R(e_j, u)e_i\rangle\\
%&=&R(e_i, u, e_j, u)-R(e_j, u, e_i, u), 
%\end{eqnarray*}
%%Note that  $R(X, Y, Z, W) = R(Z, W, X, Y)$. Thus
%i.e.
%\begin{eqnarray*}
%R(e_i, u, e_j, u)=R(e_i, u, e_i, u).
%\end{eqnarray*}
%This yields 
%\begin{eqnarray*}
%\sum^{n}_{k=1}(e^{*}_{k}\wedge R(e_k, u)u^{*})(e_i, e_j) = 0,
%\end{eqnarray*}
%which finishes our proof. 

%Thus, Eq. $(\ref{euler3})$ reads as follows
%which will be written briefly as follows
%\begin{eqnarray}
%\partial_t \omega+ \nabla_{u} \omega+\omega \wedge  \nabla u^*=-\gamma \omega.\label{euler4} 
%\end{eqnarray}

Now taking inner product with $|\omega|^{p-2}\omega$ in the both sides of $(\ref{omega-1})$, and integrating on $M$, we have
\begin{eqnarray*}
\int_M \langle {D\over \partial t}\omega, |\omega|^{p-2}\omega\rangle d\mu +\int_M\langle \omega \wedge \nabla u, |\omega|^{p-2}\omega\rangle d\mu=-\gamma \|\omega\|_p^p
\end{eqnarray*}
where $ {D\over \partial t}\omega=\partial_t \omega+u\cdot\nabla \omega$. Note that 
\begin{eqnarray*}
 \langle {D\over \partial t}\omega, |\omega|^{p-2}\omega\rangle={1\over p}{D\over \partial t}|\omega|^p.
\end{eqnarray*}
Hence
\begin{eqnarray*}
 {1\over p}{d\over dt}\|\omega(t)\|_p^p&=& -\int_M \left(\langle \omega \wedge \nabla u, |\omega|^{p-2}\omega\rangle d\mu-\gamma \|\omega(t)\|_p^p\right)d\mu\\
 &\leq& \int_M\left( |\nabla u| |\omega(t)|^p d\mu-\gamma \|\omega(t)\|_p^p\right)d\mu\\
 &\leq& \left(\|\nabla u\|_{L^\infty}-\gamma\right)\|\omega(t)\|_p^p.
\end{eqnarray*}
By Gronwall inequality, we have 
\begin{eqnarray*}
\|\omega(t)\|_p^p\leq \|\omega(0)\|_p^p\exp\left\{p\left(\|\nabla u\|_{L^\infty}-\gamma\right)t\right\}.
\end{eqnarray*}
Thus, if $\omega(0)=0$, then for all $t\geq 0$, we have $\omega(t)=0$.  Hence, if $u^*(0, \cdot)$ is a closed one-form on $M$, then $u^*(t, \cdot)$ is also a closed one-form on $M$. 

\epf

We now state the main result of this section, which gives the existence and uniqueness of the local solution to Langevin deformation of flows~\eqref{ww} and \eqref{vvvv} for any fixed $c \in (0, \infty)$.

\bthm \label{existence} Let $M=\mathbb{R}^n$ or be a compact Riemannian manifold, $c\in [0, \infty]$  Given $(\rho_0, \phi_0)\in TP_2^\infty(M, \mu)$ with  $\rho_0, \phi_0\in C^\infty(M)$, there exists $T=T_c>0$ such that the Cauchy problem of the Langevin deformation of flows~\eqref{ww} and \eqref{vvvv} has a unique solution  $(\rho, \phi)\in  C^1([0, T], C^\infty(M)^2)$. 
 \ethm
 \bpf
The cases $c=0$ and $c=\infty$ are well known.  For $c\in (0, \infty)$, consider the compressible Euler equation with damping on $M$
\begin{eqnarray}\label{EEEEE}
& &\partial_t u+u\cdot\nabla u=-{u\over c^2} - {1\over c^2}{\delta V\over \delta \rho},\ \ \ \  \left. u\right|_{t=0}=\nabla \phi_0,\\
& &\partial_t \rho+\nabla_\mu^*(\rho u)=0, \ \ \  \ \left.\rho\right|_{t=0}=\rho_0.
\end{eqnarray}
By Theorem \ref{thmm1},  if the initial dada are in $H^s(M, \mu)$ for any $s>{n\over 2}+1$, then 
there exists $T=T_c>0$ such that above system has a  unique solution $(\rho, u)\in C([0, T], H^s(M, \mu))^2$. Moreover, tame estimates show that the time of existence $T>0$ can be chosen independent of $s>{n\over 2}+1$. 

By Theorem \ref{potential},  if $u^*_0=d\phi_0$,  $u^*(t, \cdot)$ is closed on $M$. For all $(t, x)\in [0, T]\times M$, let 
\begin{eqnarray}
\label{HJHJHJ} \phi(t, x)=e^{-\gamma t}\phi_0(x)+e^{-\gamma t}\int_0^t e^{\gamma s}\left(f(\rho(s, x))-{1\over 2 }|u(s, x)|^2 \right)ds, 
 \end{eqnarray}
 where $\gamma={1\over c^2}$ and $f(\rho)= - {1\over c^2}{\delta V\over \delta \rho}$. We have
\begin{eqnarray*}
\partial_t  \phi=-\gamma \phi+f(\rho(t, x))-{1\over 2 }|u(t, x)|^2.
 \end{eqnarray*}
Note that, as $u^*(t, \cdot)$ is a closed one-form on $M$, it holds that
\begin{eqnarray*}
\nabla |u|^2=2u\cdot \nabla u.
\end{eqnarray*}
Hence we can check that
\begin{eqnarray*}
\partial_t (\nabla \phi-u)=-\gamma (\nabla \phi-u).
\end{eqnarray*}
Note that at $t=0$, $u(0)=\nabla\phi(0)$. Thus $u(t)=\nabla\phi(t)$ on $[0, T]\times M$.  Substituting this into the compressible Euler equation with damping, we have 
\begin{eqnarray*}
\nabla\left(\partial_t \phi+{1\over 2}|\nabla \phi|^2\right)=-\gamma \nabla \phi-{1\over c^2}\nabla V'(\rho).
\end{eqnarray*}
This indicates that  the Cauchy problem of \eqref{ww} and \eqref{vvvv} has a unique solution $(\rho, \phi)\in C([0, T], H^s(M)\times H^{s+1}(M))$ for all $s>{n\over 2}+1$. The proof is completed. 
 \epf

 \subsection{Entropy dissipation formulas along the deformation of flows}
 
In this subsection we prove the entropy dissipation formulas along the Langevin deformation of flows on the Wasserstein space.  In this subsection, we assume that  $M$ is Euclidean space or a compact Riemannian manifolds. By Theorem \ref{existence}, for any $c\in [0, \infty]$,  the Langevin deformation of flows on $TP_2(M, \mu)$ has a unique smooth solution (up to an additional constant) on $[0, T_c]\times P_2^\infty(M, \mu)$. 

Similarly to the case of finite dimensional manifolds, we have the following variational formulas for the Hamiltonian and Lagrangian along the Langevin deformation of  flows~\eqref{ww} and \eqref{vvvv}.
%
%\bprop \label{Th-H-c} Let $(\rho, \phi)$ be the solution to the Langevin deformation of flows~\eqref{ww} and \eqref{vvvv} on $TP_2^\infty(M, \mu)$. Define the Hamiltonian and the Lagrangian on $TP_2^\infty(M, \mu)$ by
%\begin{eqnarray}
%H(\rho, \phi)&=&{c^{2} \over 2}\int_M |\nabla \phi|^2\rho d\mu +V(\rho),\\
%L(\rho, \phi)&=&{c^{2} \over 2}\int_M |\nabla \phi|^2\rho d\mu -V(\rho).
%\end{eqnarray}
%Then we have
%\beqna
%{d\over dt} H(\rho, \phi) \leq 0.
%\eeqna
%and 
%
%\eprop
%\bpf 
%
%The second order derivative of $H(\rho, \phi)$ is 
%\beqnas
%{d^2\over dt^2} H(\rho, \phi) &=& - {d\over dt}\int_{M} |\nabla \phi|^{2}\rho d\mu\\
%&=&2\int_M \langle \nabla \phi, \nabla (\partial_t \phi + {1\over 2}|\nabla\phi|^2) \rangle \rho d\mu\\
%&=&{2\over c^2}\int_M \left\langle \nabla \phi, \nabla (\phi + {\delta V\over \delta \rho}) \right\rangle \rho d\mu.
%\eeqnas
%In our case, $V(\rho)=\int_M \rho \log\rho d\mu$, and ${\delta V\over \delta \rho}=\log\rho+1$. Hence
%\beqnas
%{d^2\over dt^2} H(\rho, \phi)={2\over c^2}\int_M \left\langle \nabla \phi, \nabla \phi +\nabla\log\rho \right\rangle \rho d\mu.
%\eeqnas
%\epf

\bthm\label{Th-LL-T1} Let $(\phi_t, \rho_t)$ be a smooth solution to Eq. $(\ref{ww})$ and Eq. $(\ref{vvv})$ on $TP_2^\infty(M, \mu)$. 
%\begin{eqnarray}
%\partial_t \rho-\nabla_\mu^*(\rho\nabla \phi)&=&0,\label{GF1}\\
%c^2\left({\partial\phi\over \partial t}+{1\over 2}|\nabla \phi|^2\right)&=&-\phi-{\delta V\over \delta \rho}.\label{GF2}
%\end{eqnarray}
%Then  $(\rho, \phi)$ is a deformation of flows on $TP_2^\infty(M, \mu)$ which interpolates the geodesic flow $(\rho, \phi)$,  which satisfies $(\ref{TA})$  together  with $(\ref{HJ})$,  and the backward gradient flow of $V$ on $P_2^\infty(M, \mu)$, which satisfies 
%\begin{eqnarray*}
%\partial_t \rho=\nabla_\mu^*\left(\rho\nabla{\delta V\over \delta \rho}\right),
%\end{eqnarray*}
Let
\begin{eqnarray*}
H(\rho, v)&=&{c^2\over 2}\int_M |\nabla\phi(x)|^2\rho d\mu+V(\rho),\\
L(\rho, v)&=&{c^2\over 2}\int_M |\nabla\phi(x)|^2\rho d\mu-V(\rho).
\end{eqnarray*}
Then
\begin{eqnarray*}
{d\over dt}H(\rho_t, \phi_t)=-\int_M |\nabla\phi_t|^2\rho_t d\mu,
\end{eqnarray*}
and
\begin{eqnarray*}
{d^2\over dt^2}L(\rho_t, \phi_t)&=&2c^{-2}\left\|\dot\rho+\nabla V(\rho)\right\|^2-2\nabla^2 V(\rho)\left(\dot\rho, \dot\rho\right)\\
&=&c^2\|\ddot{\rho}\|^2-2\nabla^2 V\left(\dot\rho, \dot\rho\right).
\end{eqnarray*}
In particular, if $-V$ is $K$-convex on $P_2^\infty(M, \mu)$, i.e., the Hessian of $-V$ on $P_2^\infty(M, \mu)$ satisfies 
$$-\nabla^2 V(\rho)=-{\rm Hess}_{P_2^\infty(M)}V(\rho)\geq K,$$
 then
\begin{eqnarray*}
{d^2\over dt^2}L(\rho_t, v_t)\geq 2K\int_M |\nabla\phi_t|^2\rho_t d\mu+2\int_M \left[\nabla\left(\partial_t \phi+{1\over 2}|\nabla\phi|^2\right)\right]^2\rho d\mu.
\end{eqnarray*}
\ethm
\bpf  Notice that this result can be directly derived from Proposition~\ref{Th-B-2}, for that on $TP_2^\infty(M, \mu)$, 
$$
{d\over dt} H(\rho, \phi) = -|v|^{2} =  - \int_{M} |\nabla \phi|^{2}\rho d\mu.
$$
Of course we could also prove it by direct computations. 
\beqnas
{d\over dt} H(\rho, \phi)=  {d\over dt}\left({c^{2} \over 2}\int_M |\nabla \phi|^2\rho d\mu \right) +  \langle \nabla V, \dot{\rho} \rangle.
\eeqnas
Note that
\beqnas
{d\over dt}\left({c^{2} \over 2}\int_M |\nabla \phi|^2\rho d\mu\right) 
&=& c^{2}\int_M \langle \nabla \phi, \nabla \partial_t \phi \rangle \rho d\mu + {c^{2} \over 2}\int_M |\nabla \phi|^2\partial_t \rho d\mu\\
&=& c^{2}\int_M \langle \nabla \phi, \nabla \partial_t \phi \rangle \rho d\mu + {c^{2} \over 2}\int_M\langle \nabla |\nabla \phi|^2, \nabla \phi \rangle \rho d\mu\\
&=& c^{2}\int_M \langle \nabla \phi, \nabla (\partial_t \phi + {1\over 2}|\nabla\phi|^2) \rangle \rho d\mu\\
&=&-\int_M \left\langle \nabla \phi, \nabla (\phi + {\delta V\over \delta \rho}) \right\rangle \rho d\mu.
\eeqnas
On the other hand
\beqnas
{d\over dt}V(\rho_t)=\langle \nabla V, \dot{\rho} \rangle=\int_M \langle { \delta V\over \delta \rho}, \nabla\phi\rangle \rho d\mu.
\eeqnas
Hence
\beqnas
{d\over dt} H(\rho, \phi) &=& -\int_M \langle \nabla \phi, \nabla (\phi + {\delta V\over \delta \rho}) \rangle \rho d\mu +  \langle \nabla V, \dot{\rho} \rangle\\
&=& - \int_{M} |\nabla \phi|^{2}\rho d\mu.
\eeqnas

The second dissipation formula in Theorem \ref{Th-LL-T1} can be proved by straightforwardly extending the argument used in the proof of Proposition \ref{Th-B-2} to the Wasserstein space $P_2^\infty(M, \mu)$ equipped with Otto's infinite dimensional Riemannian metric. To save the length of the paper, we omit the detail here. \epf

\subsection{Proof of Theorem \ref{MT0}} 
Applying Theorem \ref{Th-LL-T1} to  $V(\rho)=-{\rm Ent}(\rho)=-
\int_{M}\rho\log \rho d\mu$, we can derive Theorem \ref{MT0}.

%
%\bthm \label{Th-LL-T2} Let $c>0$. Let $\rho, \phi$ be a smooth solution of the following equations
%\begin{eqnarray*}
%\partial_t \rho+\nabla_\mu^*(\rho\nabla \phi)&=&0,\\
%c^2\left({\partial\phi\over \partial t}+{1\over 2}|\nabla \phi|^2\right)&=&-\phi+\log \rho+1.
%\end{eqnarray*}
%Let
%\begin{eqnarray*}
%H(\rho, \phi)={c^2\over 2}\int_M |\nabla\phi|^2\rho d\mu+{\rm Ent}(\rho).
%\end{eqnarray*}
%Then
%\begin{eqnarray*}
%{d\over dt}H(\rho)&=&2\int_M \nabla\phi\cdot\nabla \rho d\mu-\int_M |\nabla\phi|^2\rho d\mu\\
%&=&-\int_M (2L\phi+|\nabla \phi|^2) \rho d\mu,
%\end{eqnarray*} 
%and
%\begin{eqnarray}
%{d^2\over dt^2}H(\rho, \phi)=2\int_M \left[c^{-2}|\nabla \phi-\nabla\log\rho|^2+|{\rm Hess}\phi|^2+Ric(L)(\nabla\phi, \nabla\phi)\right]\rho d\mu.\label{2ndH}
%\end{eqnarray}
%\ethm
\bpf  Indeed, applying Theorem \ref{Th-LL-T1} to $V=-{\rm Ent}$, we have
\begin{eqnarray*}
{d^2\over dt^2}L(\rho, \phi)=2c^2\|\ddot \rho\|^2+2{\rm Hess}_{P_2^\infty(M)} {\rm Ent}(\rho)\left(\dot\rho, \dot\rho\right).
\end{eqnarray*}
Here ${\rm Hess}_{P_2^\infty(M)}{\rm Ent}$ is the Hessian of ${\rm Ent}$ on the Wasserstein space ${P}_2^\infty(M)$ equipped with Otto's infinite dimensional Riemannian metric. By Theorem \ref{Ent-Hess} , we have
\begin{eqnarray*}
{\rm Hess}_{P_2^\infty(M)}{\rm Ent}(\rho)(\dot\rho, \dot\rho)=\int_M (|{\rm Hess}\phi|^2+\Ric(L)(\nabla\phi, \nabla\phi))\rho d\mu.
\end{eqnarray*}
On the other hand, we have
\begin{eqnarray*}
c^2\|\ddot\rho\|^2&=&c^2\int_M \left|\nabla \left({\partial\phi\over \partial t}+{1\over 2}|\nabla \phi|^2\right)\right|^2\rho d\mu\\
&=&c^{-2}\int_M |\nabla \phi+\nabla\log\rho|^2\rho d\mu.
\end{eqnarray*}
This finishes the proof of Theorem \ref{MT0}. \epf

On compact manifolds with non-negative Bakry-Emery Ricci curvature, Theorem \ref{MT0} implies that the Hamiltonian function $H$ is always monotonically nonincreasing and 
the Lagrangian function $L$ is always convex along the Langevin deformation flow  $(\rho, \phi)$ which interpolates the geodesic flow and the backward gradient flow of $V$ on the Wasserstein space over a compact Riemannian manifold with weighted measure.

Note that, when $c=0$ in Theorem \ref{MT0}, we have $\phi=-\log\rho-1$, hence $\nabla\phi+\rho^{-1}\nabla\rho=0$. This yields $\partial_t \rho=L\rho$,  and 
\begin{eqnarray*}
{d^2\over dt^2}{\rm Ent}(\rho(t))= 2\int_M \left[|{\rm Hess}\phi|^2+\Ric(L)(\nabla\phi, \nabla\phi)\right]\rho d\mu.
\end{eqnarray*}
If we formally take $c=\infty$ in Theorem \ref{MT0}, we have $\ddot\rho=0$, and we obtain 
\begin{eqnarray*}
{d^2\over dt^2}H(\rho, \phi)= 2\int_M \left[|{\rm Hess}\phi|^2+\Ric(L)(\nabla\phi, \nabla\phi)\right]\rho d\mu.\label{2ndHH}
\end{eqnarray*}
However,  this formula is not correct. Indeed,  when $c=\infty$, the kinetic energy term ${1\over 2}\int_M |\nabla\phi(t)|^2\rho(t) d\mu$ is a constant along the geodesic flow $(\rho(t), \phi(t))$ on the Wasserstein space $P_2(M, \mu)$,  thus 
${c^2\over 2}\int_M |\nabla\phi(t)|^2\rho(t) d\mu=\infty$. In this case, we must replace the left hand side  of $(\ref{2ndH})$ in Theorem  \ref{MT0} by 
the second order derivative of ${\rm Ent}(\rho(t))$,  which is given by the entropy dissipation formula  (Theorem  \ref{Ent-Hess}). 
\begin{eqnarray*}
{d^2\over dt^2} {\rm Ent}(\rho(t))=\int_M [|{\rm Hess}\phi|^2+\Ric(L)(\nabla\phi, \nabla\phi)]\rho d\mu.
\end{eqnarray*}

In other words, we have the following

\bcor\label{Th-LL-C1} Let $M$ be a compact Riemannian manifold, $f\in C^2(M)$. Then \\
(i) When $c=\infty$, we have $\ddot\rho=0$, i.e., $(\rho, \phi)$ is a geodesic flow on $P_2(M, \mu)$, and satisfies the transport equation $(\ref{Geo-TE-1})$ and
the Hamilton-Jacobi equation $(\ref{Geo-HJ-1})$. Moreover
\begin{eqnarray*}
{d^2\over dt^2}{\rm Ent}(\rho(t))= \int_M \left[|{\rm Hess}\phi|^2+\Ric(L)(\nabla\phi, \nabla\phi)\right]\rho d\mu,
\end{eqnarray*}
(ii) When $c=0$, we have $\phi=-\log\rho-1$, i.e., $\rho$ is a positive solution to the heat equation
$$\partial_t \rho=L \rho.$$
Moreover
\begin{eqnarray*}
{d^2\over dt^2}{\rm Ent}((\rho(t))= 2\int_M \left[|{\rm Hess}\phi|^2+\Ric(L)(\nabla\phi, \nabla\phi)\right]\rho d\mu.
\end{eqnarray*}
\ecor

\brmq In \cite{LL-flow17}, we introduced the Langevin deformation of flows as follows  
\begin{eqnarray}
\partial_t \rho-\nabla_\mu^*(\rho\nabla \phi)&=&0, \label{ww-17} \\
c^2\left({\partial\phi\over \partial t}+{1\over 2}|\nabla \phi|^2\right)&=&-\phi + {\delta V\over \delta \rho}.\label{vvvv-17}
\end{eqnarray}
It is indeed the Langevin deformation between the backward gradient flow $\partial_t \rho=\nabla_\mu^*\left(\rho\nabla {\delta V\over \delta \rho}\right)$ on the Wasserstein space and the geodesic flow on the tangent bundle $TP_2^\infty(M, \mu)$.  Define the Hamiltonian and the Lagrangian on $TP_2^\infty(M, \mu)$ by
\begin{eqnarray*}
H(\rho, \phi)&=&{c^{2} \over 2}\int_M |\nabla \phi|^2\rho d\mu +V(\rho),\\
L(\rho, \phi)&=&{c^{2} \over 2}\int_M |\nabla \phi|^2\rho d\mu -V(\rho).
\end{eqnarray*}
Then 
\beqna
{d\over dt} L(\rho, \phi) =-\int_M |\nabla\phi|^2\rho d\mu,
\eeqna
and
\begin{eqnarray*}
{d^2\over dt^2}H(\rho, \phi)&=&2c^{-2}\left\|\dot\rho-\nabla V(\rho)\right\|^2+2\nabla^2 V(\rho)\left(\dot\rho, \dot\rho\right)\\
&=&c^2\|\ddot{\rho}\|^2+2\nabla^2 V\left(\dot\rho, \dot\rho\right).
\end{eqnarray*}
In particular, if $V$ is $K$-convex on $P_2^\infty(M, \mu)$, i.e., the Hessian of $V$ on $P_2^\infty(M, \mu)$ satisfies 
$$\nabla^2 V(\rho)={\rm Hess}_{P_2^\infty(M)}V(\rho)\geq K,$$
 then
\begin{eqnarray*}
{d^2\over dt^2}H(\rho, \phi)\geq 2K\int_M |\nabla\phi|^2\rho d\mu+2\int_M \left[\nabla\left(\partial_t \phi+{1\over 2}|\nabla\phi|^2\right)\right]^2\rho d\mu.
\end{eqnarray*}
In particular, taking $V(\rho)={\rm Ent}(\rho)=
\int_{M}\rho\log \rho d\mu$, we can derive that
\begin{eqnarray}
{d\over dt}L(\rho, \phi)&=&-\int_M |\nabla \phi|^2\rho d\mu,\label{monotoH}\\
{d^2\over dt^2}H(\rho, \phi)&=&2\int_M \left[c^{-2}|\nabla \phi-\nabla\log\rho|^2+|{\rm Hess}\phi|^2+\Ric(L)(\nabla\phi, \nabla\phi)\right]\rho d\mu.\label{2ndH}
\end{eqnarray}
Suppose that $Ric(L)=\Ric+\nabla^2 f\geq K$, then we have
\begin{eqnarray*}
{d^2\over dt^2}H(\rho, \phi)\geq 2K\int_M |\nabla\phi|^2\rho_t d\mu+{2\over n}\int_M |\Delta \phi|^2\rho d\mu+2c^{-2} \int_M |\nabla \phi-\nabla\log\rho|^2\rho d\mu.
\end{eqnarray*}
\ermq

\subsection{Convergence results}
In this subsection we show the convergence results for the Langevin deformation of flows. Note that the $L^2$-adjoint of $\nabla$ with respect to the standard volume measure on $\mathbb{R}^d$ or a compact Riemannian manifold is $\nabla^*=-\nabla\cdot$. Taking $V(\rho) = \int \rho \log \rho$ in $(\ref{ww})$ and $(\ref{vvvv})$,  and using ${\delta V\over \delta \rho} = \log \rho + 1$, we get the Langevin deformation of flows as follows 
\beqna
\label{langevin.eq1}\partial_t \rho + \nabla \cdot (\rho \nabla \phi) &=& 0, \\
\label{langevin.eq2}c^{2}(\partial_t \phi + \frac{1}{2}|\nabla \phi|^{2}) &=& -\phi - \log \rho -1,
\eeqna
on $\bbR^{d}$ or compact manifold $M$, with initial value
$$
\rho(0, x) = \rho_{0}(x), \ \ \phi(0, x) = \phi_{0}(x).
$$
Our aim is to prove that when $c \rightarrow 0$ and $c \rightarrow \infty$, the solution $(\rho, \phi)$ converges in a precise sense to that to the heat equation and to the geodesic flow respectively.

The key point of the proof still relies on the close connection between the Langevin deformation of flows and the compressible Euler equation with damping: we first prove the convergence for $(\ref{tran})$ and $(\ref{euler})$ with $V(\rho) = \int \rho \log \rho$, and then turn back to $(\ref{langevin.eq1})$ and $(\ref{langevin.eq2})$. We will first show the convergence in Euclidean space and then extend the results to compact manifolds. For simplicity, we consider the Laplace-Beltrami operator here instead of the Witten Laplacian. Under suitable assumptions on the potential function $f$, it is easy to prove that the convergence results also hold for the Witten Laplacian. 
 
 We now state the convergence theorem in a precise sense. To keep the paper to be more readable, we leave its  proof to the Section 7. 

\bthm\label{main.convergence}
Let $M$ be $\bbR^{n}$ or a compact Riemannian manifold. Let $s \in \bbN$ with $s > {d \over 2} +2$ and $c \in (0, \infty)$, $\bar{\rho} > 0$ be two constants. Let $(\rho^{c}, \phi^{c})$ be the local solution to the initial value problem of Langevin deformation of flows $(\ref{langevin.eq1})$, $(\ref{langevin.eq2})$. 
\bitem
{\item
Given the initial value $(\rho_{0}, \phi_{0}) \in TP^{\infty}_{2}(M, \mu)$ satisfying $\rho_{0} \in H^{s}(\bbR^{d})$, $\phi_{0} \in H^{s+1}(M)$,  there exists $T > 0$ which is independent of $c>1$, we have as $c \rightarrow 0$ 
$$
\sup_{t \in [0, T]}  \|\rho^{c} - \rho^{0}\|_{L^{1}}  \rightarrow 0,
$$
where $\rho^{0} \in \cC( [0, T], H^{s}(\bbR^{d}))$ is the solution to the heat equation with the initial value $\rho_{0}$. 
}
{\item
Given the initial value $(\rho_{0}, \phi_{0}) \in TP^{\infty}_{2}(M, \mu)$ satisfying $\|\rho_{0} - \bar{\rho}\|_{H^{s}(M)} + \|\phi_{0}\|_{H^{s+1}(M)} \leq \delta$ for some $\delta > 0$, there exists $T > 0$ which is independent of $c>1$, such that as $c \rightarrow 0$, $\rho^{c}$ converges to the solution to the heat equation $\rho^{0} \in \cC(\bbR^{+}, \bar{\rho} + H^{s}(\bbR^{d}))$ in $\cC([0, T], H^{s'}(B_{R}))$.
%
% \beqnas
%\partial_t \rho^{0} &=& \Delta \rho^{0},\\
%\rho^{0}(0, x) &=& \rho_{0}(x).
%\eeqnas
}
{\item
Given initial value $(\rho_{0}, \phi_{0}) \in TP^{\infty}_{2}(M, \mu)$ with $\log \rho_{0} \in H^{s}(M)$, $\phi_{0} \in H^{s+1}(M)$, there exists $T > 0$ which is independent of $c>1$, such that $\rho^{c}dx$ weakly converges to $\rho^{\infty}dx$ in $\cC([0, T], \cP(\bbR^{d}))$ and $\phi^{c}$ converges to $\phi^{\infty}$ in $\cC([0, T], \cC^{1}(B_{R}))$ as $c \rightarrow \infty$, with $(\rho^{\infty}, \phi^{\infty})$ satisfying 
\beqna
\label{eq111}\partial \rho^{\infty} + \nabla \cdot (\rho^{\infty} \nabla \phi^{\infty}) &=& 0,\\
\label{eq112}\partial_t \phi^{\infty} + \frac{1}{2}|\nabla \phi^{\infty}|^{2} &=& 0,\\
\nonumber \rho^{\infty}(0, x) = \rho_{0}(x), \ \  \phi^{\infty}(0, x) &=& \phi_{0}(x).
\eeqna
}
\eitem
\ethm
\brmq
Consider the initial value problem to the following Hamilton-Jacobi equation
\beqnas
\partial_t \phi + \frac{1}{2}|\nabla \phi|^{2} &=& 0, \\
\phi(0, \cdot) &=& \phi_{0}(\cdot). 
\eeqnas
When $\phi_{0} \in W^{2, \infty}$, by the method of characteristics, there exists $T^{*} > 0$ such that the local solution satisfies $\phi \in W^{2, \infty}([0,T^{*}], \bbR^{d})$. Therefore, if the initial value in the above theorem satisfying $\phi_{0} \in H^{s+1}$ with $s > {d \over 2} +1$, we conclude that $\phi^{\infty} \in W^{2, \infty} \cap \cC^{1}$,  since $H^{s+1}$ can be embedded into $W^{2, \infty}$. 
\ermq

\subsection{The model of deformation of flows on $TP_2^\infty(\mathbb{R}^m, dx)$} \label{ref.model}

Let $m\in \mathbb{N}$, $m\geq n$. In this subsection we introduce the model of Langevin deformation of flows~\eqref{ww} and \eqref{vvvv} for $V(\rho) = \int \rho \log \rho$ on $T^*P_2^\infty(\mathbb{R}^m)$.  As in the previous section, for simplicity we consider Laplace operator instead of Witten Laplacian and take $\mu = \nu$. 

Let $u$ be a positive solution of the following ODE on $(0, T]\subset [0, \infty)$
\begin{eqnarray}
c^2u''+u'={1\over 2u}, \label{model-1}
\end{eqnarray}
with given initial datas $u(0)>0$ and $u'(0)\in \mathbb{R}$. Note that, in the case $c=0$,  $u(t)=\sqrt{t}$ is a solution to $(\ref{model-1})$ on $(0, T]$ for any $T>0$, and in the case $c=\infty$, $u(t)=t$ is a solution to  $(\ref{model-1})$ on $(0, T]$ for any $T>0$. 

\bthm \label{model}  Let $u$ be a smooth solution to the ODE $(\ref{model-1})$. Let $\alpha(t)={u'(t)\over u(t)}$, and $\beta(t)$ be a smooth function such that 
\begin{eqnarray*}
c^2\dot \beta(t)=-\beta(t)+m\log u(t)+{m\over 2}\log(4\pi)-1, 
\end{eqnarray*}
with a given initial data  $\beta(0) \in \mathbb{R}$. For $x\in \mathbb{R}^m$ and $t>0$, let 
\begin{eqnarray*}
\rho(x, t)&=&{1\over (4\pi u^2(t))^{m/2}}e^{-{\|x\|^2\over 4u^2(t)}},\\
\phi(x, t)&=&{\alpha(t)\over 2}\|x\|^2+\beta(t).
\end{eqnarray*}
Then $(\rho(x, t), \phi(x, t))$ satisfies the Langevin deformation of flows~\eqref{langevin.eq1} and \eqref{langevin.eq2} on $TP_2^\infty(\mathbb{R}^m)$.
\ethm

\bpf  
Note that
\begin{eqnarray*}
\nabla\phi(x, t)=\alpha(t) x.
\end{eqnarray*}
The transport equation~\eqref{langevin.eq1} has a special solution given by
\begin{eqnarray*}
\rho(x, t)=\gamma^m(t)\rho_0(\gamma(t)x),
\end{eqnarray*}
where $\rho_0(x)$ is any probability density on $\mathbb{R}^m$ with respect to the Lebesgue measure,  and $\gamma$ is a smooth function in $t$ which will be determined later. Indeed, we have
\begin{eqnarray*}
\partial_t \rho
&=&m\gamma^{m-1}\dot\gamma\rho_0(\gamma x)+ \gamma^{m}\dot \gamma \langle\nabla\rho_0(\gamma x), x\rangle,\\
\nabla \cdot (\rho\nabla\phi)&=&m\gamma^m\alpha \rho_0(\gamma x)+\gamma^{m+1} \alpha \langle\nabla\rho_0(\gamma x), x\rangle.
%&=&m\gamma^m \alpha \rho_0(\gamma x)+\sum\limits_{i=1}^m\gamma^{m+1}\alpha x_i\partial_{x_i}\rho_0(\gamma x).
\end{eqnarray*}
Thus, $(\rho, \phi)$ satisfies the transport equation if and only if
\begin{eqnarray*}
\dot \gamma+\gamma\alpha=0.
\end{eqnarray*}
Substituting $\phi$ and $\rho$ into~\eqref{langevin.eq2} leads to
%\begin{eqnarray*}
%c^2\left(\partial_t\phi+{1\over 2}|\nabla\phi|^2\right)=-\phi-\log\rho-1,
%\end{eqnarray*}
%
\begin{eqnarray*}
c^2\left({\dot\alpha(t)\|x\|^2\over 2}+\dot\beta(t)+{\alpha^2(t)\|x\|^2\over 2}\right)=-{\alpha(t)\|x\|^2\over 2}-\beta(t)-m\log \gamma(t)-\log\rho_0(\gamma(t)x)-1.
\end{eqnarray*}
Changing the variable $y=\gamma(t)x$, we have
\begin{eqnarray*}
c^2\left({\dot\alpha(t)\|y\|^2\over 2\gamma^2(t)}+\dot\beta(t)+{\alpha^2(t)\|y\|^2\over 2\gamma^2(t)}\right)=-{\alpha(t)\|y\|^2\over 2\gamma^2(t)}-\beta(t)-m\log \gamma(t)-\log\rho_0(y)-1.
\end{eqnarray*}
That is
\begin{eqnarray*}
\left[c^2\left(\dot\alpha(t)+a^2(t)\right)+\alpha(t)\right]{\|y\|^2\over 2\gamma^2(t)}+c^2\dot\beta(t)=-\beta(t)-m\log \gamma(t)-\log\rho_0(y)-1.
\end{eqnarray*}
In particular, taking
\begin{eqnarray*}
\rho_0(y)={1\over (4\pi )^{m\over 2}}e^{-{\|y\|^2\over 4}},
\end{eqnarray*}
we can choose $\alpha(t)$ and $\beta(t)$ by solving  the following ODEs 
\begin{eqnarray}
& &c^2\left(\dot\alpha(t)+\alpha^2(t)\right)+\alpha(t)={\gamma^2(t)\over 2},\label{aaa}\\
& &c^2\dot \beta(t)=-\beta(t)-m\log \gamma(t)+{m\over 2}\log(4\pi)-1.\label{bbb}
\end{eqnarray}

Let $u(t)=e^{\int_0^t \alpha(s)ds}$, and assume $\gamma(0)=1$. Then
$\alpha={u'\over u}$, $\dot\alpha={u''\over u}-{u'^2\over u^2}$, $\gamma(t)={1\over u(t)}$, and Eq.$(\ref{aaa})$ is equivalent to 
\begin{eqnarray*}
c^2\left({u''\over u}-{u'^2\over u^2}+{u'^2\over u^2}\right)+{u'\over u}={1\over 2u^2}.
\end{eqnarray*}
Thus $u$ satisfies the following nonlinear ODE 
\begin{eqnarray*}
c^2u''+u'={1\over 2u}.
\end{eqnarray*}
Note that, in the case $c=0$, we can take $\alpha(t)={1\over 2t}$, $u(t)=\sqrt{t}$ and $\gamma(t)={1\over \sqrt{t}}$, $t\in [0, T)$; and in 
the case $c=\infty$, we can take $\alpha(t)={1\over t}$,  $u(t)=t$ and $\gamma(t)={1\over t}$, $t\in [0, T)$.

\epf

\brmq
It is worth to point out that the equation $(\ref{model-1})$ can be realized by finite dimensional Langevin deformation of flows~\eqref{Bis-4} and \eqref{Bis-5} on $T\mathbb{R} \setminus\{0\}$.  To this end, take $V(x) = -{1\over 2}\log |x|$ in~\eqref{Bis-5}, such that $V'(x) = -{1\over 2x}$. Then $x$ satisfies the Langevin equation
\begin{eqnarray*}
c^2\ddot x  + \dot x = {1\over 2x}.
\end{eqnarray*}
%
%Let $V$ be a smooth function on $\mathbb{R}^+\setminus\{0\}=(0, \infty)$.  Consider the following Newton-Langevin equation on $T^*\mathbb{R}^+\setminus\{0\}=(0, \infty)\times \mathbb{R}$
%\begin{eqnarray*}
%\dot x&=&{v\over c}\\
%\dot v_t&=&-{v\over c^2}+{\nabla V(x)\over c}.
%\end{eqnarray*}
%Then $x$ satisfies the Langevin equation
%\begin{eqnarray*}
%c^2\ddot x=-\dot x+\nabla V(x).
%\end{eqnarray*}
%In particular, taking
%$$V(x)=-{1\over 2}\log x, \ \ x>0,$$
%which is Lipschitz on $[\delta, +\infty)$ for any $\delta>0$, we have
%$$\nabla V(x)=-{1\over 2x}.$$
%Thus
%$$c^2\ddot x+\dot x=-{1\over 2x}$$
%can be realized as the solution of the Newton-Langevin equation on $T^*\mathbb{R}^+\setminus\{0\}=(0, \infty)\times \mathbb{R}$. 
Given any initial position $x(0)>0$ and initial velocity $x'(0)\in \mathbb{R}$, there exists a unique solution $x(t)$ to the above equations on a small interval $[0, T)\subset [0, \infty)$  with the given initial datas $x(0)$ and $x'(0)$.

\ermq
%
%\brmq
%Let $v(t)=u(T-t)$. Then $v'(t)=-u'(T-t)$, and $v''(t)=u''(T-t)$.  Moreover, $v$ satisfies the following nonlinear ODE on the positive real line
%\begin{eqnarray*}
%c^2v''-v'=-{1\over 2v}.
%\end{eqnarray*}
%Note that, in the case $c=0$, we have $v'={1\over 2v}$, for which we can take $v(t)=\sqrt{t}$, $\gamma(t)={1\over \sqrt{t}}$, $\alpha(t)=-{1\over 2t}$; and in the case $c=\infty$, we have $v''=0$, for which can take $v(t)=t$,
%$\gamma(t)={1\over t}$, $\alpha(t)=-{1\over t}$.
%\ermq
%
\begin{proposition} \label{mmm} Let $(\rho_m, \phi_m)$ be defined as Theorem \ref{model}. Then 
\begin{eqnarray*}
{\rm Ent}(\rho_m(t))&=&-{m\over 2}(1+\log(4\pi u^2(t))),\\
H(\rho_m(t), \phi_m(t))&=&{mc^2u'(t)^2\over 2}-{m\over 2}(1+\log(4\pi u^2(t))).
\end{eqnarray*}
\end{proposition}
\bpf Indeed, $\rho_m(t)dx$ is the Gaussian distribution with variance  $u^2(t)$. 
\epf

\section{The $W$-entropy formula for the Langevin deformation of flows on Wasserstein space}

In this section we prove Perelman's type $W$-entropy formulas and monotonicity results for the Langevin deformation of flows on the Wasserstein space over compact Riemannian manifolds or on Euclidean spaces. Our results can be regarded as an interpolation between Langevin deformation of flowsthe $W$-entropy formula for the geodesic flow and the the heat flow on the Wasserstein space over compact Riemannian manifolds. We also provide the rigidity models for the $W$-entropy of the  Langevin deformation of flows on the Wasserstein space over complete noncompact Riemannian manifolds. 
Langevin deformation of flows
In this section we only consider the case where the potential $V$ is the Boltzmann-Shannon entropy $V(\rho) = \Ent(\rho) = \int_M \rho \log \rho d\mu$. {We will study the case of the R\'enyi entropy $V(\rho)={1\over m-1}\int_M \rho^{m} d\mu$ with $m\neq 1$  in future}. 

We first extend Proposition~\ref{VHH} to the Langevin deformation of flows $(\rho(t), \phi(t))$ on $TP_2^\infty(M, \mu)$.

%
%\bthm \label{Th-LL-T2} Let $c>0$. Let $\rho, \phi$ be a smooth solution of the following equations
%\begin{eqnarray*}
%\partial_t \rho+\nabla_\mu^*(\rho\nabla \phi)&=&0,\\
%c^2\left({\partial\phi\over \partial t}+{1\over 2}|\nabla \phi|^2\right)&=&-\phi+\log \rho+1.
%\end{eqnarray*}
%Let
%\begin{eqnarray*}
%H(\rho, \phi)={c^2\over 2}\int_M |\nabla\phi|^2\rho d\mu+{\rm Ent}(\rho).
%\end{eqnarray*}
%Then
%\begin{eqnarray*}
%{d\over dt}H(\rho)&=&2\int_M \nabla\phi\cdot\nabla \rho d\mu-\int_M |\nabla\phi|^2\rho d\mu\\
%&=&-\int_M (2L\phi+|\nabla \phi|^2) \rho d\mu,
%\end{eqnarray*} 
%and
%\begin{eqnarray}
%{d^2\over dt^2}H(\rho, \phi)=2\int_M \left[c^{-2}|\nabla \phi-\nabla\log\rho|^2+|{\rm Hess}\phi|^2+Ric(L)(\nabla\phi, \nabla\phi)\right]\rho d\mu.\label{2ndH}
%\end{eqnarray}
%\ethm

\bprop  Let $c>0$. Let $\rho, \phi$ be a smooth solution of the  Langevin deformation of flows
\begin{eqnarray*}
\partial_t \rho-\nabla_\mu^*(\rho\nabla \phi)&=&0,\\
c^2\left({\partial\phi\over \partial t}+{1\over 2}|\nabla \phi|^2\right)&=&-\phi-\log \rho-1.
\end{eqnarray*}
Then we have
\beqna\label{went.flow.infty}
\nonumber{d^2\over dt^2}\Ent(\rho)+ {1\over c^2}{d\over dt}\Ent(\rho) + {1\over c^2}\|\nabla \Ent (\rho)\|^2  =  \int_M \left[|{\rm Hess}~\phi|^2+\Ric(L)(\nabla\phi, \nabla\phi)\right]\rho d\mu.
\\
\eeqna
In particular, if $\Ric(L)\geq K$, then for all $c>0$, we have
\beqna\label{went.nonnegative}
{d^2\over dt^2}\Ent(\rho)+  {1\over c^2}{d\over dt}\Ent(\rho) + {1\over c^2}|\nabla \Ent (\rho)|^2 
\geq {1\over n}\int_M |\Delta \phi|^2 \rho d\mu+K\int_M |\nabla\phi|^2\rho d\mu.
\eeqna
\eprop

\bpf
By the same argument as used in the proof of Proposition~\ref{VHH} and \eqref{hess.ent}, we can prove
\beqnas
& &{d^2\over dt^2}\Ent(\rho)+  {1\over c^2}{d\over dt}\Ent(\rho) + {1\over c^2}|\nabla \Ent (\rho)|^2 \\
 &=& \Hess~ \Ent (\dot{\rho}, \dot{\rho}) = \int_M \left[|{\rm Hess}~\phi|^2+\Ric(L)(\nabla\phi, \nabla\phi)\right]\rho d\mu, 
\eeqnas
which leads to \eqref{went.nonnegative} if $\Ric(L)\geq K$ by using the inequality $|{\rm Hess}~\phi|^2\geq {1\over n}|\Delta\phi|^2$.
\epf

%
%
%
%{\blue
%\bthm\label{Th-LLL-1} Under the same notation as in Theorem \ref{Th-LL-T2}, for any $c>0$, we have
%\begin{eqnarray*}
%\left({d^2\over dt^2}+{1\over c^2}{d\over dt}\right) {\rm Ent} (\rho(t))&=&{1\over c^2} \int_M {|\nabla\rho|^2\over \rho} d\mu+\int_M [|{\rm Hess}\phi |^2+Ric(L)(\nabla\phi, \nabla\phi)]\rho d\mu,\\
%\left({d^2\over dt^2}+{2\over c^2}{d\over dt}\right) H(\rho(t), \phi(t))&=&{2\over c^2} \int_M {|\nabla\rho|^2\over \rho} d\mu+\int_M [|{\rm Hess}\phi |^2+Ric(L)(\nabla\phi, \nabla\phi)] \rho d\mu.
%\end{eqnarray*}
% Let 
%\begin{eqnarray*}
%W_{H, c}(\rho(t), \phi(t))&=&H(\rho(t), \phi(t))+{c^2(1-e^{2t\over c^2})\over 2}{d\over dt}H(\rho(t), \phi(t)),\\
%W_{{\rm Ent}, c}(\rho(t))&=& {\rm Ent}(\rho(t))+c^2(1-e^{t\over c^2}){d\over dt}{\rm Ent}(\rho(t)).
%\end{eqnarray*}
%Then 
%\begin{eqnarray*}
%{d\over dt}W_{H, c}(\rho(t), \phi(t))&=&(1-e^{2t\over c^2})\int_M {|\nabla\rho|^2\over \rho}d\mu+\int_M [|{\rm Hess}\phi|^2+Ric(L)(\nabla\phi, \nabla\phi)]\rho d\mu,\\
%{d\over dt}W_{{\rm Ent}, c}(\rho(t))&=&(1-e^{t\over c^2})\int_M {|\nabla\rho|^2\over \rho}d\mu+\int_M [|{\rm Hess}\phi|^2+Ric(L)(\nabla\phi, \nabla\phi)]\rho d\mu.
%\end{eqnarray*}
%In particular, if $Ric(L)\geq 0$, then for all $c>0$, we have
%\begin{eqnarray*}
%{d\over dt}W_{H, c}(\rho(t), \phi(t))\geq 0, \ \ \ \forall t\geq 0,
%\end{eqnarray*}
%and
%\begin{eqnarray*}
%{d\over dt}W_{{\rm Ent}, c}(\rho(t))\geq 0, \ \ \ \forall t\geq 0.
%\end{eqnarray*}
%\ethm
%
%\bpf
%By Theorem \ref{Th-LL-T2}, we can prove Theorem \ref{Th-LLL-1} similarly to the one of Theorem \ref{WHV}. $ 
%\epf
%}
%

\subsection{Proof of Theorem~\ref{MT1}}
%
%\bthm\label{Main} Let $c\in [0, \infty)$, and $(\rho(t), \phi(t))$ be smooth solution to the Langevin deformation of flow. Let $\alpha(t)=(\log u)'$. 
%In the case $m=n$, and $L=\Delta$, we have
%\begin{eqnarray}\label{went.1}
%\nonumber\ \ \ \ & &{d^2\over dt^2} {\rm Ent}(\rho(t))+\left(2\alpha(t)+{1\over c^2}\right){d\over dt}{\rm Ent}(\rho(t))  + {1\over c^2}\|\nabla \Ent(\rho(t))\|^2 +n\alpha^2(t)\\
%&=& \int_M \left[\left|{\rm Hess}\phi-\alpha(t)g\right|^2+\Ric(\nabla \phi, \nabla \phi) \right]\rho d\nu,
%\end{eqnarray}
%and in the case $m>n$, we have
%\begin{eqnarray}\label{went.2b}
%& &{d^2\over dt^2} {\rm Ent}(\rho(t))+\left(2\alpha(t)+ {1\over c^2}\right){d\over dt}{\rm Ent}(\rho(t))+ {1\over c^2}\|\nabla \Ent(\rho(t))\|^2+m\alpha^2(t) \nonumber\\
%%&=&  \int_M \left[|{\rm Hess}\phi|^2+\Ric(L)(\nabla\phi, \nabla\phi)\right]\rho d\mu - 2\alpha(t)\int_M L\phi \rho d\mu+m\alpha^2(t)\\
% &=&\int_M \left[\left|{\rm Hess}\phi-\alpha(t)g\right|^2+\Ric_{m, n}(L)(\nabla \phi, \nabla \phi) \right]\rho d\mu +(m-n)\int_M \left|\alpha(t)+{\nabla \phi\cdot \nabla f\over m-n}\right|^2 \rho d\mu,\nonumber\\
%\end{eqnarray}
%where 
%\beqna
%\|\nabla \Ent (\rho(t))\|^2=\int_M {|\nabla\rho(t)|^2\over \rho(t)}d\mu.
%\eeqna 
%
%\ethm

\bpf
The proof has the same spirit as those for Theorem \ref{MT2} and Theorem \ref{MT3}. By Otto's calculation on $P_2(M, \mu)$, we have
\begin{eqnarray*}
{d\over dt} {\rm Ent} (\rho(t))=\nabla{\rm Ent} (\rho(t)) \cdot \dot \rho(t)%=\int_M \nabla (\log \rho+1)\cdot \nabla \phi \rho d\mu=\int_M \nabla \phi \cdot \nabla \rho d\mu
=-\int_M L\phi \rho d\mu.
\end{eqnarray*}
By Theorem \ref{Ent-Hess}, we have
\begin{eqnarray*}
& &{d^2\over dt^2}{\rm Ent}(\rho(t))=\nabla^2 {\rm Ent}(\rho(t))(\dot \rho(t), \dot\rho(t))+\nabla{\rm Ent}\cdot\ddot \rho(t)\\
&=&\int_M \left(|{\rm Hess}~\phi|^2+Ric(L)(\nabla\phi, \nabla \phi)\right)\rho d\mu+\int_M \nabla (\log \rho+1)\cdot \nabla  \left({\partial_t}\phi+{1\over 2}|\nabla \phi|^2\right)\rho d\mu\\
&=&\int_M \left(|{\rm Hess}~\phi|^2+Ric(L)(\nabla\phi, \nabla \phi)\right)\rho d\mu+{1\over c^2}\int_M \nabla \rho \cdot \nabla(-\phi-\log\rho-1)d\mu\\
&=&\int_M \left(|{\rm Hess}~\phi|^2+Ric(L)(\nabla\phi, \nabla \phi)\right)\rho d\mu-{1\over c^2}\int_M {|\nabla\rho|^2\over \rho}d\mu+{1\over c^2}\int_M L\phi \rho d\mu.
\end{eqnarray*}
Then, similarly to the proof of Theorem \ref{MT2}, for the case $m=n$ and $L=\Delta$, we derive that
\beqnas
& &{d^2\over dt^2} {\rm Ent}(\rho(t))+\left(2\alpha(t)+c^{-2}\right){d\over dt}{\rm Ent}(\rho(t))  + c^{-2}\|\nabla \Ent(\rho(t))\|^2 +n\alpha^2(t)\\
&=&  \int_M \left[|{\rm Hess}\phi|^2+\Ric(\nabla\phi, \nabla\phi)\right]\rho d\nu - 2\alpha(t)\int_M \Delta \phi \rho d\nu+n\alpha^2(t)\\
&=&  \int_M \left[|{\rm Hess}\phi - \alpha(t)g|^2+\Ric(\nabla\phi, \nabla\phi)\right]\rho d\nu.
\eeqnas
For the case  $m>n$, similar computations lead to
\beqnas
& &{d^2\over dt^2} {\rm Ent}(\rho(t))+\left(2\alpha(t)+ c^{-2} \right){d\over dt}{\rm Ent}(\rho(t))+ c^{-2}\|\nabla \Ent(\rho(t))\|^2+m\alpha^2(t)\\
&=&  \int_M \left[|{\rm Hess}\phi|^2+\Ric(L)(\nabla\phi, \nabla\phi)\right]\rho d\mu - 2\alpha(t)\int_M L\phi \rho d\mu+m\alpha^2(t)\\
&=&\int_M \left[\left|{\rm Hess}\phi-\alpha(t)g\right|^2+\Ric_{m, n}(L)(\nabla \phi, \nabla \phi) \right]\rho d\mu +(m-n)\int_M \left|\alpha(t)+{\nabla \phi\cdot \nabla f\over m-n}\right|^2 \rho d\mu,
\eeqnas
where the last line is due to
\begin{eqnarray*}
& &\Ric(L)(\nabla \phi, \nabla\phi)+2\alpha\nabla \phi\cdot\nabla f+(m-n)\alpha^2\\
&=&\Ric_{m, n}(L)(\nabla\phi, \nabla \phi)+{|\nabla\phi \cdot\nabla f|^2\over m-n}+2\alpha \nabla \phi\cdot\nabla f+(m-n)\alpha^2\\
&=& \Ric_{m, n}(L)(\nabla\phi, \nabla \phi)+(m-n)\left|\alpha+{\nabla \phi\cdot \nabla f\over m-n}\right|^2.
\end{eqnarray*}
This completes the proof of Theorem \ref{MT1}. \epf

\subsection{The $W$-entropy for the Langevin deformation and rigidity model}

Based on the formulas~\eqref{wentc} and \eqref{went.2b}, we now introduce the $W$-entropy for the  Langevin deformation of flows, as stated in Section 1. More precisely, we define the $W$-entropy for the  Langevin deformation of flows as follows: for all $0<t_0<t<\infty$, 
\begin{eqnarray}
W_{c}(\rho(t))- W_{c}(\rho(t_0))&=& {d\over dt}{\rm Ent}(\rho(t))+\int_{t_0}^t \left( 2\alpha(s)+{1\over c^2}\right) {d\over ds} {\rm Ent}(\rho(s))ds\nonumber\\
& &\hskip3cm  +{1\over c^2}\int_{t_0}^t \|\nabla \Ent (\rho(s))\|^2 ds,
\label{went}
\end{eqnarray}
such that its time derivative is given by
\beqna\label{went}
{d\over dt}W_{c}(\rho(t)) = {d^2\over dt^2} {\rm Ent}(\rho(t))+\left(2\alpha(t)+{1\over c^2}\right){d\over dt}{\rm Ent}(\rho(t))  + {1\over c^2}\|\nabla \Ent (\rho(t))\|^2.
\eeqna
%{\blue for the backward flow
%\begin{eqnarray*}
%{d\over dt}W_{c}(\rho(t))={d^2\over dt^2}{\rm Ent}(\rho(t))+\left(2\alpha(t)+ c^{-2}\right){d\over dt}{\rm Ent}(\rho(t))- c^{-2}\int_M {|\nabla\rho|^2\over \rho}d\mu.
%\end{eqnarray*}
%}
%
%The following result provides the rigidity model for the entropy formula in 
% Theorem \ref{Main}, which achieves the optimal constant $-m\alpha^2(t)$.  See also Section 7.2 below.
%
\bprop For the model $(\rho_m, \phi_m)$ on $(\mathbb{R}^m, dx)$, we have
\begin{eqnarray}\label{went.model}  {d\over dt}W_{c}(\rho_m(t))
%&=&{d^2\over dt^2}{\rm Ent}(\rho_m(t))+\left(2\alpha(t)+ c^{-2}\right){d\over dt}{\rm Ent}(\rho_m(t)) + c^{-2}|\nabla \Ent (\rho_{m}(t))|^{2}\\
=-m\alpha^2(t).
\end{eqnarray}
\eprop
\bpf  
Indeed\begin{eqnarray*}
{\rm Hess}\phi_m=\alpha(t) g,  \ \ \ \ \Ric=0,
 \end{eqnarray*}
from which and Theorem \ref{MT1}, we derive \eqref{went.model}.
\epf

\subsection{Comparison theorem for the $W$-entropy for Langevin deformation}

As a direct consequence of Theorem \ref{MT1} and \eqref{went.model}, we  can derive the following comparison theorem  for the $W$-entropy for Langevin deformation of flows.

\bthm \label{MainM} (i.e., Theorem \ref{MainMM})
\begin{eqnarray}\label{went.2}
\nonumber {d\over dt}(W_{c}(\rho(t))-W_{c}(\rho_m(t)))&=&\int_M  \left|{\rm Hess}\phi-\alpha(t)g\right|^2\rho d\mu+\int_M \Ric_{m, n}(L)(\nabla\phi, \nabla \phi)\rho d\mu\\
& &\hskip2cm +{1\over m-n}\int_M \left|\nabla f\cdot\nabla\phi+(m-n)\alpha(t)\right|^2\rho d\mu.\nonumber\\
\end{eqnarray}
In particular, if $\Ric_{m, n}(L)\geq 0$, then for all $t>0$, we have the comparison theorem
\begin{eqnarray}
{d\over dt}W_{c}(\rho(t))\geq {d\over dt}W_{c}(\rho_m(t)).
\end{eqnarray}
\ethm

 We would like to point out that, similarly to the proof of Theorem \ref{entropy-geo-noncompact}, we can extend the above $W$-entropy formula  to complete Riemannian manifolds with bounded geometry condition. To save the length of the paper, we omit the detail of this technical part. In view of this, the rigidity model of the $W$-entropy for the Langevin deformation of flows on complete Riemannian manifolds with non-negative Ricci curvature or complete and weighted Riemannian manifolds  with $CD(0, m)$-condition should be  $M=\mathbb{R}^m$, $m\in \mathbb{N}$, and $(\rho_m, \phi_m)$ given in Theorem \ref{model}. 
%
%
%\bprop \label{mmm} Let $(\rho_m, \phi_m)$ be the special solution in Theorem~\ref{model}. Then 
%\begin{eqnarray*}
%{\rm Ent}(\rho_m(t))&=&-{m\over 2}(1+\log(4\pi u^2(t))),\\
%H(\rho_m(t), \phi_m(t))&=&{mc^2u'(t)^2\over 2}-{m\over 2}(1+\log(4\pi u^2(t))).
%\end{eqnarray*}
%\eprop
%\bpf 
%Indeed, $\rho_m(t)$ Theorem \ref{model} is the Gaussian heat kernel at time {\blue $2u^2(t)$}. 
%{\blue
%$$
%\int_{\bbR^{m}}|x|^{2}\rho_m(t) = 2mu^2(t).
%$$
%}
%\epf
%

To end this paper, let us give the following remark which illustrate that the $W$-entropy formula~\eqref{went.2b} is an interpolation between those of the geodesic flow and of the the heat flow on the Wasserstein space. 

\brmq 
Indeed, when $c \rightarrow \infty$, we have $u(t)  \rightarrow  t$, $\alpha(t)  \rightarrow  {1\over t}$. This yields
%{\blue{\begin{eqnarray*}
%{d\over dt}{\rm Ent}(\rho(t))=\nabla\Ent(\rho(t))\cdot\dot\rho(t)=-\int_M L\phi \rho d\mu, 
%\end{eqnarray*}
%and since $\ddot \rho(t)=0$, we have
%\begin{eqnarray*}
%{d^2\over dt^2}{\rm Ent}(\rho(t))&=&\nabla^2 \Ent(\rho(t))(\dot\rho(t), \dot\rho(t))+\nabla \Ent(\rho(t))\cdot \ddot{\rho}(t)\\
%&=&\nabla^2 \Ent(\rho(t))(\dot\rho(t), \dot\rho(t)).
%\end{eqnarray*}
%This yields
\beqna\label{went}
\lim\limits_{c\rightarrow \infty}{d\over dt}W_{c}(\rho(t))
 %&=&\lim\limits_{c\rightarrow \infty}\left[{d^2\over dt^2} {\rm Ent}(\rho(t))+\left(2\alpha(t)+{1\over c^2}\right){d\over dt}{\rm Ent}(\rho(t))  + {1\over c^2}\|\nabla \Ent (\rho(t))\|^2\right]\nonumber\\
={d^2\over dt^2} {\rm Ent}(\rho(t))+{2\over t}{d\over dt}{\rm Ent}(\rho(t))\nonumber.
\eeqna
%}}
 Hence the time derivative of the $W$-entropy for the geodesic flow should  be given by
\begin{eqnarray*}
{d\over dt}W(\rho(t))={d^2\over dt^2}{\rm Ent}(\rho(t))+{2\over t}{d\over dt}{\rm Ent}(\rho(t)),
\end{eqnarray*}
and the $W$-entropy formula reads as 
(i.e., Theorem \ref{MT2})
\begin{eqnarray*}
{d\over dt}(W(\rho(t))-W(\rho_m(t)))&=&\int_M  \left|{\rm Hess}\phi-{g\over t}\right|^2\rho d\mu+\int_M \Ric_{m ,n}(L)(\nabla\phi, \nabla \phi)\rho d\mu\\
& &\hskip2cm +{1\over m-n}\int_M \left|\nabla f\cdot\nabla\phi+{m-n\over t}\right|^2\rho d\mu.
\end{eqnarray*}
On the other hand, when $c \rightarrow  0$, we have $u(t) \rightarrow \sqrt{t}$, $\alpha(t) \rightarrow {1\over 2t}$. Moreover, the convergence results lead to
\begin{eqnarray*}
{d\over dt}{\rm Ent}(\rho(t)) =\nabla \Ent(\rho(t))\cdot \dot \rho(t)=-\|\nabla \Ent(\rho(t))\|^2.
%{d^2\over dt^2}{\rm Ent}(\rho(t)) &=&2 {\rm Hess}~\Ent(\dot{\rho}, \dot{\rho}),
\end{eqnarray*}
This yields
 \begin{eqnarray*}%\label{went}
\lim\limits_{c\rightarrow 0}{d\over dt}W_{c}(\rho(t)) 
%&=&\lim\limits_{c\rightarrow 0}\left[{d^2\over dt^2} {\rm Ent}(\rho(t))+\left(2\alpha(t)+{1\over c^2}\right){d\over dt}{\rm Ent}(\rho(t))  + {1\over c^2}\|\nabla \Ent (\rho(t))\|^2\right]\nonumber\\
={d^2\over dt^2} {\rm Ent}(\rho(t))+{1\over t}{d\over dt}{\rm Ent}(\rho(t)).
\end{eqnarray*}
Hence the time derivative of the $W$-entropy for the gradient flow should be given by
\begin{eqnarray*}
{d\over dt}W(\rho(t))={d^2\over dt^2}{\rm Ent}(\rho(t))+{1\over t}{d\over dt}{\rm Ent}(\rho(t)),
\end{eqnarray*}
and the $W$-entropy formula reads as 
(i.e., Theorem \ref{MT3}) 
\begin{eqnarray*}
{d\over dt}(W(\rho(t))-W(\rho_m(t)))&=&\int_M  \left|{\rm Hess}\phi-{g\over 2t}\right|^2\rho d\mu+\int_M \Ric_{m, n}(L)(\nabla\phi, \nabla \phi)\rho d\mu\\
& &\hskip2cm +{1\over m-n}\int_M \left|\nabla f\cdot\nabla\phi+{m-n\over 2t}\right|^2\rho d\mu.
\end{eqnarray*} 
{Note that the factor in front of ${d\over dt}{\rm Ent}(\rho(t)$ in the right hand side of the $W$-entropy formula  is ${1\over t}$ for $c=\infty$ but is ${2\over t}$ for $c=0$. For its reason, see Section 3.5. }
\ermq

%
%\section{Section 7}\label{Section 7}

\section{Convergence of the Langevin deformation}

\subsection{The convergence of the deformation flow in finite dimension}
Let $(x, v)$ be the smooth solution to the following Langevin equation on $\bbR^{d} \times [0, T]$
\beqna
\label{deformation.finite1}\dot{x} &=& v,\\
\label{deformation.finite2}c^{2}\dot{v} &=& - v - \nabla V(x), 
\eeqna
with the initial data
\beqnas
x(0) = x_{0}, \ \ \  v(0) = v_{0}.
\eeqnas

The existence and uniqueness of solution to  \eqref{deformation.finite1} and \eqref{deformation.finite2} follow from the Cauchy-Lipschitz theorem in ODE 
theory. The following convergence results can be proved by standard $C^2$-estimates. To save the length of this paper, we omit the details of proof.

\bthm\label{theorem1}\label{cthm1}
Let $V \in \cC^{2}(\bbR^{d})$ and assume that $\nabla V$ is Lipschitz, i.e., there exists a constant $K > 0$ such that 
$$|\nabla V(x_1) - \nabla V(x_2)| \leq K|x_1 - x_2|, \ \ \forall~ x_1, x_2 \in \bbR^{d}.$$ Then there exists $T > 0$, such that the equations \eqref{deformation.finite1} and \eqref{deformation.finite2} have a unique solution $(x_t, v_t)\in  \cC([0, T], \cC^{2}(\bbR^{d}))\times \cC^{1}(\bbR^{d})) $. 
Moreover, when $c\rightarrow 0$, the solution to \eqref{deformation.finite1},   \eqref{deformation.finite2}  uniformly converges on $[0, T]$ to the solution to the gradient flow of $V$ 
\beqna\label{eqs.1}
 & &\dot{y}=-\nabla V(y), \\
& & y(0)= x_{0},
\eeqna
and when  $c \rightarrow \infty$, the solution to \eqref{deformation.finite1},   \eqref{deformation.finite2}  uniformly converges on $[0, T]$ to the solution to the geodesic equation 
\beqna\label{eqs.2}
\ddot{z} &=& 0, \\
z(0) &=& x_{0}, \ \ \ \dot{z}(0) = v_{0}.
\eeqna
Indeed, letting 
$y(t) \in  \cC([0, T], \cC^{2}(\bbR^{d}))$  be the unique solution to \eqref{eqs.1},  
and $z(t) \in \cC([0, T], \cC^{2}(\bbR^{d}))$ be the unique solution to \eqref{eqs.2}, then 
\beqnas
\label{eq1}\lim_{c \rightarrow 0}\sup_{t \in [0, T]}|x(t) - y(t)| &=& 0, \\
\label{eq2}\lim_{c \rightarrow \infty}\sup_{t \in [0, T]}|x(t) - z(t)|&=& 0.
\eeqnas
Moreover, if $\dot x(0)=\dot y(0)=\dot z(0)$, $\ddot x(0)=\ddot y(0)=\ddot z(0)$, and if $\nabla^2 V$ is $K_2$-Lipschitz and uniformly bounded, then 

\beqnas
\label{eq12}\lim_{c \rightarrow 0}\sup_{t \in [0, T]}[ |\dot x(t)-\dot y(t)|+|\ddot x(t) - \ddot y(t)|] &=& 0, \\
\label{eq22}\lim_{c \rightarrow \infty}\sup_{t \in [0, T]}[ |\dot x(t)-\dot z(t)| +
|\ddot x(t) |] &=& 0.
\eeqnas

  \ethm
\subsection{The convergence of the Langevin deformation of flows on Wasserstein space: $c \rightarrow 0$}
In this section we consider the Langevin deformation of flows on $\bbR^{d}$ or compact manifold $M$,
\beqna
\label{langevin.eq1}\partial_t \rho + \nabla \cdot (\rho \nabla \phi) &=& 0, \\
\label{langevin.eq2}c^{2}(\partial_t \phi + \frac{1}{2}|\nabla \phi|^{2}) &=& -\phi - \log \rho -1,
\eeqna
with initial value
$$
\rho(0, x) = \rho_{0}(x), \ \ \phi(0, x) = \phi_{0}(x).
$$
We will first discuss the convergence in Euclidean space when $c$ approachs 0 and infinity respectively, and then extend the results to compact manifolds. For simplicity, we consider 
the  Laplace-Beltrami operator here instead of the Witten Laplacian. Under suitable assumptions on the potential function $f$, it is easy to prove that the convergence results also hold for the Witten Laplacian. 

Consider the Cauchy problem of \eqref{langevin.eq1} and \eqref{langevin.eq2}. In this section we start with the convergence results for the compressible Euler equation with damping. 
Let $\nabla \phi = u$. Taking differential on both sides of \eqref{langevin.eq2}, we obtain
\beqna
\label{euler.eq1}\partial_t \rho^{c} + \nabla \cdot (\rho^{c} u^{c}) &=& 0, \\
\label{euler.eq2}c^{2}(\partial_t u^{c} + u^{c} \cdot \nabla u^{c}) &=& - u^{c} - \frac{\nabla \rho^{c}}{\rho^{c}}.
\eeqna
In the following we use the notation $(\rho^{c}, u^{c})$ instead of  $(\rho, u)$ to emphasis its relevance with $c$. 
 
\subsubsection{The relative entropy method for the convergence of local solution}
Inspired by the works of Lattanzio-Tzavaras~\cite{Lat-Tza13} and~\cite{Lat-Tza17}, we apply the relative entropy method to study the convergence of local solutions $(\rho^{c}, u^{c})$. In~\cite{Lat-Tza13}, the method was employed to study the diffusive limit of entropy weak solutions to compressible isentropic Euler equation with damping only in one dimension. In this section we apply the relative entropy method to prove the convergence of entropy solutions for Langevin deformation of flows. It has  its own interest to use the relative entropy method to prove the convergence, since for the compressible Euler equation with damping, 
$L^{\infty}$ solution is more natural (because of the shocks) than local smooth solutions. As the well posedness of entropy solution to the compressible Euler equation with damping has not yet been established in high dimension, we will focus the case on one dimension. 

For any $\rho, \bar{\rho} \in \cP(\bbR^{d})$, define the density for the Boltzmann entropy as $h(\rho) = \rho \log \rho$. Then $h''(\rho) = \frac{1}{\rho}$, and for the relative  Boltzmann entropy  $h(\rho | \bar{\rho})$ between $\rho$ and $\bar{\rho}$, we have
\beqnas
h(\rho | \bar{\rho}) &=& h(\rho) - h(\bar{\rho}) - h'(\bar{\rho})(\rho - \bar{\rho})\\
&=& (\rho - \bar{\rho})^{2}\int^{1}_{0}\int_{0}^{\tau}h''(s\rho + (1-s)\bar{\rho})dsd\tau\geq 0.
\eeqnas

Moreover, the density for the kinetic energy for $(\rho, m)$ is given by
$$
k(\rho^{c}, m^{c}) =  \frac{1}{2}\rho^{c} \|u^{c}\|^{2}.
$$
We may also define the relative kinetic energy between $(\rho, m)$  and $(\bar{\rho}, \bar{m})$ as
\beqnas
k((\rho, m) |(\bar{\rho}, \bar{m})) &=& k(\rho, m) - k(\bar{\rho}, \bar{m}) - \langle \nabla k(\rho, m), (\rho -  \bar{\rho}, m - \bar{m})\rangle\\
%&=& k(\rho, m) - k(\bar{\rho}, \bar{m}) - \langle \nabla k(\bar{\rho}, \bar{m}), (\rho -  \bar{\rho}, m - \bar{m})\rangle\\
&=& \frac{1}{2}\rho\|u - \bar{u}\|^{2}.
\eeqnas
% $\langle \nabla k(\bar{\rho}, \bar{m}), (\rho -  \bar{\rho}, m - \bar{m})\rangle =  \frac{1}{2}(\rho - \bar{\rho})|\bar{u}|^{2} + \frac{\bar{m}}{\bar{\rho}} \cdot (m - \bar{m})$.

For~\eqref{euler.eq1} and \eqref{euler.eq2},  we define the entropy as 
$$
\eta(\rho^{c}, m^{c}) =  \rho^{c} \log \rho^{c} + \frac{c^{2}}{2}\rho^{c} \|u^{c}\|^{2},
$$
and the flux 
$$
q(\rho^{c}, m^{c}) = c(1 + \log \rho^{c})m^{c} + c^{3}\frac{\|m^{c}\|^{2}}{2(\rho^{c})^{2}}m^{c}, 
$$
where $m^{c} = \rho^{c} u^{c}$. 

Let $\rho^{0}$ be the solution to the heat equation
\beqnas
\partial_t \rho^{0} &=& \Delta \rho^{0},\\
\rho^{0}(0, x) &=& \rho_{0}(x),
\eeqnas
such that $m^{0} = \nabla  \rho^{0}$.

Thus the relative entropy between $(\rho^{c}, m^{c})$ and $(\rho^{0}, m^{0})$ is given by
\beqnas
\eta(\rho^{c}, m^{c}|\rho^{0}, m^{0}) &=& \Ent(\rho^{c}|\rho^{0}) + c^{2} K(\rho^{c}, m^{c}|\rho^{0}, m^{0})\\
&=& \int_{\mathbb{R}^{d}} \rho^{c} (\log \rho^{c} - \log \rho^{0})+  \frac{c^{2}}{2}\int_{\mathbb{R}^{d}}  \rho^{c}\|u^{c} - u^{0}\|^{2}.
\eeqnas
%{\blue
%where we define the relative kinetic energy $K(\rho^{c}, m^{c}|\rho^{0}, m^{0})$ as follows
%\begin{eqnarray*}
%K(\rho^{c}, m^{c}) &=& \frac{c^{2}}{2}\int_{\mathbb{R}^{d}}  \frac{\|m^{c}\|^{2}}{\rho^{c}}, \ \ K(\rho^{0}, m^{0})={c^2\over 2} \int_{\mathbb{R}^{d}}  \frac{\|m^{0}\|^{2}}{\rho^{0}},\\
%K(\rho^{c}, m^{c}|\rho^{0}, m^{0}) &=& K(\rho^{c}, m^{c}) - K(\rho^{0}, m^{0}) - \langle \nabla K(\rho^{0}, m^{0}), (\rho^{c} -  \rho^{0}, m^{c} - m^{0})\rangle\\
%&=& \frac{c^{2}}{2}\int_{\mathbb{R}^{d}}  \rho^{c}\|u^{c} - u^{0}\|^{2}.
%\end{eqnarray*}
%}
%
%{\red
%which can be considered as the relative kinetic energy, for that we may define
%$$
%K(\rho, m) &=& \frac{c^{2}}{2}\int_{\mathbb{R}^{d}}  \frac{\|m\|^{2}}{\rho}
%$$
%as the kinetic energy of $(\rho, m)$. Moreover, in view of Riemannian metric on $\cP^{2, \infty}(\bbR^{d})$ according to Otto's calculus, we have
%$$
%\nabla K(\rho, m)  = \rho \|u\|^{2}.
%$$
%
%is viewed as a function of $(\rho^{c}, m^{c})$, and  $\langle \nabla K(\rho^{0}, m^{0}), (\rho^{c} -  \rho^{0}, m^{c} - m^{0})\rangle$ is defined  with respect to the inner product in the sense of Otto's calculus on $\cP^{2, \infty}(\bbR^{d})$.
%
%
%\beqnas
%\langle \nabla K(\rho^{0}, m^{0}), (\rho^{c} -  \rho^{0}, m^{c} - m^{0})\rangle = 
%\eeqnas
%}
%
We have the following results for the relative entropy $\eta(\rho^{c}, m^{c}|\rho^{0}, m^{0})$.
 
\bprop For smooth solution to the compressible Euler equation, we have
 \beqna\label{entropyineq}
& &  \eta(\rho^{c}, m^{c}|\rho^{0}, m^{0})(t)\nonumber\\
   &=& \int_{\mathbb{R}^{d}}  \rho (\log \rho - \log \rho^{0})+  \frac{c^{2}}{2}\int_{\mathbb{R}^{d}} \rho\|u^{c} - u^{0}\|^{2} \nonumber\\
 &\leq& \big( \eta(\rho^{c}, m^{c}|\rho^{0}, m^{0})(0) + c^{4}\int_{\mathbb{R}^{d}} \rho^{c} |\Delta u^{0} - u^{0} \cdot \nabla u^{0}|^{2} \big)e^{2\|\nabla u^{0}\|_{\infty}t}.
 \eeqna
\eprop

\bpf
Notice that for smooth solutions, we have
 \beqnas
 \partial_t \eta(\rho^{c}, m^{c}) + \frac{1}{c}\nabla_{x} \cdot q(\rho^{c}, m^{c}) = - \frac{|m^{c}|^{2}}{\rho^{c}},
 \eeqnas
and  $u^{c}$ satisfies 
$$
\partial_t u^{c} = -u^{c} \cdot \nabla u^{c} - c^{-2}u^{c} - c^{-2}\nabla \log \rho^{c}.
$$ 
%and for  entropy solutions, we have
% \beqnas
% \partial_t \eta(\rho, m) + \frac{1}{c}\nabla_{x} \cdot q(\rho, m) \leq  - \frac{|m|^{2}}{\rho}.
% \eeqnas
Thus we derive that
\begin{eqnarray}\label{eeee1}
\frac{d}{dt}\left( \int_{\mathbb{R}^{d}}  \rho^{c} \log \rho^{c} + \frac{c^{2}}{2}\int_{\mathbb{R}^{d}}  \rho^{c} |u^{c}|^{2}\right) = - \int_{\mathbb{R}^{d}}  \rho^{c}|u^{c}|^{2}.
%\nonumber  &=& \int_{\mathbb{R}^{d}}  \left( - \frac{1}{c}\nabla_{x} \cdot q(\rho^{c}, m^{c})  - \rho^{c}|u^{c}|^{2}\right)
 \end{eqnarray}

On the other hand, by the integration by parts formula, we have
%$\Delta \log \rho_0={\Delta\rho_0\over \rho_0}-|\nabla \rho_0|^2$ and integration by parts, we have 
%As for  $u^{0} = - \nabla \log \rho^{0}$, we have
%$$
%\partial_t u^{0} = - 2 u^{0} \cdot \nabla u^{0} +  \Delta u^{0}.
%$$
\beqnas
 \frac{d}{dt}\left(- \int_{\mathbb{R}^{d}}  \rho^{c} \log \rho^{0}\right)
%  &=& - \int_{\mathbb{R}^{d}}  \partial_t \rho^{c} \log \rho^{0} - \int \rho^{c} \frac{\partial_t \rho^{0}}{\rho^{0}}\\
% &=&\int \nabla\cdot (\rho^{c} u^{c})\log \rho^{0} - \int \rho^{c} \frac{\partial_t \rho^{0}}{\rho^{0}}\\
% &=& - \int_{\mathbb{R}^{d}}   \rho^{c} \langle u^{c} ,  \nabla\log \rho^{0} \rangle - \int \rho^{c} \frac{\Delta \rho^{0}}{\rho^{0}}\\
 &=& - \int_{\mathbb{R}^{d}}   \rho^{c} \langle u^{c} ,  \nabla\log \rho^{0} \rangle - \int_{\mathbb{R}^{d}}  \langle \nabla \rho^{c} , u_0\rangle+ \int_{\mathbb{R}^{d}}  \rho^{c}\langle u^{0},  \nabla \log \rho^{0}\rangle.
 \eeqnas

Moreover, since $u^{0} = - \nabla \log \rho^{0}$ satisfies
$$
\partial_t u^{0} = - 2 u^{0} \cdot \nabla u^{0} +  \Delta u^{0}, 
$$
we derive that
 \beqnas
& &\frac{d}{dt}\int_{\mathbb{R}^{d}}  \rho^{c}( \|u^{0}\|^{2} - 2u^{c} \cdot u^{0}) \\
%&=&\int_{\mathbb{R}^{d}}  \partial_t \rho^{c} ( \|u^{0}\|^{2} - 2u^{c} \cdot u^{0}) +\int_{\mathbb{R}^{d}}  \rho^{c} \partial_t ( \|u^{0}\|^{2} - 2u^{c} \cdot u^{0}) \\
%&=&-\int_{\mathbb{R}^{d}}  \nabla \cdot (\rho^{c} u^{c})( \|u^{0}\|^{2} - 2u^{c} \cdot u^{0}) +\int_{\mathbb{R}^{d}}  \rho^{c} \partial_t ( \|u^{0}\|^{2} - 2u^{c} \cdot u^{0}) \\
%&=& \int_{\mathbb{R}^{d}}  \rho^{c} u^{c} \cdot \nabla ( \|u^{0}\|^{2} - 2u^{c} \cdot u^{0}) +\int_{\mathbb{R}^{d}}  \rho^{c} \partial_t ( \|u^{0}\|^{2} - 2u^{c} \cdot u^{0}) \\
% &=&  \int_{\mathbb{R}^{d}}  \rho^{c}\langle u^{c}, 2(u^{0} - u^{c})\cdot  \nabla u^{0} - 2\nabla u^{0} \cdot u^{0}\rangle  \\
% & & +\int_{\mathbb{R}^{d}}  \rho( 2(u^{0} - u^{c}) \cdot  (\Delta u^{0} - 2\nabla u^{0} \cdot u^{0}) + 2(u^{c} \cdot \nabla u^{c} + c^{-2}u^{c} + c^{-2}\nabla \log \rho^{c}) \cdot u^{0} )\\
& =&  2 \int_{\mathbb{R}^{d}}  \rho\langle u^{c}, (u^{0} - u^{c})\cdot  \nabla u^{0}\rangle  +2 \int_{\mathbb{R}^{d}}  \rho( (u^{0} - u^{c}) \cdot  (\Delta u^{0}-2\nabla u^{0} \cdot u^{0}) + c^{-2}(u^{c} + \nabla \log \rho^{c}) \cdot u^{0} ).
 \eeqnas
 
 Combining the above terms together leads to
\begin{eqnarray*}
\frac{d}{dt}\eta(\rho^{c}, m^{c}|\rho^{0}, m^{0})&=& \frac{d}{dt}\left( \int_{\mathbb{R}^{d}}  \rho^{c} (\log \rho^{c} - \log \rho^{0})+  \frac{c^{2}}{2}\int_{\mathbb{R}^{d}}  \rho^{c}\|u^{c} - u^{0}\|^{2}\right)  \nonumber \\
 %&=& \frac{d}{dt}\left( \int_{\mathbb{R}^{d}}  \rho^{c} \log \rho^{c} + \frac{c^{2}}{2}\int_{\mathbb{R}^{d}}  \rho^{c} |u^{c}|^{2} - \int_{\mathbb{R}^{d}}  \rho^{c} \log \rho^{0}+  \frac{c^{2}}{2}\int \rho^{c}( |u^{0}|^{2} - 2u^{c} \cdot u^{0})\right)  \nonumber \\
% &\leq& - \int_{\mathbb{R}^{d}}  |u^{c}|^{2}\rho^{c} d\mu - \int_{\mathbb{R}^{d}}  \rho^{c} \langle u^{c}, \nabla \log \rho^{0}\rangle - \int_{\mathbb{R}^{d}}  \langle \nabla \rho^{c}, u^{0}\rangle + \int_{\mathbb{R}^{d}}  \rho^{c}\langle \nabla \log \rho^{0}, u^{0}\rangle +  \frac{c^{2}}{2}\frac{d}{dt}\int_{\mathbb{R}^{d}}  \rho^{c}( |u^{0}|^{2} - 2u \cdot u^{0})  \nonumber \\
%  &\leq& - \int_{\mathbb{R}^{d}}  \|u^{c}\|^{2}\rho^{c} d\mu  - \int_{\mathbb{R}^{d}}  \rho^{c} \langle u^{c}, \nabla \log \rho^{0}\rangle - \int_{\mathbb{R}^{d}}  \langle \nabla \rho^{c}, u^{0}\rangle + \int_{\mathbb{R}^{d}}  \rho^{c}\langle \nabla \log \rho^{0}, u^{0}\rangle  \nonumber  \\
%  & & +  c^{2}\big( \int_{\mathbb{R}^{d}}  \rho^{c}\langle u^{c}, (u^{0} - u^{c})\cdot  \nabla u^{0}\rangle  + \int_{\mathbb{R}^{d}}  \rho^{c}( (u^{0} - u^{c}) \cdot  (\Delta u^{0} - 2\nabla u^{0} \cdot u^{0}) \\
%  & & + c^{-2}(u^{c} + \nabla \log \rho) \cdot u^{0})\big)  \nonumber \\
 % &\leq& - \int_{\mathbb{R}^{d}}  \rho^{c}\|u^{c} - u^{0}\|^{2}  +  c^{2}\int_{\mathbb{R}^{d}}  \rho^{c}\langle u^{c}, (u^{0} - u^{c})\cdot  \nabla u^{0}\rangle  + c^{2}\int_{\mathbb{R}^{d}}  \rho (u^{0} - u^{c}) \cdot  (\Delta u^{0} - 2\nabla u^{0} \cdot u^{0}) \nonumber \\
  &\leq&- \int_{\mathbb{R}^{d}}  \rho^{c}\|u^{c} - u^{0}\|^{2}   -  c^{2}\int_{\mathbb{R}^{d}}  \rho^{c}\langle u^{c} - u^{0}, (u^{c} - u^{0}) \cdot  \nabla u^{0}\rangle \\
  & &\hskip1cm  + c^{2}\int_{\mathbb{R}^{d}}  \rho^{c} (u^{0} - u^{c}) \cdot  (\Delta u^{0} - u^{0} \cdot \nabla u^{0}).   \label{eeee2}
   \end{eqnarray*}
%  That is to say 
% \begin{eqnarray}
%& & \frac{d}{dt}\eta(\rho^{c}, m^{c}|\rho^{0}, m^{0}) \\
%&\leq&- \int_{\mathbb{R}^{d}}  \rho^{c}\|u^{c} - u^{0}\|^{2}   -  c^{2}\int_{\mathbb{R}^{d}}  \rho^{c}\langle u^{c} - u^{0}, (u^{c} - u^{0}) \cdot  \nabla u^{0}\rangle  + c^{2}\int_{\mathbb{R}^{d}}  \rho (u^{0} - u^{c}) \cdot  (\Delta u^{0} - u^{0} \cdot \nabla u^{0}).   \nonumber\\
% \label{eeee2}
%  \end{eqnarray}
   
%In the case of entropy solution,  we may replace the equality in $(\ref{eeee1})$ by the entropy inequality
%
% \begin{eqnarray}
% \frac{d}{dt}\left( \int_{\mathbb{R}^{d}}  \rho \log \rho + \frac{c^{2}}{2}\int_{\mathbb{R}^{d}}  \rho |u|^{2}\right)\leq  \int_{\mathbb{R}^{d}}  \left( - \frac{1}{c}\nabla_{x} \cdot q(\rho, m)  - |u|^{2}\rho\right)d\mu= - \int_{\mathbb{R}^{d}}  |u|^{2}\rho d\mu. \label{eeee3}
% \end{eqnarray}
%Then we can prove that the inequality  $(\ref{eeee2})$ still holds for entropy solution. 
%

By the Cauchy-Schwarz inequality, we have
\begin{eqnarray*}
(u^0-u^{c})\cdot  (\Delta u^{0} - u^{0} \cdot \nabla u^{0})\leq \frac{1}{2}(\|u^0-u^{c}\|^2+\|\Delta u^{0} - u^{0} \cdot \nabla u^{0}\|^2).
\end{eqnarray*} 
Together with the fact that $  \int_{\mathbb{R}^{d}} \rho^{c}  (\log \rho^{c}  - \log \rho^{0}) \geq 0$ we derive that
\beqnas
& &\frac{d}{dt}\eta(\rho^{c}, m^{c}|\rho^{0}, m^{0})\\
    %&\leq& c^{2}\|\nabla u^{0}\|_{\infty}\int_{\mathbb{R}^{d}}  \rho| u - u^{0}|^{2}  +  2c^{2}\int_{\mathbb{R}^{d}}  \rho |u^{0} - u|^2 + 2c^{2}\int_{\mathbb{R}^{d}}  \rho |\Delta u^{0} - u^{0} \cdot \nabla u^{0}|^{2}\\
 &\leq&  (4+2\|\nabla u^{0}\|_{\infty})\left(\frac{c^{2}}{2}\int_{\mathbb{R}^{d}}  \rho^{c} \| u^{c} - u^{0}\|^{2}  +  \int_{\mathbb{R}^{d}}  \rho^{c}  (\log \rho^{c}  - \log \rho^{0})\right)+2 c^{2}\|\Delta u^{0} - u^{0} \cdot \nabla u^{0}\|_{\infty}^{2},
 \eeqnas
 i.e.
%Thus 
% \beqnas
% \frac{d}{dt}\eta(\rho, m|\rho^{0}, m^{0}) \leq  (4+2\|\nabla u^{0}\|_{\infty})\eta(\rho, m|\rho^{0}, m^{0})+2 c^{2}\|\Delta u^{0} - u^{0} \cdot \nabla u^{0}\|_{\infty}^{2}.
% \eeqnas
%which yields
\beqnas
& &\eta(\rho^{c}, m^{c}|\rho^{0}, m^{0})(t) - \eta(\rho^{c}, m^{c}|\rho^{0}, m^{0})(0) \\
& & \leq \int^{t}_{0}(4+2\|\nabla u^{0}\|_{\infty}) \eta(\rho^{c}, m^{c}|\rho^{0}, m^{0})(s)ds + 2c^{2}\int^{t}_{0} \|\Delta u^{0} - u^{0} \cdot \nabla u^{0}\|_{\infty}^{2} (s)ds.
\eeqnas

Therefore by Gronwall's inequality, we conclude that
  \beqnas
& & \eta(\rho^{c}, m^{c}|\rho^{0}, m^{0})(t) \\%&=& \int \rho (\log \rho - \log \rho^{0})+  \frac{c^{2}}{2}\int \rho|u - u^{0}|^{2} \\
 &\leq& \left( \eta(\rho^{c}, m^{c}|\rho^{0}, m^{0})(0) + 2 c^{2}\int^{t}_{0} \|\Delta u^{0} - u^{0} \cdot \nabla u^{0}\|_{\infty}^{2} (s)ds\right)e^{\int^{t}_{0}(4+2\|\nabla u^{0}\|_{\infty})ds}.  
 \eeqnas
 \epf

 In our case, $\eta(\rho^{c}, m^{c}|\rho^{0}, m^{0})(0) = 0$, therefore \eqref{entropyineq} implies that when $c \rightarrow 0$, $h(\rho^{c}|\rho^{0}) \rightarrow 0$. Applying classical Csisz\'ar-Kullback inequality, we obtain
 $$
 \int_{\mathbb{R}^{d}}  |\rho^{c} - \rho_{0}|dx \leq h(\rho^{c}|\rho^{0}) \rightarrow 0. 
 $$
 thus $\rho^{c}$ converges to $\rho_{0}$ in $L^{1}([0, T] \times \bbR^{d})$. 

\bthm
Let $\rho_{0}, u_{0} \in H^{s}(\bbR^{d})$ be the initial value, and there exists local solution $(\rho, u)$ in $C([0, T], H^{s}) \in C^{1}([0, T], H^{s-1})$ to \eqref{euler.eq1} and \eqref{euler.eq2}. Then as  $c \rightarrow 0$, we have
$$
\sup_{t \in [0, T]}  \|\rho - \rho_{0}\|_{L^{1}}  \rightarrow 0.
$$
\ethm

% 
% Recall the log-Sobolev inequality for the heat semigroup $P_t=e^{t\Delta}$ on $\bbR^{d}$ is the following (Bakry-Gentil-Ledoux 5.5.1-5.5.2)
% \beq\label{LSIheatsemigroup}
% P_{t}(f\log f) - P_{t}f\log P_{t}f \leq tP_{t}\left(\frac{\Gamma(f)}{f}\right
% ),
% \eeq
% and with a $CD(\rho, \infty)$ condition we have
%  $$
% P_{t}(f\log f) - P_{t}f\log P_{t}f \leq \frac{1 - e^{-2\rho t}}{2\rho}P_{t}\left(\frac{\Gamma(f)}{f}\right)
% $$
% holds for diffusion semigroup $P_{t}$. 
% 
% 
% According to Villani \cite{V1}, the Talagrand inequality reads:  there exists a constant $\lambda>0$ such that for any $\rho \in \cP(\bbR^{n})$, 
% $$
% W^{2}_{2}(\rho, \nu) \leq \frac{2}{\lambda}H(\rho | \nu). 
% $$
% An LSI with  $\lambda$  is given by 
% $$
% H(\rho | \nu) \leq \frac{1}{2\lambda} I(\rho | \nu), 
% $$
% where $I(\rho | \nu) = \int \frac{\Gamma(h)}{h} d\nu$, where $h = \frac{\rho}{\nu} $. 
% 
% So the above results can not lead to the inequality we want. In fact, if it holds, then at $t = 0$, we should have
% $$
% W^{2}_{2}(\mu, f(x)dx) \leq CH(\mu |  f(x)dx),
% $$
% and this force a condition on the initial value. 

\subsubsection{The convergence of strong solutions with small initial value}
Inspired by the work of Coulombel-Goudon~\cite{CG}, we show the convergence of strong solutions in this part. Indeed, when $c \rightarrow 0$,  the convergence of \eqref{euler.eq1}, \eqref{euler.eq2} is equivalent to the strong relaxation limit of the isothermal Euler equations. Recall that Coulombel and Goudon~\cite{CG}  consider the isothermal Euler equations
\beqna
\label{isothermal.euler1}\partial_t \varrho^{c} + \nabla \cdot m^{c} &=& 0, \\
\label{isothermal.euler2}\partial_t m^{c} + \nabla \cdot (\frac{m^{c} \otimes m^{c}}{\varrho^{c}}) + a^{2}\nabla \varrho^{c} &=&  -\frac{1}{c}m^{c},
\eeqna
with initial values 
$$
\varrho^{c}(0, x) = \varrho_{0}(x), \ \ \ m^{c}(0, x) = m_{0}(x), 
$$
where $\varrho^{c}: \bbR^{+} \times \bbR^{d} \rightarrow (0, +\infty)$ is the density, $m^{c}: \bbR^{+} \times \bbR^{d} \rightarrow \bbR^{d}$ is the momentum, $a > 0$ is the speed of sound and $0 < c < 1$ is the relaxation time. We may write $m^{c} = \varrho^{c}v^{c}$, where $v^{c}$ is the velocity. 
According to the proof in \cite{CG}, the value of $a$ has no impact on the results, thus we take $a = 1$ in what follows.

 Now define
\beq
\rho^{c}(t, x) = \varrho^{c}({t\over c}, x), \ \ \ \ u^{c}(t, x) = \frac{1}{c} v^{c}({t\over c}, x), 
\eeq
then  $\rho^{c}$ and $u^{c}$ satisfy
\begin{eqnarray}
\partial_t \rho^{c} + \nabla \cdot (\rho^{c} u^{c}) &=& 0, \label{PQ1}\\
\partial_t (\rho^{c} u^{c}) + \nabla \cdot (\rho^{c} u^{c} \otimes u^{c}) + \frac{1}{c^{2}}(\nabla \rho^{c} + \rho^{c} u^{c}) &=& 0, \label{PQ2}
\end{eqnarray}
which are exactly Equations \eqref{euler.eq1}, \eqref{euler.eq2}.

%More precisely, they proved the following 
%\bthm[Coulombel-Goudon2014]\label{CGestimate} Let $\bar{\varrho} > 0$, and let $s \in \bbN$ with $s > {d \over 2} +1$. There exist two constants $\delta > 0$ and $C> 0$ such that, for all $c \in (0, 1]$ and for all initial data $(\varrho_{0}, m_{0})$ verifying $\|\varrho_{0} - \bar{\varrho}\|_{H^{s}(\bbR^{d})} + \|m_{0}\|_{H^{s}(\bbR^{d})} \leq \delta$, there exists a unique global solution $(\varrho^{c}, m^{c})$ to \eqref{isothermal.euler1}, \eqref{isothermal.euler2} such that $(\varrho^{c} - \bar{\varrho}, m^{c}) \in \cC(\bbR^{+}, H^{s}(\bbR^{d}))$. Moreover, this solution satisfies
%\beqna\label{estimate}
%& &\sup_{t \geq 0}(\|\varrho^{c} - \bar{\varrho}\|^{2}_{H^{s}(\bbR^{d})} + \|m^{c}\|^{2}_{H^{s}(\bbR^{d})}) + \frac{1}{c}\int^{\infty}_{0}\|m^{c}\|^{2}_{H^{s}(\bbR^{d})}ds\\
%\nonumber & &\hskip2cm  \leq C(\|\varrho_{0} - \bar{\varrho}\|^{2}_{H^{s}(\bbR^{d})} + \|m_{0}\|^{2}_{H^{s}(\bbR^{d})}),
%\eeqna
%where the constant $C$ is independent of $c$.
%\ethm
%Direct computations lead to
%\beqnas
%\|\varrho^{c} - \bar{\varrho}\|^{2}_{H^{s}(\bbR^{d})} &=& \|\rho^{c} - \bar{\varrho}\|^{2}_{H^{s}(\bbR^{d})},\\
%\|v^{c}\|^{2}_{H^{s}(\bbR^{d})} &=& c^{2}\|u^{c}\|^{2}_{H^{s}(\bbR^{d})},\\
%\|\varrho^{c} v^{c}\|^{2}_{H^{s}(\bbR^{d})} &=& c^{2}\|\rho^{c} u^{c}\|^{2}_{H^{s}(\bbR^{d})},\\
%\frac{1}{c}\int^{\infty}_{0}\|\varrho^{c}v^{c}\|^{2}_{H^{s}(\bbR^{d})}ds &=& \int^{\infty}_{0}\|\rho^{c}u^{c}\|^{2}_{H^{s}(\bbR^{d})}ds.
%\eeqnas

In \cite{CG}, the authors proved that for small initial data, Equations \eqref{isothermal.euler1}, \eqref{isothermal.euler2} admit a global solution, which satisfies a uniformly energy estimate. With the above scaling procedure, their uniform estimate leads to the following result for $(\rho^{c}, u^{c})$.
\bthm\label{thm.estimate1}
Let $\bar{\rho} > 0$, and $s \in \bbN$ with $s > {d \over 2} +1$. Then there exists two constants $\delta > 0$ and $C> 0$ such that for any $c \in (0, 1)$ and for initial data $(\rho_{0}, u_{0})$ satifying $\|\rho_{0} - \bar{\rho}\|_{H^{s}(\bbR^{d})} + \|\rho_{0}u_{0}\|_{H^{s}(\bbR^{d})} \leq \delta$,  there exists a unique global solution $(\rho^{c}, u^{c})$ to \eqref{euler.eq1}, \eqref{euler.eq2}, such that $(\rho^{c} - \bar{\rho}, u^{c}) \in \cC(\bbR^{+}, H^{s}(\bbR^{d}))$ satisfies
\beqna\label{estimate111}
& &\sup_{t \geq 0}(\|\rho^{c} - \bar{\rho}\|^{2}_{H^{s}(\bbR^{d})} +c^{2}\|\rho^{c}u^{c}\|^{2}_{H^{s}(\bbR^{d})}) + \int^{\infty}_{0}\|\rho^{c}u^{c}\|^{2}_{H^{s}(\bbR^{d})}ds\\
\nonumber & & \hskip2cm \leq C(\|\rho_{0} - \bar{\rho}\|^{2}_{H^{s}(\bbR^{d})} + \|\rho_{0}u_{0}\|^{2}_{H^{s}(\bbR^{d})}).
\eeqna
\ethm

By the same argument as used in~\cite{CG}, we can prove the following convergence result for the strong solution with small initial value. To save the length of the paper, we omit the proof. 
\bthm
Let $(\rho^{c}, u^{c})$ be the solution to Equations \eqref{euler.eq1} and \eqref{euler.eq2} with initial data $(\rho_{0}, u_{0})$. Let $B_{r}$ be the ball in $\bbR^{d}$ centered at the origin of radius $r$. Under the same assumption as in Theorem \ref{thm.estimate1}, for any $0 < T, R < \infty$, $0 < s' < s$,  $\rho^{c}$ converges in $\cC([0, T], H^{s'}(B_{R}))$ towards $\rho^{0} \in \cC(\bbR^{+}, \bar{\rho} + H^{s}(\bbR^{d}))$, which is the unique solution to the Cauchy problem of the heat equation
\beqnas
\partial_t \rho^{0} &=& \Delta \rho^{0},\\
\rho^{0}(0, x) &=& \rho_{0}(x).
\eeqnas
\ethm

\subsection{The convergence of the Langevin deformation of flows on Wasserstein space: $c \rightarrow \infty$}
In this part we turn to the case when $c \rightarrow \infty$. First rewrite \eqref{euler.eq1}, \eqref{euler.eq2} into the following symmetric hyperbolic equations
\beqna
\label{hyper1}\partial_t \log \rho^{c} + u^{c} \cdot \nabla \log \rho^{c} + \nabla \cdot u^{c}&=& 0, \\
\label{hyper2}c^{2}\partial_t u^{c} + c^{2}u^{c} \cdot \nabla u^{c} + \nabla \log \rho^{c} + u^{c} &=& 0.
\eeqna
Let $U^{c} = \begin{pmatrix}
\log \rho^{c} \\
u^{c}
\end{pmatrix}$, then the above equations turn into
%\beqnas
%\begin{pmatrix}
%1 & 0 \\
%0     & c^{2}\Id_{n}
%\end{pmatrix}\partial _t \begin{pmatrix}
%\log \rho \\
%u
%\end{pmatrix} +  \begin{pmatrix}
%u \cdot \nabla & \nabla \cdot \\
%\nabla                 & c^{2}u \cdot \nabla
%\end{pmatrix}\begin{pmatrix}
%\log \rho \\
%u
%\end{pmatrix} +  \begin{pmatrix}
%0 & 0 \\
%0     & \Id_{n}
%\end{pmatrix}\begin{pmatrix}
%\log \rho \\
%u
%\end{pmatrix}  = 0.
%\eeqnas
\beq\label{eq2.2}
A_{0}(c)\partial_t U^{c} + \sum^{n}_{j = 1}A_{j}(U^{c}, c)\partial_j U^{c} + B(0)U^{c} = 0,
\eeq
where $A_{0}(c) = \begin{pmatrix}
1 & 0 \\
0     & c^{2}\Id_{d}
\end{pmatrix}$, $B(0) = \begin{pmatrix}
0 & 0\\
0     & \Id_{d}
\end{pmatrix}$, 
and 
$$A_{j}(U, c) =  \begin{pmatrix}
u^{j}  &  0        & \cdots & 1 & \cdots & 0 \\
0        & c^{2}u^{j}   &   \cdots  &   0  & \cdots & 0 \\
\vdots & \vdots & \ddots  &   0 & \cdots & 0 \\
1 &    0       &         &          &          &                   \\
\vdots &  \vdots  &         &    \vdots       &   \ddots        &    \vdots       \\
0 & 0 &  \cdots &   0 & \cdots & c^{2}u^{j} \\
\end{pmatrix}_{(d+1) \times (d+1)}.$$

%\beqnas
%(u \cdot \nabla u)^{i} = \sum^{n}_{j=1}u^{j}\nabla_{j}u^{i}.
%\eeqnas
%
%And the limit equation is 
%\beqnas
%\begin{pmatrix}
%1 & 0 \\
%0     & 0
%\end{pmatrix}\partial _t \begin{pmatrix}
%\log \rho \\
%u
%\end{pmatrix} +  \begin{pmatrix}
%u \cdot \nabla & \nabla \cdot \\
%\nabla                 & 0
%\end{pmatrix}\begin{pmatrix}
%\log \rho \\
%u
%\end{pmatrix} +  \begin{pmatrix}
%0 & 0 \\
%0     & \Id_{n}
%\end{pmatrix}\begin{pmatrix}
%\log \rho \\
%u
%\end{pmatrix}  = 0.
%\eeqnas

Following the proof in Klainerman-Majda \cite{KM1, KM2, Maj}, we obtain the convergence result.
\bthm\label{convergence2}
Let $s \in \bbN$ with $s > {d \over 2} +1$. Let $(\rho^{c}, u^{c})$ be the solution to \eqref{euler.eq1},  \eqref{euler.eq2}, and there exists a constant $M > 0$, such that the initial data $(\rho_{0}, u_{0})$ satisfies
$$
\|\log \rho_{0}\|_{H^{s}}^{2} + \|u_{0}\|_{H^{s}}^{2} \leq M.
$$
Then there exists $T > 0$, which is independent of $c>1$, such that when $c \rightarrow \infty$,  $\rho^{c}dx$ weakly converges to $\rho^{\infty}dx$ in $\cC([0, T], \cP(\bbR^{d}))$,  and $u^{c}$ converges to $u^{\infty}(x, t)\in \cC([0, T], H^{s}) \cap \cC^{1}([0, T], H^{s-1})$ in $\cC([0, T], H^{s'-1}(B_{R}))$ for any $R > 0$, $s' < s$. Moreover,  $(\rho^{\infty}, u^{\infty})$ satisfies
\beqna\label{geodesic2}
\label{geoeq1}\partial_t \rho^{\infty} + \nabla \cdot (\rho^{\infty}u^{\infty}) &=& 0,\\
\label{geoeq2}\partial_t u^{\infty}  + u^{\infty}  \cdot \nabla u^{\infty}  &=& 0, \\
\nonumber \rho^{\infty}(0, x) = \rho_{0}(x), \ \ \ u^{\infty}(0, x) &=& u_{0}(x).
\eeqna
\ethm

%\brmq
%When $n \rightarrow \infty$, we say a family of probability measure $\mu_{n}$ on the sample space $X$  weakly converges to the probability measure $\mu$, if
%\beqnas
%\lim_{n \rightarrow \infty}\int_{X}f(x)d\mu_{n}(x) = \int_{X}f(x)d\mu(x),
%\eeqnas
%holds for all the real continuous, bounded function $f(x)$ on $X$. 
%\ermq
\brmq
 The limit equation \eqref{geoeq2} is a pressureless Euler equation, and  the equations \eqref{geoeq1}, \eqref{geoeq2} are called sticky particle system. There are a lot of research regarding this topic, see \cite{BreGre}, \cite{ERS} for the one-dimensional case, and \cite{sever} for the existence of measure solution in high dimension. As pointed out in \cite{BreGre}, this system can formally be obtained through the compressible Euler equation with the pressure term approaching zero or the Boltzmann equation with the temperature going to zero. In \cite{chen-liu} the authors rigorously proved that the entropy solution to the compressible Euler equation converges to that to the pressureless Euler equation when the  the pressure term goes to zero in one dimension. 
\ermq

%
%\brmq
%In the later section we will prove that the equation \eqref{geoeq2} can be derived from Hamilton-Jacobi equation
%$$
%\partial_{t}\phi + \frac{1}{2}|\nabla \phi|^{2} = 0
%$$
%by taking gradient operator on both sides. The solution to this Hamilton-Jacobi equation can be obtained by the Hopf-Lax formula,  i.e. for any $t > 0$, a.s. $x \in \bbR^{d}$, 
%$$
%\phi(t, x) = \min_{y \in \bbR^{d}}\left\{ \frac{|x - y|^{2}}{2t} + \phi_{0}(y) \right\}.
%$$
%One can prove that the solution $\phi(t, x)$ is Lipschitz continuous.  According to Theorem 12.1 in P.L. Lions \cite{Lions},  when the initial value satisfies $\phi_{0} \in W^{2, \infty}(\bbR^{d})$, and $\|D\phi_{0}\|_{L^{\infty}} < \frac{1}{T}$, we have $\phi \in W^{2, \infty}([0, T] \times \bbR^{d})$.  
%Therefore shocks do not appear in the solution to the equation \eqref{geoeq2} on $[0, T] \times \bbR^{d}$. 
%\ermq
%

\brmq
As pointed out by Villani in \cite{V1}, shocks do not arise in the case of optimal transportation. One could prove that for the geodesic flow, the velocity field $u_{t}$ is locally Lipschitz in both $t$ and $x$, see Theorem 5.51 in \cite{V1}.
\ermq

According to Klainerman-Majda \cite{KM1, KM2, Maj}, to prove Theorem \ref{convergence2}, we need 
first to prove the following uniform a priori estimate.
\bprop\label{uniformestimate}
Assume that for initial value $(\log \rho_{0}, u_{0})$, there exist two constants $M >0$ and $s > {d \over 2} +1$, such that
\beqna\label{initial}
\|\log \rho_{0}\|_{H^{s}}^{2} + \|u_{0}\|_{H^{s}}^{2} \leq M.
\eeqna
%{\red this is in fact 
%\beqnas\label{initial}
%c^{-2}\|\log \rho_{0}\|_{H^{s}}^{2} + \|u_{0}\|_{H^{s}}^{2} \leq M.
%\eeqnas
%by taking $c =1$.
%}
Then there exists $T> 0$, which is independent of $c$, such that for any $c \geq 1$,  the equation \eqref{eq2.2} admits classical solution $\log \rho^{c}, u^{c} \in \cC([0, T], H^{s}) \cap \cC^{1}([0, T], H^{s-1})$, satisfying
\beqna\label{uniformestimate1}
\sup_{[0, T]} c^{-2}\|\log \rho^{c}\|_{H^{s}}^{2} + \|u^{c}\|_{H^{s}}^{2} \leq R,
\eeqna
where $R$ is a constant independent of $c$.
\eprop

To see this, we first prove the following a priori estimate as in \cite{Maj} .
\blem\label{lemma}
Assume that the conditions in Theorem \ref{uniformestimate} hold, and the equation \eqref{eq2.2} admits local classical solution $(\log \rho^{c}, u^{c})$ on  $[0, T^{*}_{c}]$, satisfying   
\beqna\label{assumpM}
\sup_{[0, T_{c}^{*}]}  c^{-2}\|\log \rho^{c}\|_{H^{s}}^{2} + \|u^{c}\|_{H^{s}}^{2} \leq 2M,
\eeqna
for all $c \geq 1$. Then there exists a constant $C>0$, which is independent of $c$, such that
\beqna\label{estimate11}
\sup_{[0, T_{c}^{*}]}  c^{-2}\|\log \rho^{c}\|_{H^{s}}^{2} + \|u^{c}\|_{H^{s}}^{2} \leq (c^{-2}\|\log \rho_{0}\|^{2}_{H^{s}} + \|u_{0}\|^{2}_{H^{s}})e^{C(2M+ 1)T_{c}^{*}}.
\eeqna

\elem
%{\red An explanation is about iteration scheme. Hyperbolic equation of this kind satisfies the fixed point theorem. not enough\\
%integration by parts formula: can be checked with compact supported functions first. } 
%
\bpf
We apply the classical energy estimate method.  First, multiply $U$ on both sides of the equation \eqref{eq2.2} and integrate on $[0, t]$, then we obtain
\beqnas
\langle A_{0}(c)U, U\rangle \big{|}^{t}_{0} -\sum_{j}\int^{t}_{0}\langle \partial_{j}A_{j} U, U\rangle + 2\int^{t}_{0}\langle B(0)U, U\rangle = 0.
\eeqnas
Thus, 
 \beqnas
& &\|\log \rho^{c}\|_{L^{2}}^{2}(t) + c^{2}\|u^{c}\|_{L^{2}}^{2}(t) \\
%&=&  \|\log \rho_{0}\|_{L^{2}}^{2} + c^{2}\|u_{0}\|_{L^{2}}^{2} + \sum_{j}\int^{t}_{0}\langle \partial_{j}A_{j} U, U\rangle - 2\int^{t}_{0}\langle B(0)U, U\rangle\\
%&\leq&  \|\log \rho_{0}\|_{L^{2}}^{2} + c^{2}\|u_{0}\|_{L^{2}}^{2} + \sum_{j}\big(\int^{t}_{0}\langle \partial_{j}u^{c}_{j} \log \rho^{c}, \log \rho^{c} \rangle + c^{2}\sum_{i}\int^{t}_{0}\langle \partial_{j}u^{c}_{j} u^{c}_{i}, u^{c}_{i} \rangle \big) \\
&\leq&  \|\log \rho_{0}\|_{L^{2}}^{2} + c^{2}\|u_{0}\|_{L^{2}}^{2} + \int^{t}_{0}\|\partial_j u^{c}_{j}\|_{L^{\infty}}(\|\log \rho^{c}\|_{L^{2}}^{2}(s) + c^{2}\|u^{c}\|_{L^{2}}^{2}(s))ds.
 \eeqnas
 
By Sobolev inequality, for $s > {d \over 2}+1$ we have
$$
\|\nabla u\|_{L^{\infty}} < C_{s}\| u\|_{H^{s}},
$$
which leads to
\beqna\label{estimate0}
c^{-2}\|\log \rho^{c}\|_{L^{2}}^{2}(t) + \|u^{c}\|_{L^{2}}^{2}(t)&\leq&  c^{-2}\|\log \rho_{0}\|_{L^{2}}^{2} + \|u_{0}\|_{L^{2}}^{2} \\
& &+ C_{s} \int^{t}_{0}\| u^{c}\|_{H^{s}}(c^{-2}\|\log \rho^{c}\|_{L^{2}}^{2}(s) + \|u^{c}\|_{L^{2}}^{2}(s))ds. \nonumber 
\eeqna

Now applying differential operator $D^{\alpha}$ on both sides of the equation \eqref{eq2.2}, where $1 \leq \alpha \leq s$, we obtain
\beqnas
A_{0}(c)\partial_t D^{\alpha}U + \sum^{n}_{j = 1}A_{j}(u, c)\partial_j D^{\alpha}U + B(0)D^{\alpha}U = - \sum^{\alpha -1}_{\beta = 0}\sum^{n}_{j = 1}D^{\alpha - \beta}(A_{j}(u, c))D^{\beta}\partial_j U. 
\eeqnas
Again using the  energy estimate method, we derive that
%
%\beqnas
%& &\partial_t \langle A_{0}(c)D^{\alpha}U, D^{\alpha}U \rangle\\
% &=&2 \langle A_{0}(c)\partial_t D^{\alpha}U, D^{\alpha}U \rangle\\
%&\leq& \langle \sum^{n}_{j = 1}\partial_j A_{j}(u, c)D^{\alpha}U, D^{\alpha}U \rangle - 2\langle \sum^{n}_{j = 1}D^{\alpha}(A_{j}(u, c))\partial_j U, D^{\alpha}U \rangle\\
%&\leq&  \langle \sum^{n}_{j = 1}\partial_j A_{j}(u, c)D^{\alpha}U, D^{\alpha}U \rangle - 2\langle  \sum^{n}_{j = 1}D^{\alpha}(A_{j}(u, c))\partial_j U, D^{\alpha}U \rangle\\
%\eeqnas
%'

\begin{eqnarray}
& &\langle A_{0}(c)D^{\alpha}U, D^{\alpha}U \rangle(t) \label{estimate1}\\
&\leq&  \langle A_{0}(c)D^{\alpha}U, D^{\alpha}U \rangle(0) + \int^{t}_{0} \big(\langle \sum^{n}_{j = 1}\partial_j A_{j}(u, c)D^{\alpha}U, D^{\alpha}U \rangle(s) \nonumber \\
& &  \hskip2cm - 2\int^{t}_{0} \langle \sum^{\alpha -1}_{\beta = 0}\sum^{n}_{j = 1}D^{\alpha - \beta}(A_{j}(u, c))D^{\beta}\partial_j U, D^{\alpha}U \rangle(s) ds \nonumber\\
&\leq&  \langle A_{0}(c)D^{\alpha}U, D^{\alpha}U \rangle(0) + \int^{t}_{0}  \sum^{n}_{j = 1}\|\partial_j u^{c}_{j}\|_{L^{\infty}}(\|D^{\alpha}\log \rho^{c}\|_{L^{2}}^{2} + c^{2}\|D^{\alpha}u^{c}\|_{L^{2}}^{2})ds\nonumber \\
&  & +  \sum^{\alpha -1}_{\beta = 0}\sum^{n}_{j = 1}  \int^{t}_{0}  \big(   \|D^{\alpha - \beta}u^{c}_{j} \|_{L^{2}}^{2}(\|D^{\beta}\partial_j \log \rho^{c}\|_{L^{2}}^{2} +c^{2}\|D^{\beta}\partial_j u^{c}\|_{L^{2}}^{2}) \nonumber \\
& &  \hskip2cm  + ( \|D^{\alpha} \log \rho^{c} \|_{L^{2}}^{2}+ c^{2}\|D^{\alpha}u^{c}\|_{L^{2}}^{2}) \big) ds, \nonumber
\end{eqnarray}
where the second inequality is due to the Cauchy-Schwarz inequality 
\begin{eqnarray*}
\|ab\|_{L^1}\leq \|a\|_{L^2}\|b\|_{L^2},
\end{eqnarray*}
and also
\beqnas
& &- 2\int^{t}_{0} \langle \sum^{\alpha -1}_{\beta = 0}\sum^{n}_{j = 1}D^{\alpha - \beta}(A_{j}(u, c))D^{\beta}\partial_j U, D^{\alpha}U \rangle(s) ds\\
&=& - 2\int^{t}_{0}   \sum^{\alpha -1}_{\beta = 0}\sum^{n}_{j = 1}\big( \langle D^{\alpha - \beta}u^{c}_{j}D^{\beta}\partial_j \log \rho^{c}, D^{\alpha} \log \rho^{c} \rangle + c^{2}\langle D^{\alpha-\beta}u^{c}_{j}D^{\beta}\partial_j u^{c}, D^{\alpha}u^{c} \rangle \big) ds\\
%&\leq&  \int^{t}_{0}  \sum^{\alpha -1}_{\beta = 0}\sum^{n}_{j = 1} \big( \|D^{\beta}\partial_j \log \rho^{c}\|_{L^{2}}^{2}\|D^{\alpha - \beta}u^{c}_{j} \|_{L^{2}}^{2} + \|D^{\alpha} \log \rho^{c} \|_{L^{2}}^{2}\\
%& & \ \ \ \  \ \ \ \ \ \ \ \ \  + c^{2}\|D^{\beta}\partial_j u^{c}\|_{L^{2}}^{2}\|D^{\alpha - \beta}u^{c}_{j}\|_{L^{2}}^{2} + c^{2}\|D^{\alpha}u^{c}\|_{L^{2}}^{2} \big) ds\\
&\leq&  \int^{t}_{0}   \sum^{\alpha -1}_{\beta = 0}\sum^{n}_{j = 1} \big(   \|D^{\alpha - \beta}u^{c}_{j} \|_{L^{2}}^{2}(\|D^{\beta}\partial_j \log \rho^{c}\|_{L^{2}}^{2} +c^{2}\|D^{\beta}\partial_j u^{c}\|_{L^{2}}^{2}) + ( \|D^{\alpha} \log \rho^{c} \|_{L^{2}}^{2}+ c^{2}\|D^{\alpha}u^{c}\|_{L^{2}}^{2}) \big) ds.
\eeqnas
Therefore we have
\beqna
& &c^{-2}\|D^{\alpha}\log \rho^{c}\|_{L^{2}}^{2}(t) + \|D^{\alpha}u^{c}\|_{L^{2}}^{2}(t) \\
\nonumber &\leq& c^{-2}\|D^{\alpha}\log \rho_{0}\|_{L^{2}}^{2} + \|D^{\alpha}u_{0}\|_{L^{2}}^{2} +  \int^{t}_{0}  \sum^{n}_{j = 1}\|\partial_j u^{c}_{j}\|_{L^{\infty}}(c^{-2}\|D^{\alpha}\log \rho^{c}\|_{L^{2}}^{2} + \|D^{\alpha}u^{c}\|_{L^{2}}^{2})ds \\
\nonumber & & \ \ \ \ +  C  \int^{t}_{0}  \sum^{n}_{j = 1}\big(   \|u^{c}_{j} \|_{H^{s}}(c^{-2}\| \log \rho^{c}\|_{H^{s}}^{2} +\|u^{c}\|_{H^{s}}^{2}) + ( c^{-2}\|D^{\alpha} \log \rho^{c} \|_{L^{2}}^{2}+ \|D^{\alpha}u^{c}\|_{L^{2}}^{2}) \big) ds. \label{estimate2}
\eeqna

Combining \eqref{estimate1} and \eqref{estimate2}, then summing over $1 \leq \alpha \leq s$, we get
\beqnas
& &(c^{-2}\|\log \rho^{c}\|^{2}_{H^{s}} + \|u^{c}\|^{2}_{H^{s}})(t) \\
&\leq& (c^{-2}\|\log \rho_{0}\|^{2}_{H^{s}} + \|u_{0}\|^{2}_{H^{s}}) +  C_{1}\int^{t}_{0}\|u^{c}\|_{H^{s}}(c^{-2}\|\log \rho^{c}\|^{2}_{H^{s}} + \|u^{c}\|^{2}_{H^{s}})ds \\
& & + C_{2}\int^{t}_{0}(\|u^{c}\|^{2}_{H^{s}} + 1)(c^{-2}\|\log \rho^{c}\|^{2}_{H^{s}} + \|u^{c}\|^{2}_{H^{s}})ds,
\eeqnas
where $C_{1}, C_{2}$ only depend on Sobolev constants. By Gronwall's  inequality and the assumption \eqref{assumpM}, we deduce that
\beqnas
(c^{-2}\|\log \rho^{c}\|^{2}_{H^{s}} + \|u^{c}\|^{2}_{H^{s}})(t) \leq (c^{-2}\|\log \rho_{0}\|^{2}_{H^{s}} + \|u_{0}\|^{2}_{H^{s}})e^{C(M + 1)t},
\eeqnas
where the constant $C>0$ is independent of $c$. This finishes the proof of the lemma. 
\epf

Now we are ready to prove Proposition~\ref{uniformestimate}. 

\bpf
We choose $T$ such that
$$
e^{C(2M + 1)T} =2.
$$
Notice that under the initial condition \eqref{initial}, \eqref{estimate11} may directly lead to \eqref{assumpM},  where $T$ is independent of $c$. By taking $R=2M$ we finish the proof. 
\epf

Now we are in  the position to prove Theorem~\ref{convergence2}. Since the uniform estimate~\eqref{uniformestimate1} does not guarantee that the convergence of the equations, we need to prove uniform estimate on the time derivative. 

\bpf
We first prove the uniform estimate of time derivative. Set $\partial_t U_{c} = V_{c}$, and take derivative on both sides of  \eqref{eq2.2}, then we obtain
\beqnas
A_{0}(c)\partial_t V_{c} + \sum_{j}A_{j}(U, c)\partial_j V_{c} + B(0)V_{c} = - \sum_{j}\partial_tA_{j}(U, c)\partial_{j}U_{c}.
\eeqnas
Thus
\beqnas
\frac{1}{2}\partial_t \langle A_{0}(c)V_{c}, V_{c}\rangle = \frac{1}{2}\sum_{j}\langle \partial_{j}A_{j}V_{c}, V_{c}\rangle - \langle B(0)V_{c}, V_{c}\rangle - \langle \partial_tA_{j}(U, c)\partial_{j}U_{c}, V_{c}\rangle.
\eeqnas
By Sobolev inequality, we derive that
\beqnas
 & &\langle A_{0}(c)V_{c}, V_{c}\rangle(t) -  \langle A_{0}(c)V_{c}, V_{c}\rangle(0) \\
% &\leq& \int^{t}_{0} \big(\sum_{j}\|\partial_j u^{c}_{j}\|_{L^{\infty}}(\|\partial_t \log \rho^{c}\|_{L^{2}}^{2} + c^{2}\|\partial_t u^{c}\|_{L^{2}}^{2}) + \sum_j\langle \partial_t (u^{c})_{j}\partial_j \log \rho, \partial_t \log \rho \rangle\\
 %& & + c^{2}\sum_{ij}\langle \partial_t u^{c}_{j}\partial_j u^{c}_{i}, \partial_t u^{c}_{i}\rangle \big) ds\\
&\leq& C\int^{t}_{0}  \big(\|u^{c}\|_{H^{s}}(\|\partial_t \log \rho^{c}\|_{L^{2}}^{2} + c^{2}\|\partial_t u^{c}\|_{L^{2}}^{2}) + \|c^{-1}\nabla \log \rho^{c}\|_{L^{\infty}}(\|c\partial_t u^{c}\|_{L^{2}}^{2} + \|\partial_t \log \rho^{c}\|_{L^{2}}^{2}) ds\\
& & \hskip2cm+ C\int_0^t c^{2}\|\nabla u^{c}\|_{L^{\infty}}\|\partial_t u^{c}\|_{L^{2}}^{2} ds\\
%&\leq&  C\int^{t}_{0} \|u^{c}\|_{H^{s}}(\|\partial_t \log \rho^{c}\|_{L^{2}}^{2} + c^{2}\|\partial_t u^{c}\|_{L^{2}}^{2}) \\
%& & \hskip2cm + (\|c^{-1}\nabla \log \rho^{c}\|_{L^{\infty}} + \|\nabla u^{c}\|_{L^{\infty}})(c^{2}\|\partial_t u^{c}\|_{L^{2}}^{2} + \|\partial_t \log \rho^{c}\|_{L^{2}}^{2})\\
&\leq&  C\int^{t}_{0}  (\|c^{-1}\nabla \log \rho^{c}\|_{H^{s}} + \|\nabla u^{c}\|_{H^{s}})(\|\partial_t \log \rho^{c}\|_{L^{2}}^{2} + c^{2}\|\partial_t u^{c}\|_{L^{2}}^{2})ds.
\eeqnas
Hence 
\beqnas
& &(\|\partial_t \log \rho^{c}\|_{L^{2}}^{2} + c^{2}\|\partial_t u^{c}\|_{L^{2}}^{2})(t) \leq  (\|\partial_t \log \rho^{c}\|_{L^{2}}^{2} + c^{2}\|\partial_t u^{c}\|_{L^{2}}^{2})(0)\\
& & \hskip2cm + C\int^{t}_{0}(\|c^{-1}\nabla \log \rho^{c}\|_{H^{s}} + \|\nabla u^{c}\|_{H^{s}})(\|\partial_t \log \rho^{c}\|_{L^{2}}^{2} + c^{2}\|\partial_t u^{c}\|_{L^{2}}^{2})ds.
\eeqnas
Note that $\|c^{-1}\nabla \log \rho^{c}\|_{H^{s}} + \|\nabla u^{c}\|_{H^{s}}$ is uniformly bounded by a constant $R>0$ which is independent of $c$. Applying Gronwall's inequality, we obtain
\beqnas
(c^{-2}\|\partial_t \log \rho^{c}\|_{L^{2}}^{2} + \|\partial_t u^{c}\|_{L^{2}}^{2})(t)  \leq (c^{-2}\|\partial_t \log \rho^{c}\|_{L^{2}}^{2}(0) + \|\partial_t u^{c}\|_{L^{2}}^{2}(0))e^{CRt},
\eeqnas
which implies that
\beq\label{estimate4}
\sup_{[0, T]}\|\partial_t u^{c}\|_{L^{2}} \leq C,
\eeq
where $C$ is a constant which is independent of $c$. 

According to the uniform estimate~\eqref{uniformestimate1}, $u^{c}$ is uniformly bounded in $L^{\infty}([0, T], H^{s})$, then there exists $u^{\infty} \in L^{\infty}([0, T], H^{s})$ and a sequence  $\{u^{n}\}$ in $L^{\infty}([0, T], H^{s})$, such that $u^{n}$ weakly converges to $u^{\infty}$. 
%
%{\red The convergence holds in the sense of $\sigma(L^\infty, L^{1})$. See Cov 3.26 in H. Brezis. Notice that $L^{1}$ is separable while  $L^{\infty}$ is not.}
%

In view of Arzela-Ascoli theorem and \eqref{estimate4}, we get a subsequence $\{u^{c}\}$ of  $u^{n}$(still denoted by $\{u^{c}\}$) which converges to $u^{\infty}$ in $\cC([0, T], L^{2}(B_{R}))$.
% in the sense of strong topology. 

%{\red 
%Arzela-Ascoli Theorem: Let $X$ be a compact Hausdorff space and $Y$ a metric space.  Then $F \in C(X, Y)$ is COMPACT in the compact open topology if and only if it is equicontinuous, pointwise relatively compact and closed. 
%
%Now $X = [0, T]$, $Y = L^{2}(\bbR^{d})$. \eqref{estimate4} ensures that $u^{c}$ is equicontinuous in $\cC([0, T], L^{2}(\bbR^{d}))$. Now we need to check that any fixed $t \in [0, T]$, $u^{c}(t)$ is relatively compact (in the sense of strong topology in $L^{2}$) and closed in $L^{2}(\bbR^{d})$. 
%
%What we know is that $u^{c}(t)$ is uniformly bounded in  $L^{2}(\bbR^{d})$, which results in its relatively compactness in  $L^{2}(B_{R})$. 
%
%The closedness comes from the fact that $u^{c}(t)$ is uniformly bounded in $L^{2}$.}

By the Nash inequality, for $s' < s$ and any function $v \in H^{s}(\bbR^{d})$, we have
$$
\|v\|_{H^{s'}} \leq C_{s}\|v\|^{1 - \frac{s'}{s}}_{L^{2}}\|v\|^{{s' \over s}}_{H^{s}},
$$
which ensures that $u^{c}$ converges to $u^{\infty}$ in $\cC([0, T], H^{s'}(B_{R}))$.

Since $c^{-1}\nabla \log \rho^{c}$ and $u^{c}$ are uniformly bounded in $\cC([0, T], H^{s-1})$ (Proposition~\ref{uniformestimate}), we deduce that $\lim_{c \rightarrow \infty}c^{-2}(\nabla \log \rho^{c} + u^{c}) = 0$ holds in $\cC([0, T], H^{s})$. By the fact $u^c\in \cC([0, T], H^{s})$, we have 
$\nabla u^c\in \cC([0, T], H^{s-1})$. This yields {{$u^c\cdot \nabla u^c \in \cC([0, T], H^{s-1})$}} as 
$$
\|u^c\cdot \nabla u^c\|_{H^{s-1}} \leq C_s \|u^c\|_{H^{s-1}}\|\nabla u^c\|_{H^{s-1}}\leq C_s \|u^c\|_{H^{s}}\|\nabla u^c\|_{H^{s-1}}.
$$ 
Moreover,  Eq.~\eqref{hyper2} yields $\partial_t u^{c}=-u^{c} \cdot \nabla u^{c}- c^{-2}(\nabla \log \rho^{c} + u^{c})\in \cC([0, T], H^{s-1})$ and $\lim_{c \rightarrow \infty}(\partial_t u^{c} + u^{c} \cdot \nabla u^{c})  = -\lim_{c \rightarrow \infty}c^{-2}(\nabla \log \rho^{c} + u^{c})=
0$ in $\cC([0, T], H^{s-1})$. 

On the other hand, as a consequence of the Sobolev inequality: if $u, v \in H^{s}$, $s > \frac{d}{2}$, then $uv \in H^{s}$, and 
$$
\|uv\|_{H^{s}} \leq C_s \|u\|_{H^{s}}\|v\|_{H^{s}},
$$ 
we can derive that $u^{c} \cdot \nabla u^{c}$ converges to $u^{\infty} \cdot \nabla u^{\infty}$ in $\cC([0, T], H^{s'-1}(B_{R}))$. 

Moreover, by \eqref{estimate4} there exists $v^{\infty} \in L^{\infty}([0, T], L^{2})$ such that $\partial_t u^{c}$ weakly converges to $v^{\infty}$.  Also by the fact that $u^{c}$ converges to $u^{\infty}$ in $\cC([0, T], H^{s'}(B_{R}))$, we have $\partial_t u^{c}$ converges to $\partial_t u^{\infty}$ in distribution, implying that $\partial_t u^{\infty} = v^{\infty} \in  L^{\infty}([0, T], L^{2})$.
% convergence in distribution means that being coupled with functions in  $\cC([0, T] \times B_{R})$, while weak convergence in  $L^{\infty}([0, T], L^{2})$ means that being coupled with functions in  $\L^{1}([0, T] \times H^{s}(B_{R}))$. Since  $\cC([0, T] \times B_{R})$ is dense in  $\L^{1}([0, T] \times H^{s}(B_{R}))$, we deduce that $v^{\infty} = \partial_t u^{\infty}$.

As a result $\partial_t u^{\infty} + u^{\infty} \cdot \nabla u^{\infty} = 0$ holds in the weak-$*$ topology of $L^{\infty}([0, T], L^{2})$. Also by $u^{\infty} \cdot \nabla u^{\infty} \in L^{\infty}([0, T], H^{s-1})  \cap \cC([0, T], H^{s'-1}(B_{R}))$, we deduce that $\partial_t u^{\infty} \in  L^{\infty}([0, T], H^{s-1})  \cap \cC([0, T], H^{s'-1}(B_{R}))$, and thus $u^{\infty}  \in  Lip([0, T], H^{s-1}) \cap \cC^{1}([0, T],  H^{s'-1}(B_{R}))$. 
%i.e. $u^{\infty} \in Lip([0, T], H^{s-1}) \cap \cC^{1}([0, T] \times \bbR^{d})$ due to Sobolev embedding theorem.
%{\red  should take $\frac{d}{2} +1 < s' < s$.}
% The argument is that: for some generalized function, if its derivative is a continuous function, then it is itself a  differentiable  continuous function.  (do not know how to prove it yet)

Meanwhile, by Sobolev embedding theorem, for $s > s' > \frac{d}{2} +1$ we have $H^{s'-1}, H^{s-1} \subset C^{0}$, $H^{s'}, H^{s} \subset C^{1}$, which implies that ${u^{c}}, u^{\infty} \in \cC^{1}([0, T] \times \bbR^{d})$, and hence $u^{c} \in Lip([0, T] \times \bbR^{d})$ converges to $u^\infty \in Lip([0, T], H^{s-1}) \cap \cC^{1}([0, T] \times \bbR^{d})$ in  $\cC([0, T], H^{s'}(B_{R}))$ for any $R > 0$. 

Now we claim that for the solution $\rho^c$ to the continuity equation $\partial_t \rho^{c} + \nabla \cdot (\rho^{c}u^{c}) = 0$, $\rho^{c}dx$ weakly converges to $\rho^{\infty}dx$ in $\cC([0, T], \cP(\bbR^{d}))$, which satisfies \eqref{geodesic2}.

Since $u^{c}$ is Lipschitz in $t$ and $x$, according to Ambrosio-Gigli-Savar\'e~\cite{AGS} or Theorem 5.34 in Villani~\cite{V1}, we 
have $\rho^{c} = X^{c}_\sharp \rho_{0}$, where $X^{c}$ is the integral curve of $u^{c}$,   
%Villani Theorem 5.34, need u to be both lipschitz in t and x
\beqnas
\frac{d}{dt}X^{c}(t, x) &=& u^{c}(t, X^{c}(t, x)), \\
X(0, x) &=& x.
\eeqnas
Now we prove that as $c \rightarrow \infty$, $X^{c}$ uniformly converges to $X^{\infty}$ in $\cC([0, T] \times B_{R})$, for any $R > 0$ and $X^{\infty}$ is the flow generated by the vector field $u^{\infty}$ on $\mathbb{R}^d$. Now take $x \in B_{R}$, and 
%since
%\beqnas
%X^{c}(t, x) &=& x +  \int^{t}_{0}u^{c}(s, X^{c}(s, x))ds, \\
%X^{\infty}(t, x) &=& x +  \int^{t}_{0}u^{\infty}(s, X^{\infty}(s, x))ds, 
%\eeqnas
we derive that 
\beqnas
& &|X^{c}(t, x) - X^{\infty}(t, x) | \\
%&\leq& \int^{t}_{0}|u^{c}(s, X^{c}(s, x)) - u^{\infty}(s, X^{\infty}(s, x))|ds\\
&\leq&  \int^{t}_{0}|u^{c}(s, X^{c}(s, x)) - u^{\infty}(s, X^{c}(s, x))|ds + \int^{t}_{0}|u^{\infty}(s, X^{c}(s, x)) - u^{\infty}(s, X^{\infty}(s, x))|ds\\
&\leq&  \int^{t}_{0}|u^{c}(s, X^{c}(s, x)) - u^{\infty}(s, X^{c}(s, x))|ds + L\int^{t}_{0}|X^{c}(s, x) - X^{\infty}(s, x)|ds,
\eeqnas
where $L$ is the Lipschitz constant of $u^{\infty}$ in $x$ direction on $B_{R}$. Thus
\beqnas
\sup_{t \in [0, T], x \in B_{R}}|X^{c}(t, x) - X^{\infty}(t, x) | &\leq& \sup_{t \in [0, T], x \in  B_{R}} e^{Lt} \int^{t}_{0}|u^{c}(s, X^{c}(s, x)) - u^{\infty}(s, X^{c}(s, x))|ds \\
%&\leq& Te^{LT} \sup_{s \in [0, T], x \in  B_{R}} |u^{c}(s, X^{c}(s, x)) - u^{\infty}(s, X^{c}(s, x))|\\
&\leq& CTe^{LT} \sup_{s \in [0, T]}\|u^{c} - u^{\infty}\|_{H^{s'}( B_{R})} \rightarrow 0.
\eeqnas

Then the weak convergence of  $\rho^{c}(x)$ to $\rho^{\infty} = X^{\infty}_\sharp \rho_{0}$ in $\cC([0, T], \cP(\bbR^{d}))$ is a direct consequence of Lemma 5.2.1 in Ambrosio-Gigli-Savar\'e~\cite{AGS}. 

%which we quote here for the completeness of this paper. 
%\blem
%Let $X_{1}$, $X_{2}$ be two seperable metric spaces,  and $r: X_{1} \rightarrow X_{2}$ is a Borel map. Let $r_{n}: X_{1} \rightarrow X_{2}$ be a family of Borel maps, which uniformly converges to $r$ on compact sets in $X_{1}$ .  Also let $(\mu_{n})\in \cP(X_{1})$ be a compact sequence, which weakly converges to $\mu \in \cP(X_{1})$. Then if $r$ is continuous, we have $(r_{n}) \sharp \mu_{n}$ weakly converges to $r \sharp \mu$. 
%\elem
%
\epf

\subsection{The convergence of the Langevin deformation of flows on compact manifolds}
In this section we give the proof of Theorem~\ref{main.convergence}.

\bpf
It suffices to prove that when $u^{\infty}(0, x) = \nabla \phi_{0}(x)$, there exits a function $\phi^{\infty}$, such that $u^{\infty} = \nabla \phi^{\infty}$. Following the proof of Theorem~\ref{existence}, we construct $\phi^{c}$ by $u^{c}$. More precisely, define
\beqnas
\phi^{c}(t, x) = \phi_{0}(t, x) -  e^{-\frac{1}{c^{2}}t}\int^{t}_{0}  e^{\frac{1}{c^{2}}s}\big(\frac{1}{c^{2}}\nabla \log \rho^{c}(s, x)  + \frac{1}{2}|u^{c}(s, x)|^{2}\big)ds.
\eeqnas
Following the proof of~\ref{convergence2}, we derive that $c^{-1}\nabla \log \rho^{c}$ is uniformly bounded in $\cC([0, T], H^{s-1})$, thus $e^{-\frac{1}{c^{2}}t}\int^{t}_{0}  \frac{1}{c^{2}}e^{\frac{1}{c^{2}}s}\nabla \log \rho^{c}(s, x) ds$ converges to $0$ in $\cC([0, T], H^{s-1})$. Since $u^{c}$ converges in \\
$\cC([0, T], H^{s'}(B_{R}))$ for any $R > 0$, $s' < s$, we deduce that the convergence also holds in $\cC([0, T], \cC^{1}(B_{R}))$ according to Sobolev embedding theorem. Therefore we get the convergence of $ e^{-\frac{1}{c^{2}}t}\int^{t}_{0}  e^{\frac{1}{c^{2}}s}|u^{c}(s, x)|^{2}ds$ to $\int^{t}_{0}|u^{\infty}(s, x)|^{2}ds$ in $\cC([0, T], \cC^{1}(B_{R}))$. Now define
$$
\phi^{\infty}(t, x) = \phi_{0}(t, x) - \frac{1}{2} \int^{t}_{0} |u^{\infty}(s, x)|^{2}ds.
$$
So we conclude that $\phi^{c}$ converges to $\phi^{\infty}$ in $\cC([0, T], \cC^{1}(B_{R}))$. The fact that $u^{\infty} \in \cC^{1}([0, T] \times \bbR^{d})$ implies that $\phi^{\infty} \in \cC^{1}([0, T] \times \bbR^{d})$. Also by Theorem~\ref{potential}, we know that when $u^{\infty}(0, x) = \nabla \phi_{0}(x)$, $u^{\infty}(s, x) = \nabla \phi^{\infty}(s, x)$ holds on $[0, T] \times  \bbR^{d}$, i.e. $\phi^{\infty}$ satisfies the equation~\eqref{eq112}, which ends our proof. 
\epf

\section{Further remarks}

To end this paper, let us give some further remarks.

\brmq In \cite{EKS},  Erbar-Kuwada-Sturm introduced a new definition of the curvature-dimension condition on metric-measure spaces, called the entropic curvature-dimension condition. By \cite{EKS},  we say that the  entropic curvature-dimension condition, denoted by  $CD_{\rm Ent}(K, m)$, holds if the Boltzmann entropy ${\rm Ent}$ satisfies 
\begin{eqnarray*}
{\rm Hess} {\rm Ent}-{1\over N}\nabla {\rm Ent}^{\otimes 2}\geq K,
\end{eqnarray*}
where $K\in \mathbb{R}$, $N\geq n$ are two constants.  As was pointed out in \cite{EKS}, when $M$ is a complete Riemannian manifold, the $CD_{\rm Ent}(K, m)$ is equivalent to the $CD(K, m)$-condition. 

Based on Erbar-Kuwada-Sturm's  new definition of  the entropic curvature-dimension condition $CD_{\rm Ent}(K, m)$, we can also prove the $W$-entropy inequalities for the geodesic flow, the gradient flow as well as the Langevin deformation of flows (of the Bolztmann-Shannon entropy) on the Wasserstein space over complete Riemannian manifolds. More precisely, we have the following

\bthm Let $M$ be a complete Riemannian manifold of dimension $n$. Suppose that Erbar-Kawada-Sturm's $CD_{\rm Ent}(K, N)$ condition holds, i.e., 
\begin{eqnarray*}
{\rm Hess} {\rm Ent}-{1\over N}\nabla {\rm Ent}^{\otimes 2}\geq K,
\end{eqnarray*}
where $K\in \mathbb{R}$, $N\geq n$ are two constants. 
Then \\
$(i)$ for geodesic flow $(\rho(t), \phi(t))$ on $TP_2(M, \mu)$, we have
\begin{eqnarray*}
{d^2\over dt^2}{\rm Ent}(\rho(t))+{2\over t}{d\over dt}{\rm Ent}(\rho(t))+{N\over t^2}
\geq{1\over N}\left|\langle \nabla {\rm Ent}(\rho(t)), \dot\rho(t))\rangle+{N\over t}\right|^2+K\|\dot\rho(t)\|^2.
\end{eqnarray*}
$(ii)$ for the gradient flow $\dot\rho(t)=-\nabla{\rm Ent}(\rho(t))$  on $P_2(M, \mu)$, we have
%\footnote{Note that $\langle \nabla {\rm Ent}(\rho(t)), \dot\rho(t))\rangle=\int_M |\nabla\log\rho|^2\rho d\mu$. Compare to $(\ref{wcdkm2})$ in  Remark \ref{wcdkm3}.}
\begin{eqnarray*}
{d^2\over dt^2}{\rm Ent}(\rho(t))+{1\over t}{d\over dt}{\rm Ent}(\rho(t))+{N\over 2}\left(K+{1\over t}\right)^2\geq {2\over N}\left|\langle \nabla {\rm Ent}(\rho(t)), \dot\rho(t))\rangle+{N\over 2}\left(K+{1\over t}\right)g\right|^2.
\end{eqnarray*}
\ethm

\bthm \label{Main-EKS} Let $c\in [0, \infty)$, and let $(\rho(t), \phi(t))$ be the Langevin deformation of flows on $TP_2(M, \mu)$. Suppose that Erbar-Kawada-Sturm's $CD_{\rm Ent}(K, N)$-condition holds for some constants $K\in \mathbb{R}$ and $N\in \mathbb{N}$ with $N\geq n$, i.e.,
\begin{eqnarray*}
{\rm Hess} {\rm Ent}-{1\over N}\nabla {\rm Ent}^{\otimes 2}\geq K.
\end{eqnarray*}
Let $\alpha(t)=(\log u)'$ be as in Section $6$ with $m=N$. 
Then
\begin{eqnarray*}
& &{d^2\over dt^2} {\rm Ent}(\rho(t))+\left(2\alpha(t)+{1\over c^2}\right){d\over dt}{\rm Ent}(\rho(t))+N\alpha^2(t)+{1\over c^2}\|\nabla {\rm Ent}(\rho(t))\|^2\\
& &\hskip2cm  \geq  {1\over N}\left|\langle \nabla{\rm Ent}(\rho(t)), \dot\rho(t)\rangle+N \alpha(t)\right|^2+K\|\dot \rho(t)\|^2.
\end{eqnarray*}
\ethm

The above results might bring some new insights to the study of geometric analysis on non smooth metric measure spaces. To save the length of this paper, we omit the proof of the above results. See Section $9$ of our previous preprint \cite{LL-flow17}. 

\ermq

\brmq As we have mentioned in Remark \ref{remrigidity} and Section $6$, it is possible to extend the $W$-entropy formula to the Langevin deformation of flows on the Wasserstein space over complete Riemannian manifolds with bounded geometry condition. This will implies the rigidity theorem as mentioned in Section $6$. To this end, we need only to verify some technical condition to exchange the differentiation and integration of the Boltzmann-Shannon entropy along the Langevin deformation of flows. To save the length of the paper, we omit the detail of the technical part. Moreover, as mentioned also in Section $6$, this paper we only we consider the Langevin deformation of flows for the Boltzmann-Shannon entropy  $\Ent(\rho) = \int_M \rho \log \rho d\mu$ on the Wasserstein space over Euclidean spaces and Riemannian manifolds. We can also consider the  Langevin deformation of flows for the R\'enyi entropy $V(\rho)={1\over m-1}\int_M \rho^{m} d\mu$ with $m\neq 1$. This will be studied in future. 

\ermq

\brmq As we have mentioned in Section $1$,  McCann and Topping \cite{MT} proved the contraction property of the $L^2$-Wasserstein distance between solutions of the backward heat equation on closed manifolds equipped with the Ricci flow. See also \cite{T1, T2}. In \cite{Lo2}, Lott proved the convexity  of the Boltzmann-Shannon entropy along the geodesics on the Wasserstein space over closed manifolds equipped with Ricci flow, which is  closely related to Perelman's results on the monotonicity of the $\mathcal{F}$ and $\mathcal{W}$-entropy functionals for Ricci flow. In \cite{LL13b}, the authors extended Lott's convexity results to the Wasserstein space on compact Riemannian manifolds equipped with Perelman's Ricci flow.  Inspired by these works, we can introduce the Langevin deformation of flows on the Wasserstein space over manifolds equipped with time dependent metrics and potentials, in particular, to manifolds  equipped with the Ricci flow and the $(K, m)$-Ricci flow  as introduced in our previous papers  \cite{LL15, LL17, LL18SCM}. This has been done and will be included in a forthcoming paper. 

\ermq

\brmq In \cite{KL21}, Kuwada and the second named author of this paper introduced the $W$-entropy for the heat flow on metric measure space and proved its monotonicity on metric measure space with RCD$(0, N)$ condition. It would be interesting to extend the main results of this paper to the $W$-entropy for the geodesic flow and Langevin deformation of flows on the Wasserstein space over metric measure spaces and prove their monotonicity on metric measure spaces with RCD$(0, N)$ condition. Moreover,  we can raise the question how to extend the monotonicity and rigidity theorems to the $W$-entropy along the Langevin deformation of flows on the Wasserstein space over super Ricci flows on metric measure spaces as introduced by Sturm \cite{St18} and Kopfer-Sturm \cite{KS18}.

\ermq

\brmq Finally, let us mention that one can also formally introduce the Langevin diffusions on the Wasserstein space  over Euclidean space and Riemannian manifolds as Bismut \cite{Bis05, Bis10} did on tangent bundle over Riemannian manifolds. See Section $4$. This gives an interpolation between the stochastic gradient flows (in particular, the Brownian motion) on the  Wasserstein space and the geodesic flow (or more general the Hamiltonian flow) on the tangent bundle over the  Wasserstein space. However, the question how to rigorously define the  infinite dimensional Brownian motion or  more general  infinite dimensional stochastic gradient flows on the  Wasserstein space remains open. 

\ermq

 \noindent\textbf{Acknowledgement}.  The authors would like to thank Professors S. Aida, D. Bakry, J.-M. Bismut,  G.-Q. Chen,  F.-M. Huang, K. Kuwada, K. Kuwae, K.-T. Sturm, C. Villani and  F.-Y. Wang   for helpful discussions and their interests on this work.

\begin{flushleft}
%\medskip\noindent

Songzi Li, School of Mathematics,  Renmin University of China,  Beijing, 100872, China
Email: sli@ruc.edu.cn
\medskip

Xiang-Dong Li, Academy of Mathematics and Systems Science, Chinese
Academy of Sciences, No. 55, Zhongguancun East Road, Beijing, 100190,  China\\
E-mail: xdli@amt.ac.cn

and

School of Mathematical Sciences, University of Chinese Academy of Sciences, Beijing, 100049, China
\end{flushleft}

\begin{thebibliography}{99}
\bibitem{AGS}
L.~Ambrosio, N.~Gigli, G.~Savar\'e, Gradient flows in metric spaces
  and in the space of probability measures, Lectures in Mathematics ETH
  Z\"urich. Birkh\"auser Verlag, Basel, 2005.
\bibitem{Au} T. Aubin, Nonlinear Analysis on Manifolds, Monge-Amp\`ere Equations, Springer-Verlag, 1980
\bibitem{BE} D. Bakry, M. Emery, Diffusion hypercontractives,
S\'em. Prob. XIX, Lect. Notes in Maths. 1123 (1985) 177-206.
\bibitem{Bes} A. L. Besse, Einstein manifolds, Ergebnisse der Mathematik (3) 10, Springer, Berlin, 1987. 
\bibitem{Bis05} J.-M. Bismut, The hypoelliptic Laplacian on the cotangent bundle. J. Amer. Math. Soc. 18 (2005), no. 2, 379-476.
\bibitem{Bis10} J.-M. Bismut, Hypoelliptic Laplacian and orbital integrals. Annals of Mathematics Studies, 177. Princeton University Press, Princeton, NJ, 2011.
%\bibitem{BK} B. Ben Moussa, G. T. Kossioris, On the system of Hamilton-Jacobi and transport equations %arising in geometrical optics, Comm. in Partial Diff. Equations, Vol 28, No. 5 and 6, 2003, 1085-1111.
%\bibitem{Bot}L. Boltzmann, Weitere Studien ber das Wsarmegleichgewicht unter Gasmoleksulen,
%Wiener Berichte, 66 (1872), 275-370.
%\bibitem{Br} Y. Brenier, Polar factorization and monotone rearrangement of vector-valued functions, Comm. Pure Appl. Math. 44, 4 (1991), 375-417.
\bibitem{BreGre}
Yann Brenier, Emmanuel Grenier, Sticky particles and scalar
  conservation laws, SIAM J. Numer. Anal. 35(1998), no.~6,
  2317--2328.
\bibitem{BB}  J.-D. Benamou, Y. Brenier, A computational fluid mechanics solution to the Monge-Kantorovich mass
transfer problem. Numer. Math. 84, 375-393 (2000)
\bibitem{Car} R. Carles, Semi-classical analysis for nonlinear Schr\"odingier equations, World Scientific Publishing, 2008.  
\bibitem{chen-liu}
Gui-Qiang Chen, Hai-Liang Liu, Formation of $\delta$-shocks and vacuum
  states in the vanishing pressure limit of solutions to the euler equations
  for isentropic fluids, SIAM J. Math. Anal. 34(2003), no.~4,
  925--938.
  \bibitem{CG}
Jean-Francoise Coulombel, Thierry Goudon, The strong relaxation limit
  of the multidimensional isothermal euler equations, Transactions of the
  American Mathematical Society 359(2007), no.~2, 637--648.
%\bibitem{CL} M. G. Crandall, P. L. Lions, Viscosity solutions of Hamilton-Jacobi equations, Trans. Amer. Math. Soc. 277 (1983), 1-42.
\bibitem{EKS} M. Erbar, K. Kuwada, K.-T. Sturm, On the equivalence of the entropic curvature- dimension condition and Bochner's inequality on metric measure spaces. Preprint. Available at: arXiv:1303.4382.
\bibitem{ERS}
Weinan E, Yu.G. Rykov, Ya.G. Sinai, Generalized variational
  principles, global weak solutions and behavior with random initial data for
  systems of conservation laws arising in adhesion particle dynamics, Commun.
  Math. Phys 177(1996), 349--380.
\bibitem{Ka} T. Kato, The Cauchy problem for quasi-linear symmetric hyperbolic systems. Arch Rational Mech Anal. 58 (1975), 181-205. 
\bibitem{KM1}
S.~Klainerman, A.~Majda, Singular limits of quasilinear hyperbolic
  systems with large parameters and the incompressible limit of compressible
  fluids, Comm. Pure Appl. Math. 34(1981), 481--524.
 \bibitem{KM2}
S.~Klainerman, A.~Majda, Compressible and incompressible fluids, Comm. Pure Appl. Math.
35(1982), 629--653.
\bibitem{KS18} E.~Kopfer, K.-T.~Sturm, Heat flow on time-dependent metric measure spaces and 
 super-Ricci flows. Comm. Pure Appl. Math. 71 (2018), no. 12, 
 2500-2608. 
\bibitem{KL21} K.~Kuwada, X.-D. ~Li, Monotonicity and rigidity of the $W$-entropy on
RCD$(0, N)$ spaces, Manuscripta Math. 164 (2021), 119-149.
 \bibitem{Lat-Tza13}
Corrado Lattanzio, Athanasios~E. Tzavaras, Relative entropy in
  diffusive relaxation, SIAM Journal on Mathematical Analysis 45(2013), no.~3, 1563--1584.
  \bibitem{Lat-Tza17}
Corrado Lattanzio, Athanasios~E. Tzavaras, From gas dynamics with large friction to gradient flows
  describing diffusion theories, Communications in Partial Differential
  Equations 42(2017), no.~2, 261--290.
\bibitem{Li07}X.-D. Li, On the $W$-entropy formula for the Witten Laplacian over weighted Riemannian manifolds, preprint 2006, and included in Th\`ese d'Habilitation \`a Diriger des Recherches, Universit\'e Paul Sabatier, 2007.
\bibitem{Li12}X.-D. Li, Perelman's entropy formula for the Witten Laplacian
on Riemannian manifolds via Bakry-Emery Ricci
curvature, Math. Ann.  353 (2012), 403-437.
\bibitem{Li16}X.-D. Li, Hamilton's Harnack inequality and the $W$-entropy formula on complete Riemannian manifolds, Stoch. Processes and  Appl. 126 (2016) 1264-1283
\bibitem{LL15} S. Li, X.-D. Li, $W$-entropy formula for the Witten Laplacian on manifolds with time dependent metrics and potentials, Pacific J. Math.  Vol. 278, No. 1, 2015, 173-199. 
\bibitem{LL13b} S. Li, X.-D. Li, On the convexity of Boltzmann type entropy on compact Riemannian manifolds equipped with Perelman's Ricci flow, preprint, 2013.
\bibitem{LL17} S. Li, X.-D. Li,  Hamilton differential Harnack inequality and W-entropy for Witten Laplacian on Riemannian manifolds, J Funct Anal, 2018, 274: 3263-3290, https://doi.org/10.1016/j.jfa.2017.09.017.
\bibitem{LL18SCM} S. Li, X.-D. Li, $W$-entropy formulas on super Ricci flows and
Langevin deformation on Wasserstein space over
Riemannian manifolds, Science China Mathematics, 2018 Vol. 61 No. 8: 1385-1406, 
https://doi.org/10.1007/s11425-017-9227-7.
\bibitem{LL-flow17} S. Li, X.-D. Li, $W$-entropy formulas and Langevin deformation of flows on Wasserstein space over Riemannian manifolds,  arXiv:1604.02596, 2016.
\bibitem{LL-JGA20} S. Li, X.-D. Li, $W$-Entropy, Super Perelman Ricci Flows, and $(K, m)$-Ricci
Solitons, J.  Geom. Anal. (2020) 30, 3149-3180, https://doi.org/10.1007/s12220-019-00193-4
%\bibitem{Lions} P. L. Lions, Generalized Solutions of Hamilton-Jacobi Equations, Research Notes in Mathematics, Boston, Pitman, pp. 69, 1982.
\bibitem{Lot} J. Lott, Some geometric properties of the Bakry-Emery Ricci tensor, Comment. Math. Helv. 78 30 (2003) 865-883.
\bibitem{Lo1} J. Lott, Some geometric calculations on Wasserstein space, Comm. Math. Phys. 277 (2008), no. 2, 423-437.
\bibitem{Lo2} J. Lott, Optimal transport and Perelman's reduced volume, Calc. Var. and Partial Differential Equations 36, 49-84 (2009).
\bibitem{LoV} J. Lott, C. Villani, Ricci curvature for metric-measure spaces via optimal transport, Annals of Math. 169, 903-991, 2009.
\bibitem{Maj} A. Majda, Compressible fluid flow and systems of conservation laws in several space variables. Applied Mathematical Sciences. New York, Springer-Verlag, 1984.
\bibitem{Mc} R. McCann, Polar factorization of maps on Riemannian manifolds, Geometric and Functional Analysis, 11, 3(2001), 589-608.
\bibitem{MT}  R. McCann, P. Topping, Ricci flow, entropy and optimal transpotation, American Journal of Math., 132 (2010) 711-730.
\bibitem{Nier} F. Nier, A semi-classical picture of quantum scattering, Ann. Sci. \'Ecole Norm. Sup. (4) 29 (1996), 2, 149-183.  
\bibitem{Ot}  F. Otto, The geometry of dissipative evolution equations: the porous medium equation, Commun. in Parial Differential Equations 26 (1 and 2), 101-174 (2001).
%\bibitem{Na} J. Nash, Continuity of solutions of parabolic and elliptic equations, Amer. J. Math. 80 (1958), 931-954.  
\bibitem{N1} L. Ni, The entropy formula for linear equation, J. Geom. Anal. 14 (1), 87-100, (2004).
\bibitem{OtV} F. Otto, C. Villani, Generalization of an inequality by Talagrand and links with the
logarithmic Sobolev inequality, J. Funct. Anal. 173 (2000) 361-400.
\bibitem{OtW} F. Otto, M. Westdickenberg, Eulerian calculus for the contraction in the Wasserstein distance. SIAM J. Math. Anal., 37(4):1227-1255 (electronic), 2005.
\bibitem{P1} G. Perelman, The entropy formula for the Ricci flow and its geometric applications, http://arXiv.org/abs/maths0211159.
\bibitem{sever}
Michael Sever, An existence theorem in the large for zero-pressure gas
  dynamics, Differential and Integral Equations 14(2001), no.~9,
  1077--1092.
\bibitem{STW} T. Sideris, B. Thomases, D.H. Wang, Long time behavior of solutions to the $3$D compressible Euler equations with damping, Comm. PDE. 28 (2003), no.3-4, 795-816. 
%\bibitem{Sou}P. E. Souganidis, Existence of viscosity solutions of Hamilton-Jacobi equations, J. Diff. Equations, 56 (1985), 345-390.
\bibitem{St1}K.-T. Sturm, Convex functionals of probability measures and nonlinear diffusions on
manifolds, J. Math. Pures Appl. (9) 84 (2005), no. 2, 149-168.
\bibitem{St2}K.-T. Sturm, On the geometry of metric measure spaces. I. Acta Math., 196(1): 65-131, 2006.
\bibitem{St3}K.-T. Sturm, On the geometry of metric measure spaces. II. Acta Math., 196(1), 133-177, 2006.
\bibitem{St18}K.-T. Sturm, Super-Ricci flows for metric measure spaces. J. Funct. Anal. 275 (2018), no. 12, 3504-3569. 
\bibitem{StR} K.-T. Sturm, M.-K. von Renesse,
Transport inequalities, gradient estimates, entropy, and Ricci
curvature Comm. Pure Appl. Math., 58(7), 923--940, 2005.
\bibitem{T1} P. Topping, Lectures on the Ricci Flow,  London Mathematical Society Lecture Note Series 325, Cambridge University Presse, 2006.
\bibitem{T2} P. Topping, $L$-optimal transportation for Ricci flow, J. Reine Angew. Math., 636 (2009) 93--122.
\bibitem{V1} C. Villani, Topics in Mass Transportation, Grad. Stud. Math., Amer. Math. Soc., Providence,
RI, 2003.
\bibitem{V2}C. Villani, Optimal Transport, Old and New, Springer, 2008.
\bibitem{WY} W. Wang, T. Yang, The pointwise estimates of solutions for Euler equations with damping in multi-dimensions, J. Diff. Equa. 173 (2001), 410-450.
\end{thebibliography}
\end{document}